\documentclass[11pt,a4paper]{amsart}
\usepackage[margin=2.5cm]{geometry}
\usepackage{contactgmt-macros}
\usepackage{xcolor}
\usepackage[colorlinks=true,linkcolor=blue!45!black,citecolor=green!35!black,urlcolor=blue!45!black]{hyperref}
\usepackage[expansion=false]{microtype}
\theoremstyle{plain}
\newtheorem*{claim}{Claim}

\begin{document}

\title[Height-localised monotonicity for Legendrian mean curvature flow]{Height-localised
monotonicity and tangent structures for Legendrian mean curvature flow in Sasakian
manifolds}

\author[Chang]{Shu-Cheng Chang}
\address{Shanghai Institute for Mathematics and Interdisciplinary Sciences\\
Shanghai 200433\\ China}
\email{scchang@math.ntu.edu.tw}

\author[Qiu]{Hongbing Qiu}
\address{School of Mathematics and Statistics\\ Wuhan University\\
Wuhan 430072\\ China}
\email{hbqiu@whu.edu.cn}

\author[Zhang]{Liuyang Zhang}
\address{School of Mathematical Science\\ Chongqing University of Technology,\newline
Chongqing 400054\\ China}
\email{zhangliuyang@cqut.edu.cn}

\thanks{The first author is partially supported by the Startup Foundation for Advanced
Talents of the Shanghai Institute for Mathematics and Interdisciplinary Sciences
(No.~\mbox{2302-SRFP-2024-0049}); the second author by the National Natural Science
Foundation of
China (No.~12471050) and the Hubei Provincial Natural Science Foundation of China
(No.~2024AFB746); the third author by the National Natural Science Foundation of China
(No.~12601100).}


\begin{abstract}
We introduce a Gaussian density for Legendrian mean curvature flow that is
localised in the Reeb direction by the angle-corrected height. At a smooth time,
the small-scale limit of this density, centred at a point in local contact
coordinates, counts the sheets through that point, with multiplicity. Its almost
monotonicity holds without any Type~I hypothesis, for closed Legendrian submanifolds with
exact mean curvature form in compact $\eta$-Einstein Sasakian manifolds, and more
generally under a hypothesis on the Legendrian angle. The construction applies even when
the Reeb foliation has no global K\"ahler quotient. As applications, near a finite-time
singularity every weak limit of the parabolic rescalings is a finite union of stationary
integral Lagrangian cones lifted to finitely many Reeb heights, and the set of Legendrian
angles carried by these cones does not depend on the sequence of rescalings. At a Type~I
singularity in the $\eta$-Einstein case, every tangent flow is, along every sequence of
scales, a finite union of horizontal Legendrian planes.
\end{abstract}

\maketitle

\tableofcontents

\section{Introduction}\label{sec:intro}

\subsection{The problem}
It was proven by Huisken \cite{hu90} that, along the mean curvature flow of hypersurfaces
in Euclidean space, the integral of the backward heat kernel over the evolving
hypersurface is non-increasing in time. This monotonicity formula is an important tool in
the study of singularities by rescaling, and this paper is concerned with a
height-localised counterpart of it for the Legendrian mean curvature flow.

Let $(M^{2n+1},\eta,\xi,\Phi,g)$ be a compact $\eta$-Einstein Sasakian manifold with
constant $K$, and let $F\colon\Sigma\times[0,T)\to M$ be a Legendrian mean curvature flow,
in the sense of Smoczyk \cite{s1}, of a closed Legendrian submanifold with exact mean
curvature form, written in the gauge of \cite[(1.1),
(2.2)]{ChangWuZhang2026ContactBlowup}. Then there is a function $\theta$ on
$\Sigma\times[0,T)$, the \emph{Legendrian angle}, with
\begin{equation}\label{eq:LMCF}
H=\Phi\nabla\theta,\qquad \partial_tF=H-2\theta\xi,\qquad
(\partial_t-\Delta)\theta=(K+2)\theta .
\end{equation}
This flow differs from the one of \cite[Definition~3.5]{s1}, which is written with the
Legendre pseudo-angle, by a rotation along the Reeb flow. In a general compact Sasakian
manifold the system \eqref{eq:LMCF} is taken as a hypothesis, (F) of
Section~\ref{subsec:hypotheses}. The flow moves horizontally by its mean curvature $H$ and
vertically, along the Reeb field $\xi$, with the speed $-2\theta$.

Recently, the first and third authors, together with Wu, studied the singularities of this
flow in \cite{ChangWuZhang2026ContactBlowup} by contact parabolic rescalings, the
$\cD$-homothetic deformations of the Sasakian structure, whose limit space is the
Heisenberg group $\Hn$ with its standard Sasakian structure. In the rescaled variables the
Reeb direction is stretched quadratically. The limits of the rescaled flows, the
\emph{contact tangent flows}, are described in \cite{ChangWuZhang2026ContactBlowup} by a
transverse Gaussian density, built from the horizontal potential $f=\frac14|X^{\cH}|^2$,
the lift of the Euclidean potential of the leaf space. However, this density is constant
along the Reeb orbits. In fact, it counts all the sheets of the flow which approach a Reeb
orbit, at whatever height \cite[Remarks~4.4 and~4.13]{ChangWuZhang2026ContactBlowup}, so
it cannot separate sheets that lie over the same transverse point at different Reeb
heights. The simplest example is two horizontal Legendrian planes $\{y=0,\ z=0\}$ and
$\{y=0,\ z=c\}$ in $\Hn$, $c\ne0$: the transverse density is $2$ at every point of either
plane, although the planes are disjoint. This is why, in the blow-up analysis of
\cite{ChangWuZhang2026ContactBlowup}, the vertical position of the limits has to be
determined by a separate argument. It is also why the local regularity theorem of White
\cite{wh}, applied to the projected flow, only gives a statement in terms of the
transverse density.

In \cite[Section~1 and Remark~4.13]{ChangWuZhang2026ContactBlowup} a local regularity
theorem in terms of a density localised in the Reeb direction is identified as the missing
tool for separating the sheets over one Reeb orbit, for instance in order to exclude
Type~I singularities at points where several sheets meet. The natural candidates for such
a density do not work \cite[Remark~4.4]{ChangWuZhang2026ContactBlowup}. The function
$\frac14|X|^2=f+\frac14z^2$, where $X$ is the full generator of the dilations, is not
constant along the Reeb orbits. Hence the vertical velocity $-2\theta\xi$ enters its
evolution, and the resulting expression is not a square. The Riemannian Gaussian kernel,
on the other hand, is adapted to rescalings whose limits are Lagrangian.

To address this difficulty, we keep track of the Reeb height. The tool is a Gaussian
kernel with a weight in the height direction. In the model $\Hn$, with coordinates
$(w,z)\in\C^n\times\R$, horizontal potential $f=|w|^2/16$ and $\tau=T-t$, it is
\begin{equation}\label{eq:kernel-intro}
\rho^\chi=(4\pi\tau)^{-n/2}e^{-f/\tau}\,\chi(v),\qquad v=\frac{z-4\tau\theta}{\tau},\qquad
\chi(v)=e^{-av^2},\quad 0<a<\tfrac1{32}.
\end{equation}
The Legendrian angle enters through the \emph{caloric height} $u=z-4\tau\theta$, the
height corrected for the vertical drift of the flow (Section~\ref{subsec:relation-M}
explains its relation to an evolution formula of Neves). The weight $\chi(v)$ sees only the
sheets whose height, at the time $t$, is within $O(\tau)$ of the height $4\tau\theta$
predicted by their Legendrian angle. In particular, for the two planes above, the weighted
density at a point of the first plane at the scale $\tau$ is $1+e^{-ac^2/\tau^2}$, which
tends to $1$ as $\tau\to0$ (Example~\ref{ex:two-sheets}). More generally, for $t_*<T$ the
localised density of this kernel centred at $(q,t_*)$ tends, as the scale goes to $0$, to
the number of sheets of $\Sigma_{t_*}$ through $q$ (Theorem~\ref{thm:A}(ii)). The
transverse density, by contrast, counts all the sheets which approach the Reeb orbit of
$q$.

On the other hand, in a curved Sasakian manifold the kernel is written in contact normal
coordinates, which are the coordinates of \cite[Lemma~4.1]{ChangWuZhang2026ContactBlowup}
made uniform in the centre, with $a_K(\tau)=4(e^{(K+2)\tau}-1)/(K+2)$ (and
$a_K(\tau)=4\tau$ if $K+2=0$) in place of $4\tau$. The construction uses the
transverse K\"ahler geometry of the local contact charts and retains the Reeb height of
each sheet. It applies even when the Reeb foliation has no global K\"ahler quotient.
Our main technical result,
Theorem~\ref{thm:A}, states that the kernel is almost monotone, with an error that is
integrable up to the singular time, under the standing hypotheses alone and without a
Type~I hypothesis. Thus the kernel provides, at every centre, a density localised in the
Reeb direction, as called for above. The proof of Theorem~\ref{thm:A} rests on an idea of
Neves \cite[Lemma~5.4]{nev07}, who used it to obtain stationary limits for the Lagrangian
mean curvature flow with bounded Lagrangian angle.

This paper opens a series on the local structure of singularities in Legendrian
mean curvature flow. Our starting point is a Gaussian density localised by the
angle-corrected Reeb height. We use its almost monotonicity
(Theorem~\ref{thm:A}) to describe weak contact tangent flows near finite-time
singularities and the Legendrian angles they carry.

\subsection{Main results}\label{subsec:main}
In this introduction $(M,\eta,\xi,\Phi,g)$ denotes a compact $\eta$-Einstein Sasakian
manifold with constant $K$, $\Sigma$ is a closed $n$-manifold, and
$F\colon\Sigma\times[0,T)\to M$ with $T<\infty$ is a smooth Legendrian mean curvature flow
with exact mean curvature form, written as in \eqref{eq:LMCF} with the normalisation of
\cite[(2.2)]{ChangWuZhang2026ContactBlowup}. In the body of the paper
Theorems~\ref{thm:A}, \ref{thm:B} and~\ref{thm:D} are proved under the hypotheses (S) and
(F) of Section~\ref{sec:setting}, together with the hypotheses written in their
statements. These only require that $M$ be a compact Sasakian manifold and that the flow
admit a Legendrian angle satisfying \eqref{eq:LMCF} for some constant $K$; in the
$\eta$-Einstein case this is automatic. Theorem~\ref{thm:C} is proved under (S), (F), the
Type~I bound and the assumption that $M$ is $\eta$-Einstein with the constant $K$ of (F);
for general compact Sasakian manifolds Section~\ref{subsec:typeI} gives a version in which
the compactness hypothesis (C) is assumed. By the maximum principle $|\theta|\le C_\theta$
on $\Sigma\times[0,T)$, and we write $\mu_t$ for the induced area measure. The notation is
fixed in Sections~\ref{sec:setting} and~\ref{sec:mono}; in particular $\psi_q$ is a
cut-off supported in the contact normal chart of $q$, $\rho^\chi_{q,t_*}$ is the kernel
\eqref{eq:kernel-intro} centred at $(q,t_*)$ and written in that chart,
$\cS_q=H+(\overline\nabla f_q)^\perp/\tau$ is Huisken's defect of self-similarity, and
\[
D^\chi=\frac{|v\chi'(v)|}{\tau}+(\chi+16\chi'')(v)\,|\cS_q|^2\ \ge\ \frac{|v\chi'(v)|}\tau+c_\chi\chi(v)|\cS_q|^2
\]
is the defect of the weighted kernel. Our first main result is a height-localised version
of Huisken's monotonicity formula.

\begin{mainthm}[Height-localised almost monotonicity]\label{thm:A}
There are $\tau_0>0$ and $C$, depending only on $n$, $a$, $K$, $C_\theta$, $T$, the
initial immersion and $(M,\eta,g)$, with the following properties, for every $q\in M$ and
every $t_*\in(0,T]$.
\begin{enumerate}[(i)]
\item On $[t_*-\min\{\tau_0,t_*\},t_*)$,
\[
\frac{d}{dt}\int\psi_q\,\rho^\chi_{q,t_*}\,d\mu_t\le-\frac34\int\psi_q\,D^\chi\rho_{q,t_*}\,d\mu_t+C\Big(1+\tau^2\int\psi_q|H|^2\rho_{q,t_*}\,d\mu_t\Big),
\]
and the last term is integrable in time. Consequently $\lim_{t\to
t_*}\int\psi_q\rho^\chi_{q,t_*}\,d\mu_t$ exists and
$\int\!\!\int\psi_qD^\chi\rho_{q,t_*}\,d\mu_t\,dt<\infty$.
\item If $t_*<T$, the localised contact density
$\Theta^\chi_{q,t_*}(\sigma)=\int\psi_q\rho^\chi_{q,t_*}\,d\mu_{t_*-\sigma}$ satisfies
\[
\lim_{\sigma\to0}\Theta^\chi_{q,t_*}(\sigma)=\#F_{t_*}^{-1}(q),\qquad
\Theta^\chi_{q,t_*}(\sigma)\ge\#F_{t_*}^{-1}(q)-C\sigma\quad\text{for
}0<\sigma\le\min\{\tau_0,t_*\} .
\]
\item There are $c_{\rm cl},R_*,Z_*>0$ such that the following holds. If $t_*<T$ and, for
some $0<\sigma\le\min\{\tau_0,t_*,1/(4C)\}$, the part of $\Sigma_{t_*-\sigma}$ in the box
$\{f_q<R_*^2\sigma,\ |z_q|<Z_*\sigma\}$ of the chart of $q$ has area at most $c_{\rm
cl}\sigma^{n/2}$, then $q\notin F_{t_*}(\Sigma)$.
\end{enumerate}
\end{mainthm}

Theorem~\ref{thm:A} is Propositions~\ref{prop:theta2} and~\ref{prop:mono} together with
Lemmas~\ref{lem:density} and~\ref{lem:clearing}. In particular, by (ii), the sheets of
$\Sigma_{t_*}$ which cross the Reeb orbit of $q$ at other heights do not contribute to the
limit. On the other hand, the height-blind version of (i), with $\chi\equiv1$, is a chart
version, of the same strength, of the localised monotonicity of
\cite[Proposition~4.16]{ChangWuZhang2026ContactBlowup} (Remark~\ref{rem:mono-M}).

Our second result describes the weak limits of the contact parabolic rescalings
$F^i_s=\delta_{\lambda_i}\circ\Upsilon_p\circ F_{T+\lambda_i^{-2}s}$ at a point $(p,T)$,
$\lambda_i\to\infty$. Here $\Upsilon_p$ is the contact normal chart of $p$,
$\delta_\lambda(w,z)=(\lambda w,\lambda^2z)$, and $F^i_s$ is defined on the part of
$\Sigma$ that $F_{T+\lambda_i^{-2}s}$ maps into the domain of the chart. Let $\mu$ denote
the area measure of $F^i_s$ for the rescaled metric, which is $\lambda_i^n$ times the area
measure of $F_{T+\lambda_i^{-2}s}$ on that part of $\Sigma$. We follow the position and
the Legendrian angle together: $\fm^i_s=(F^i_s,\theta^i_s)_\#\mu$ is a measure on
$\Hn\times\R$. A \emph{weak contact tangent flow} is a limit $(\mu_s,\fm_s)$ of
$(F^i_s)_\#\mu$ and $\fm^i_s$ for every $s$ in a fixed interval $[-S_0,0)$ along one
subsequence (Definition~\ref{def:weak}); a diagonal argument gives limits, along one
subsequence, for every $s<0$ (Section~\ref{sec:weak}). Let $g_0=\frac14|dw|^2$ be the
transverse metric of the model.

\begin{mainthm}[Weak contact tangent flows]\label{thm:B}
Let $p\in M$ and $\lambda_i\to\infty$. A subsequence of the rescalings at $(p,T)$
converges, for every $s<0$, to a weak contact tangent flow $(\mu_s,\fm_s)_{s<0}$ for all
negative times (Definition~\ref{def:weak}), and there are finitely many numbers
$\bar\theta_1<\dots<\bar\theta_N$ in $[-C_\theta,C_\theta]$, with $N$ bounded by a
constant depending only on the data of Theorem~\ref{thm:A}, and nonzero stationary
integral Lagrangian cones $\cC_1,\dots,\cC_N$ in $(\C^n,g_0)$ with vertex $0$ such that,
for \emph{every} $s<0$,
\[
\fm_s=\sum_{j=1}^N\big(\|\cC_j\|\times\delta_{4|s|\bar\theta_j}\big)\otimes\delta_{\bar\theta_j},\qquad
\mu_s=\sum_{j=1}^N\|\cC_j\|\times\delta_{4|s|\bar\theta_j}.
\]
Here $N\ge0$; for $N=0$ every sum in this theorem is empty, and the limit is empty:
$\mu_s=0$ and $\fm_s=0$ for every $s<0$. Thus $\mu_s$ is the weight of the stationary
Legendrian integral varifold $\sum_j\cC_j\times\{4|s|\bar\theta_j\}$ in $\Hn$, the
Legendrian angle on the $j$-th lifted cone equals $\bar\theta_j$, and each lifted cone
translates vertically with speed $-4\bar\theta_j$, reaching the origin at $s=0$. The
measure
\[
\nu_{p,T}:=\sum_jm_j\,\delta_{\bar\theta_j},\qquad m_j:=\Theta(\cC_j,0)\ \text{the vertex
density of }\cC_j,
\]
does not depend on the sequence $(\lambda_i)$. If $\Sigma$ is orientable, each $\cC_j$ is
a finite sum of special Lagrangian cones with respect to the geometric phases of a fixed
holomorphic volume form on $\C^n$.
\end{mainthm}

Theorem~\ref{thm:B} is Theorem~\ref{thm:weak}, extended to all $s<0$ by the diagonal
argument given after its proof, together with Proposition~\ref{prop:transport} and
Corollary~\ref{cor:angle-measure}. Two lifted cones with different Legendrian angles lie
at the height distance $4|s||\bar\theta_j-\bar\theta_k|$; the vertical segment between
them over their common vertex $w=0$, and over any other common projected point, is a Reeb
chord of $\eta$-action $2|s||\bar\theta_j-\bar\theta_k|$, which shrinks to zero at the
singular time (Corollary~\ref{cor:heights}). When a tangent flow consists of two planes
with different Legendrian angles, Proposition~\ref{prop:two-planes} confines the flow near
the singular point to two height bands; we do not show that the flow forms a Legendrian
double point there.

The statement that the angles $\bar\theta_j$ do not depend on the
sequence is the Legendrian form of the last assertion of Neves' theorem
\cite[Theorem~A]{nev07}, and we prove it by following his argument; in addition, the
measure $\nu_{p,T}$ records the densities carried by the angles. However, the uniqueness
of $\nu_{p,T}$ does not imply the uniqueness of the cones $\cC_j$ or of the tangent flow.

Our third result concerns Type~I singularities, $|B|^2\le C_I/(T-t)$, the setting of
\cite[Theorem~1.1]{ChangWuZhang2026ContactBlowup}. In this setting it was shown in
\cite{ChangWuZhang2026ContactBlowup} that the contact tangent flows are, up to vertical
translations of their components, self-similar flows of Legendrian self-shrinkers in
$\Hn$, and that along a suitable sequence of scales the translations vanish. Through the
rigidity theorem \cite[Theorem~A]{cqz} for Legendrian self-shrinkers, the components were
then shown to be planes. Under the Type~I bound the maps $F_t$ converge uniformly, as
$t\to T$, to a continuous map $F_T\colon\Sigma\to M$ (Lemma~\ref{lem:confinement});
$F_T(\Sigma)$ is the set of points reached by the flow at the singular time.

\begin{mainthm}[Type~I tangent flows]\label{thm:C}
Assume $|B|^2\le C_I/(T-t)$ on $\Sigma\times[0,T)$, and let $p\in M$ and
$\lambda_i\to\infty$. Then a subsequence of the rescalings at $(p,T)$ converges smoothly
on compact subsets of $\Hn\times(-\infty,0)$ to a flow $F^\infty$ with the following
properties. For every $s<0$, each component of $F^\infty_s$ is a horizontal Legendrian
plane $P_j\times\{4|s|\bar\theta_j\}$ over a Lagrangian subspace $P_j\subset\C^n$, where
$\bar\theta_j$ is the limit of the Legendrian angle on that component; there are at most
$2^nC_1/\omega_n$ components, counted with multiplicity, with $C_1$ the area constant of
\cite[Lemma~4.17]{ChangWuZhang2026ContactBlowup} and $\omega_n$ the volume of the unit
ball of $\R^n$; and $F^\infty$ is nonempty if and only if $p\in F_T(\Sigma)$.
\end{mainthm}

Theorem~\ref{thm:C} is Theorem~\ref{thm:planes-eta}. It sharpens \cite[Theorem~1.1 and
Corollary~1.2]{ChangWuZhang2026ContactBlowup} in three respects. First, the vertical
translations vanish for every component and along every sequence of scales. Second, the
identification of the components as planes uses neither the rigidity theorem
\cite[Theorem~A]{cqz} nor the orientability of the components, which in
\cite{ChangWuZhang2026ContactBlowup} is established through the orientation double cover
before that theorem is applied. Third, the centre $p$ need not be a singular point. The
compactness of the rescalings is taken from \cite[Lemmas~4.8, 4.9
and~4.17]{ChangWuZhang2026ContactBlowup}, and the centring is a consequence of
Theorem~\ref{thm:A}. Furthermore, Section~\ref{subsec:typeI} contains a version of this
result for general compact Sasakian manifolds, in which the smooth compactness of the
rescalings is a hypothesis.

Finally, we turn to Legendrian curves. For $n=1$ the closed manifold $\Sigma$ is
orientable, and the special Lagrangian decomposition of Theorem~\ref{thm:B} determines the
cones.

\begin{mainthm}[Legendrian curves]\label{thm:D}
Let $n=1$.
\begin{enumerate}[(a)]
\item In Theorem~\ref{thm:B} each cone $\cC_j$ is a finite union of lines through the
origin of $\C$, with integer multiplicities. Thus every weak contact tangent flow at
$(p,T)$ for all negative times is, for every $s<0$, a finite union of horizontal
Legendrian lines through the Reeb axis, with multiplicities, at the heights
$4|s|\bar\theta_j$.
\item Let $R>0$, let $s_0<0$ be a good time in the sense of
Proposition~\ref{prop:compactness}(c), and let $\Sigma_i\subset\Sigma$ be open sets such
that $F^i_{s_0}(\partial\Sigma_i)$ does not meet the ball $B_{3R}\subset\R^3$ and such
that $\Sigma_i\cap(F^i_{s_0})^{-1}(B_{2R})$ has exactly one connected component meeting
$(F^i_{s_0})^{-1}(B_R)$. Then, after passing to a subsequence, there is $\bar z$ such that
$z\circ F^i_{s_0}\to\bar z$ uniformly on $\Sigma_i\cap(F^i_{s_0})^{-1}(B_R)$ and
$\theta^i_{s_0}\to\bar z/(4|s_0|)$ in measure there: a local component of the rescaled
flow has a single height and a single Legendrian angle in the limit.
\end{enumerate}
\end{mainthm}

Theorem~\ref{thm:D} is Theorem~\ref{thm:curves}. Part (b) is the Legendrian counterpart,
for curves, of Neves' theorem \cite[Theorem~B]{nev07}, which gives a single Lagrangian
angle on each local component of the rescaled flow for almost calibrated and rational
Lagrangians. In our setting the hypotheses are only (S) and (F), and the statement
concerns the height and the Legendrian angle. For $n\ge2$ the corresponding statement is
not proved here.

\subsection{Relation to the Lagrangian theory}\label{subsec:lagrangian}
Our strategy follows closely the analysis by Neves \cite{nev07,nev11} of the zero-Maslov
class Lagrangian mean curvature flow in $\C^n$, which builds on Huisken's monotonicity
formula \cite{hu90} and on the work of Wang \cite{w4} and Chen--Li \cite{cl04} on almost
calibrated Lagrangians. In this analysis Huisken's monotonicity formula gives weak limits
that are self-similar. The same formula, weighted by the square of the bounded Lagrangian
angle, gives $\int\!\!\int|H|^2\rho<\infty$, so that the limits are stationary, hence
stationary cones. Finally, the evolution of the Lagrangian angle gives a finite set of
angles in the limit. However, the Legendrian setting differs in three respects, and each
of them causes a difficulty that does not arise in $\C^n$.

\emph{Curved background.} The flow lives in a compact Sasakian manifold, and the rescaled
structures converge to the Heisenberg group only in $C^\infty_{\rm loc}$. Hence the
identities of the model hold only up to error terms, and these errors have to be absorbed
up to the singular time without a Type~I bound. The orders of these errors are attained on
the round three-sphere (Remark~\ref{rem:S3}).

\emph{The height as a primitive.} On a Legendrian submanifold of $\Hn$ the function $z/2$
is a primitive of the Liouville form of $\C^n$, so that the Reeb height plays the role of
the primitive of the Liouville form which Neves uses for exact Lagrangians. In the outline
of the proof of \cite[Theorem~B]{nev07} given in \cite[Section~3]{nev07}, under the
hypotheses of that theorem, this primitive converges, on the $j$-th cone of a tangent flow
of an exact Lagrangian, to $\bar\beta_j-2s\bar\theta_j$, where the constants $\bar\beta_j$
are not controlled and may differ from cone to cone. In $\C^n$ these constants only change
the primitive. In $\Hn$, however, a change of the constant is a vertical translation, so
that components with different constants are different submanifolds over the same
transverse point, and the transverse density, which is constant along the Reeb orbits,
cannot separate them.

\emph{Legendrian angle and geometric phase.} The function $\theta$ of \eqref{eq:LMCF} is a
primitive of the mean curvature form, normalised by its evolution equation. For an
oriented Legendrian immersion in a contact chart there is also the geometric phase
$\vartheta$ of the projection with respect to a holomorphic volume form, the Lagrangian
angle of the projected immersion. Only the differential $d(\theta-\vartheta)$ is
controlled (Lemma~\ref{lem:phase}); the locally constant part is not.
Example~\ref{ex:potential-not-phase} shows that $\theta$ is in general not a lift of
$\vartheta$, and that two components over the same plane can have Legendrian angles whose
difference is not a multiple of $\pi$, so that no single holomorphic volume form makes
both of them lifts of the geometric phases. In $\C^n$, by contrast, the Lagrangian angles
of components with the same support differ by multiples of $\pi$ (cf.\
\cite[Remark~3.1]{nev07}).
Consequently the estimates for $\theta$ cannot replace conditions on $\vartheta$ such as
almost calibration, which Neves uses for the single-angle theorem \cite[Theorem~B]{nev07}.
The orientation hypothesis in Theorem~\ref{thm:B} gives the decomposition of the cones by
their geometric phases; it does not identify the two angles.

\subsection{Relation to the contact blow-up analysis}\label{subsec:relation-M}
The present paper builds on \cite{ChangWuZhang2026ContactBlowup}. From there we take the
contact normal coordinates at a point \cite[Lemma~4.1]{ChangWuZhang2026ContactBlowup}, the
area ratio bounds \cite[Lemma~4.15 and Proposition~4.16]{ChangWuZhang2026ContactBlowup},
from which the projected area ratio bound (D) of Section~\ref{subsec:D} follows, and, in
the $\eta$-Einstein Type~I case, the smooth compactness of the rescalings
\cite[Lemmas~4.8, 4.9 and~4.17]{ChangWuZhang2026ContactBlowup}. The basic model identities
of Section~\ref{subsec:model}, namely Lemma~\ref{lem:model-exact}(a)--(c) and the first
two identities of Lemma~\ref{lem:model-exact}(d), are those of \cite[Section~2 and
Proposition~4.3]{ChangWuZhang2026ContactBlowup}. On each time slice in the Heisenberg
model, half the Reeb height is a Liouville primitive for the projected Lagrangian
immersion. Under this identification, the combination used by Neves
\cite[Section~6]{nev07} is half the angle-corrected Reeb height. We use this height to
localise the transverse Gaussian in the Reeb direction. However, the transverse density
of \cite[Section~4.2]{ChangWuZhang2026ContactBlowup} cannot give the centring: being
constant along the Reeb orbits, it does not distinguish a self-similar flow from its
vertical translates. Relative to \cite{ChangWuZhang2026ContactBlowup}, the present paper
adds the weight in the height direction,
with its consequences: the almost monotonicity without Type~I (Theorem~\ref{thm:A}), the
weak contact tangent structure (Theorem~\ref{thm:B}), and the exact centring of the Type~I
tangent flows along every sequence (Theorem~\ref{thm:C}). On the other hand, the results
of \cite[Corollary~1.3 and Lemma~4.12]{ChangWuZhang2026ContactBlowup} concern the
transverse density of a regular Sasakian manifold and White's local regularity theorem for
the projected flow; we do not prove the corresponding statements for the height-localised
density here.

Monotonicity-type results that see the Reeb direction have been obtained for
\emph{stationary} Legendrian surfaces in the Heisenberg group. Rivi\`ere \cite{riv24}
proved an almost monotonicity formula for smooth Hamiltonian stationary (H-minimal)
Legendrian immersions, and Pigati and Rivi\`ere \cite{pr24} a monotonicity formula for a
modified density of Hamiltonian stationary Legendrian varifolds. Both are localised in
Kor\'anyi balls, and the second extends, up to error terms, to closed Sasakian manifolds
\cite[Remark~III.7]{pr24}. For Legendrian area-minimising currents in $\Hn$, Orriols
\cite{or26} develops the partial regularity theory without a monotonicity formula, with a
lower density bound obtained from an isoperimetric inequality in Kor\'anyi balls. The same
paper shows that the monotonicity of the area density fails for Hamiltonian stationary
Legendrian submanifolds of $\Hn$ in dimensions $n>2$.

All these results concern the stationary case. Along the flow \eqref{eq:LMCF} the sheets
also move along the Reeb orbits, and a monotone quantity has to follow them in the
vertical direction. The kernel of this paper is parabolic and Gaussian in both directions,
and it is adapted to the flow through the drift term $-4\tau\theta$. For a minimal
Legendrian submanifold, on which $\theta$ is locally constant and the flow is a rotation
along the Reeb field, the weighted kernel reduces to a static one and gives an almost
monotone density whose limit counts the sheets through the point
(Remark~\ref{rem:static}).

\subsection{Organisation}
The paper is organised as follows. In Section~\ref{sec:setting} we fix the conventions,
the hypotheses, the contact normal coordinates and the rescalings, and derive the area
ratio bound (D). Section~\ref{sec:mono} is devoted to the kernel: we introduce it in the
model, prove the almost monotonicity in a Sasakian manifold and the properties of the
localised contact density, and record the scale invariance of the quantities involved. In
Section~\ref{sec:weak} we prove Theorem~\ref{thm:B}. Section~\ref{sec:applications}
contains the Type~I analysis (Theorem~\ref{thm:C}), the structure of the tangent flows at
a point where one weak tangent flow consists of two planes with different Legendrian
angles, and the analysis of Legendrian curves (Theorem~\ref{thm:D}). Finally,
Appendix~\ref{app:A} proves a finite-phase compactness theorem for Lagrangian immersions
in metrics close to a flat one, which is used for the special Lagrangian decomposition in
Theorem~\ref{thm:B}, and Appendix~\ref{app:B} verifies the curvature errors on the round
three-sphere.

\section{Conventions, hypotheses and local geometry}\label{sec:setting}

In this section we fix the conventions and the hypotheses, and we construct the local
coordinates in which the kernel of Section~\ref{sec:mono} is written. Most of the material
comes from \cite{ChangWuZhang2026ContactBlowup}, and in the model the identities of
Lemma~\ref{lem:model} are those proved there (Remark~\ref{rem:charts-M}). Two points are
added here. First, the kernel of Section~\ref{sec:mono} depends on the height, so the
curvature errors in the height direction have to vanish at the centre to the orders of
Lemma~\ref{lem:model}(c), (d) and of the term $e_z$ in (e). This is achieved by the radial
gauge of Lemma~\ref{lem:charts}(ii), and the orders of Lemma~\ref{lem:model}(b)--(d) are
attained on the round three-sphere (Remark~\ref{rem:S3}). Second, the coordinates are
constructed uniformly in the centre, on boxes of a fixed size, which is needed for centres
that move with the rescalings (Lemma~\ref{lem:moving}(c)). We also derive the projected
area ratio bound (D) from \cite[Proposition~4.16]{ChangWuZhang2026ContactBlowup}
(Lemma~\ref{lem:D}), relate the Legendrian angle to the geometric phase
(Lemma~\ref{lem:phase}), and define the contact parabolic rescalings.

\subsection{Sasakian manifolds and the Heisenberg model}\label{subsec:sasakian}
We use the conventions of \cite[\S2.1]{ChangWuZhang2026ContactBlowup}. Let
$(M^{2n+1},\eta,\xi,\Phi,g)$ be a Sasakian manifold with contact distribution
$\cD=\ker\eta$, so that
\[
\Phi^2=-I+\eta\otimes\xi,\qquad d\eta(U,W)=2g(U,\Phi W),\qquad \overline\nabla_U\xi=-\Phi
U,\qquad g=g^T+\eta\otimes\eta ,
\]
where $\overline\nabla$ is the Levi-Civita connection of $g$. On a local leaf space $Z$ of
the Reeb foliation the transverse metric $g_Z$, $\pi^*g_Z=g^T$, is K\"ahler for the
complex structure $J$ induced by $\Phi$, with K\"ahler form $\omega_Z=2g_Z(\cdot,J\cdot)$,
and $d\eta=\pi^*\omega_Z$ \cite[\S2.3]{ChangWuZhang2026ContactBlowup}. The structure is
$\eta$-Einstein if $\mathrm{Ric}=Kg+(2n-K)\eta\otimes\eta$ for a constant $K$; then
$\mathrm{Ric}^T=(K+2)g^T$. For $\kappa>0$ the $\cD$-homothetic structure is
$(\eta_\kappa,\xi_\kappa,\Phi,g_\kappa)=(\kappa\eta,\kappa^{-1}\xi,\Phi,\kappa
g+(\kappa^2-\kappa)\eta\otimes\eta)$ \cite[(2.4)]{ChangWuZhang2026ContactBlowup}.

The model is the Heisenberg group $\R^{2n+1}$ with the structure
$(\overline\eta,\overline\xi,\overline\Phi,\overline g)$ of
\cite[\S2.4]{ChangWuZhang2026ContactBlowup}. In coordinates $(w,z)\in\C^n\times\R$,
$w=(x,y)$,
\begin{gather*}
\overline\eta=\tfrac12dz+\overline\beta,\qquad
\overline\beta=-\tfrac14\sum_j(y_j\,dx_j-x_j\,dy_j),\qquad \overline\xi=2\partial_z,\\
\overline g=g_0+\overline\eta\otimes\overline\eta,\qquad
g_0:=\tfrac14\sum_j(dx_j^2+dy_j^2) .
\end{gather*}
Thus $\overline g$ is not the Euclidean metric $|dw|^2+dz^2$: the horizontal factor
$\frac14$ corresponds to the coefficient $\frac1{16}$ in the Gaussian below, and
$\overline\xi=2\partial_z$ means that the $\overline\eta$-length $A_\eta$ of a Reeb
segment is half its $z$-length. The left translations are
\[
\ell_{(a,b,c)}(x,y,z)=\big(x+a,\ y+b,\ z+c+\tfrac12(b\cdot x-a\cdot y)\big);
\]
they preserve $\overline\eta,\overline\xi,\overline\Phi,\overline g$. The vertical
translations are $\tau_b:=\ell_{(0,0,b)}$, $\tau_b(w,z)=(w,z+b)$, and the contact
dilations are $\delta_\lambda(w,z)=(\lambda w,\lambda^2z)$, with
$\delta_\lambda\tau_b=\tau_{\lambda^2b}\delta_\lambda$
\cite[\S2.4]{ChangWuZhang2026ContactBlowup}. On horizontal lifts $\overline\Phi$ induces
the complex structure $J_0$ of $\C^n$ with $J_0\partial_{x_j}=-\partial_{y_j}$ and
$J_0\partial_{y_j}=\partial_{x_j}$, since $d\overline\eta=\frac12\sum_jdx_j\wedge
dy_j=:\omega_0$ and $d\overline\eta(U,W)=2\overline g(U,\overline\Phi W)$. The functions
$\zeta_j:=\frac12(x_j-i\,y_j)$ are $J_0$-holomorphic, and we put
$\Omega_0:=d\zeta_1\wedge\dots\wedge d\zeta_n$. For a Hermitian metric $h$ on $\C^n$
compatible with $J_0$ we write $|\Omega_0|_h:=|\Omega_0(E_1,\dots,E_n)|$, where
$(E_1,\dots,E_n)$ is any $h$-unitary frame, that is, $(E_1,J_0E_1,\dots,E_n,J_0E_n)$ is
$h$-orthonormal. The value does not depend on the frame, since $\Omega_0$ is of type
$(n,0)$ and a change of unitary frame multiplies it by a complex number of modulus one.
This is $2^{-n/2}$ times the norm induced by $h$ on complex $n$-forms. For
$g_0=\frac14|dw|^2$ the vectors $E_j=2\partial_{x_j}$ form a unitary frame with
$d\zeta_j(E_k)=\delta_{jk}$, so $|\Omega_0|_{g_0}=1$. For $A\subset\C^n$ and $\bar z\in\R$
we write $A\times\{\bar z\}:=\{(w,\bar z):w\in A\}$, and for a measure $m$ on $\C^n$ we
write $m\times\delta_{\bar z}$ for its push-forward under $w\mapsto(w,\bar z)$.

Let $X=\sum_j(x_j\partial_{x_j}+y_j\partial_{y_j})+2z\partial_z$ be the generator of
$\delta_\lambda$ and $X^{\mathcal H}=X-z\overline\xi$ its horizontal part, the horizontal
lift of the Euler field $\mathfrak r=\sum_j(x_j\partial_{x_j}+y_j\partial_{y_j})$. The
\emph{horizontal potential} is
\[
f=\tfrac14|X^{\mathcal H}|^2=\tfrac1{16}|w|^2,\qquad \overline\nabla f=\tfrac12X^{\mathcal
H},\qquad |\overline\nabla f|^2=f,\qquad \overline\nabla^2f(U,U)=\tfrac12|U|^2\ \text{for
horizontal }U
\]
\cite[(2.5), (2.6)]{ChangWuZhang2026ContactBlowup}. For $\hat y\in\R^{2n+1}$ put $f_{\hat
y}:=f\circ\ell_{\hat y}^{-1}$ and $z_{\hat y}:=z\circ\ell_{\hat y}^{-1}$. As in
\cite[\S\S2.1--2.2]{ChangWuZhang2026ContactBlowup}, $\overline\nabla$ denotes ambient
covariant derivatives (gradients, Hessians), taken for $g$, or for $\overline g$ in the
model; for a submanifold, $\nabla$, $\Delta$ and $\perp$ refer to the induced structure,
$B(U,W)=(\overline\nabla_UW)^\perp$ is the second fundamental form and $H=\tr B$. A vector
is called horizontal if it lies in $\cD$.

\subsection{The flow and the standing hypotheses}\label{subsec:hypotheses}
We write $\tau=T-t$ and, for $t\in[0,T)$, $\Sigma_t=F_t(\Sigma)$ with the induced area
measure $\mu_t$. The hypotheses (S) and (F) below are in force throughout the paper; the
hypotheses (I) and (O) are assumed only in the statements that invoke them explicitly.
\begin{enumerate}
\item[(S)] $M$ is a compact Sasakian manifold.
\item[(F)] $\Sigma$ is a closed $n$-manifold, and $F\colon\Sigma\times[0,T)\to M$, with
$T<\infty$, is a smooth, normally parametrised family of Legendrian immersions. There is a
smooth $\theta\colon\Sigma\times[0,T)\to\R$, the \emph{Legendrian angle}, with
\[
H=\Phi\nabla\theta,\qquad \partial_tF=H-2\theta\xi,\qquad
(\partial_t-\Delta)\theta=(K+2)\theta
\]
for a constant $K$.
\item[(I)] \emph{Type~I}: $|B|^2\le C_I/(T-t)$ on $\Sigma\times[0,T)$.
\item[(O)] \emph{Orientation}: $\Sigma$ is orientable.
\end{enumerate}
If $M$ is $\eta$-Einstein with constant $K$ and the mean curvature form of $\Sigma$ is
exact, then (F) is the Legendrian mean curvature flow
\cite[(1.1)]{ChangWuZhang2026ContactBlowup}, with $\theta$ normalised by
\cite[(2.2)]{ChangWuZhang2026ContactBlowup}, cf.\ \cite[Corollary~3.3]{s1}. It differs
from the flow of \cite[Definition~3.5]{s1} by a rotation along the Reeb flow
\cite[Section~2.2]{ChangWuZhang2026ContactBlowup}. In a general Sasakian manifold (F) is
an assumption. By the maximum principle $C_\theta:=\sup|\theta|\le
e^{|K+2|T}\max|\theta(\cdot,0)|<\infty$. Since $\partial_t\,d\mu=-|H|^2d\mu$ (the vertical
velocity is normal and $H\perp\xi$), $\mu_t(\Sigma)\le\mu_0(\Sigma)$.

\begin{remark}[Legendrian angle and geometric phase]\label{rem:potential-phase}
The function $\theta$ in (F) is the Legendrian angle of \cite[Definition~3.4]{s1} in the
time gauge fixed by its evolution equation, as in
\cite[(2.2)]{ChangWuZhang2026ContactBlowup}: a primitive of the mean curvature form,
$H=\Phi\nabla\theta$. Replacing $\theta$ by $\theta+c_0e^{(K+2)t}$ composes the flow with
the Reeb flow \cite[\S2.2]{ChangWuZhang2026ContactBlowup}. For an oriented Legendrian
immersion in a contact chart one also has the \emph{geometric phase} $\vartheta$ of the
projection with respect to a fixed holomorphic volume form (Lemma~\ref{lem:phase}). In the
model it is, up to an additive constant, the Legendrian angle $\theta^{\rm cqz}$ of
\cite{cqz}, which is defined through a holomorphic volume form; see
\cite[Remark~4.10]{ChangWuZhang2026ContactBlowup}. Only the differential
$d(\theta-\vartheta)$ is controlled (Lemma~\ref{lem:phase}); in a curved chart it is in
general nonzero, and $\theta-\vartheta$ is determined by it only up to an additive locally
constant function, which is not controlled. Example~\ref{ex:potential-not-phase} shows
that this locally constant part is not determined by (S), (F), (I), (C). Statements about
Legendrian angles in this paper concern $\theta$; geometric phases enter only where this
is said explicitly.
\end{remark}

\begin{example}\label{ex:potential-not-phase}
Let $\Gamma=\{(a,b,c):a,b\in\Z^n,\ c\in\frac12\Z\}$. By the formula for $\ell_{(a,b,c)}$,
$\Gamma$ is a discrete cocompact subgroup of $\Hn$. Since
$\overline\eta,\overline\xi,\overline\Phi,\overline g$ are left-invariant,
$M=\Gamma\backslash\Hn$ is a compact Sasakian manifold. For $\bar z\in\R$ the plane
$\{y=0,\ z=\bar z\}$ is Legendrian, since $\overline\eta$ vanishes on it, and minimal,
since it is the horizontal lift of the minimal Lagrangian plane $\R^n\subset\C^n$. It is
invariant exactly under the elements $(a,0,0)$ of $\Gamma$, so it descends to an embedded
Legendrian torus $T^n\subset M$.

Fix $\bar c\in\R$ and put $F_t(x)=[(x,0,4\bar c(T-t))]$, $\theta\equiv\bar c$, $K+2=0$.
Then $H=0=\Phi\nabla\theta$, $(\partial_t-\Delta)\theta=0$, and $\partial_tF=-4\bar
c\,\partial_z=-2\bar c\,\overline\xi$, so (F) holds. As contact normal chart at $p=[0]$ we
take the Heisenberg coordinates $(w,z)$ on $\{|w|<r_0,\ |z|<z_0\}$, with $r_0\le\frac12$
and $z_0\le\frac14$, so that this box injects into $M$; it satisfies
Lemma~\ref{lem:charts}(i)--(iv) with $h_{[0]}=g_0$ and $\beta_{[0]}=\overline\beta$. In
another admissible chart the limit plane below and its phase are rotated by an element of
$U(n)$. $B\equiv0$ gives (I), and the rescalings at $p=[0]$ converge smoothly to the plane
$\{y=0,\ z=4|s|\bar c\}$, on which $z=4|s|\theta$ with $\theta=\bar c$; this is the smooth
compactness (C) of Section~\ref{subsec:typeI}, and the exact centring of
Proposition~\ref{prop:centring}. With the coordinates $\zeta_k$ of
Section~\ref{subsec:sasakian}, the form $\Omega_0$ restricts to $\R^n$ as
$2^{-n}dx_1\wedge\dots\wedge dx_n$, so the geometric phase of the limit plane is $0$ or
$\pi$ mod $2\pi$, depending on the orientation. Since $\bar c$ is arbitrary, $\theta$ is
in general not a lift of the geometric phase.

Taking for $\Sigma$ two copies of $T^n$ with $\theta\equiv0$ and $\theta\equiv1$, both
reaching the same torus at time $T$, gives a limit with two components over the
\emph{same} plane. Their Legendrian angles differ by $1$. With any single choice of
holomorphic volume form, two orientations of the same plane have phases differing by a
multiple of $\pi$, so no such choice makes both Legendrian angles lifts of the geometric
phases.
\end{example}

\subsection{Contact normal coordinates}\label{subsec:charts}

\begin{lemma}\label{lem:charts}
There are $r_0,z_0>0$ and constants $C_k$, depending only on $(M,\eta,g)$, with the
following property. Every $p\in M$ has a neighbourhood $U_p$ and a diffeomorphism
\[
\Upsilon_p\colon U_p\to\{(w,z)\in\C^n\times\R:\ |w|<r_0,\ |z|<z_0\},\qquad
\Upsilon_p(p)=0,
\]
such that, identifying functions on $U_p$ with functions of $(w,z)$:
\begin{enumerate}[(i)]
\item $(\Upsilon_p)_*\xi=2\partial_z$;
\item $\eta=\frac12dz+\beta_p$, where $\beta_p$ is a $1$-form in $w$ alone in \emph{radial
gauge}, $\iota_{\mathfrak r}\beta_p=0$;
\item $w$ are holomorphic coordinates ($J=J_0$) of a local leaf space, in which the
transverse metric $g_Z$, written $h_p$, satisfies $h_p=g_0+O(|w|^2)$ and $\partial
h_p=O(|w|)$ (K\"ahler normal coordinates); $g(U^h,W^h)=h_p(U,W)$ for the horizontal lifts
$U^h,W^h$;
\item $|\beta_p-\overline\beta|\le C|w|^3$ and $|\partial(\beta_p-\overline\beta)|\le
C|w|^2$. All derivatives of $h_p,\beta_p$ are bounded by $C_k$, uniformly in $p$.
\end{enumerate}
\end{lemma}
\begin{proof}
For a single point $p$ these are the coordinates of
\cite[Lemma~4.1]{ChangWuZhang2026ContactBlowup}: K\"ahler normal coordinates $w$ on a
local leaf space centred at the image of $p$, the fibre coordinate $z$ with
$\xi=2\partial_z$ and $z(p)=0$, and the section chosen so that
$\eta-\frac12dz=\pi^*\beta_p$ with $\beta_p$ the radial primitive
\[
\beta_p|_w:=\int_0^1t\,(\omega_Z)_{tw}(w,\cdot)\,dt
\]
of $\omega_Z$. Here the first slot is the fixed vector $w=\mathfrak r(w)$, not $\mathfrak
r(tw)=tw$; then $d\beta_p=\omega_Z$ and $\iota_{\mathfrak r}\beta_p=0$.
By \cite[Lemma~4.1]{ChangWuZhang2026ContactBlowup}, $\eta-\overline\eta=O(|w|^3)|dw|$
and $g-\overline g=O(|w|^2)$, with the corresponding bounds for the derivatives. What has
to be added is the uniformity in $p$ and the shape of the domain, which are used for
moving centres (Lemma~\ref{lem:moving}(c)) and for constants that do not depend on the
centre.

Since $M$ is Sasakian, the Reeb field is a unit Killing field preserving $\eta$ and
$\Phi$, so the Reeb foliation is Riemannian. The normality of the structure makes the
transverse almost complex structure integrable, so that the foliation is transversally
K\"ahler \cite[Ch.~7]{bg}; this is the transverse K\"ahler structure recalled in
Section~\ref{subsec:sasakian}. For a slice $S$ transverse to $\xi$, the map
$(s,z)\mapsto\mathrm{Fl}^\xi_{z/2}(s)$, with $\mathrm{Fl}^\xi$ the flow of $\xi$, gives
coordinates in which $\xi=2\partial_z$. By compactness there are $z_0>0$ and finitely many
slices $S_\alpha$ with the following properties. The map
$(s,z)\mapsto\mathrm{Fl}^\xi_{z/2}(s)$ is a diffeomorphism on
$S_\alpha\times(-8z_0,8z_0)$. Each $S_\alpha$ carries holomorphic coordinates for the
transverse complex structure, in which the transverse metric is uniformly equivalent to
the Euclidean one and has bounded derivatives. Every $p\in M$ lies on the Reeb orbit of a
point $s\in S_\alpha$ at a uniform distance from $\partial S_\alpha$, and is reached from
$s$ by the Reeb flow in time less than $2z_0$. The Reeb flow preserves
$(\eta,\xi,\Phi,g)$; moving $S_\alpha$ by it, we obtain a slice through $p$ with the same
properties, and a flow box of height $(-4z_0,4z_0)$ around it.

In such a flow box, since $\mathcal L_\xi\eta=0$ and $\eta(\xi)=1$, the form
$\eta-\frac12dz$ is basic. Translate the holomorphic coordinates so that $p$ corresponds
to $0$, and compose with a complex linear map taking $g_Z(0)$ to $g_0$; this map and its
inverse are bounded independently of $p$. A holomorphic change of coordinates, adding to
$w$ terms quadratic in $w$, removes the linear part of $g_Z-g_Z(0)$; this uses
$d\omega_Z=0$, which makes the linear coefficients symmetric in the holomorphic indices.
These coefficients are bounded by the $C^1$ norm of $g_Z$, so by the inverse function
theorem the change is biholomorphic on a ball whose radius, like all the $C^k$ bounds,
does not depend on $p$. Hence $h_p=g_0+O(|w|^2)$, and $d\eta=\omega_Z=\omega_0+O(|w|^2)$
with all derivatives, uniformly in $p$.

With $\beta_p$ as above, $d(\eta-\frac12dz-\beta_p)=0$, so there is a function $G_p$ on
the ball with $\eta-\frac12dz-\beta_p=-dG_p$ and $G_p(0)=0$. As $|dG_p|\le C_G$ with $C_G$
independent of $p$, we have $|G_p|\le C_Gr_0$ on $\{|w|<r_0\}$; we choose $r_0$ so small
that also $2C_Gr_0\le z_0$, and that $h_p\ge\frac12g_0$ on $\{|w|<r_0\}$. Let $U_p$ be the
preimage of $\{|w|<r_0,\ |z'|<z_0\}$ under $(w,z)\mapsto(w,z')=(w,z-2G_p(w))$. It lies in
the flow box, since $|z|<z_0+2C_Gr_0\le2z_0$ on it. Replacing $z$ by $z'$ (a change of
slice) gives (i)--(ii), still with $\xi=2\partial_z$ and $z(p)=0$. Finally, (iv) follows
from the same homotopy formula applied to
$d\eta-d\overline\eta=\omega_Z-\omega_0=O(|w|^2)$, since
$\overline\beta=\frac12\iota_{\mathfrak r}\omega_0$ is the radial primitive of
$\omega_0=d\overline\eta$. This gives
$(\beta_p-\overline\beta)|_w=\int_0^1t\,(\omega_Z-\omega_0)_{tw}(w,\cdot)\,dt=O(|w|^3)$,
and differentiating under the integral gives the derivative bounds.
\end{proof}

We record from the proof that, by the choice of $r_0$,
\begin{equation}\label{eq:hp-lower}
h_p\ge\tfrac12g_0\quad\text{on }\{|w|<r_0\},\ \text{for every }p\in M;
\end{equation}
this is preserved when $r_0$ is decreased.

\emph{Cut-offs.} Fix once and for all $\zeta_*\in C^\infty_c(\{|w|<r_0,\ |z|<z_0\})$ with
$0\le\zeta_*\le1$ and $\zeta_*=1$ on $\{|w|<\frac{r_0}2,\ |z|<\frac{z_0}2\}$. For $q\in M$
let $\zeta_q:=\zeta_*\circ\Upsilon_q$ on $U_q$, extended by $0$, and $\psi_q:=\zeta_q^2$.
In the chart $\Upsilon_q$ the metric $g$ is given by $h_q$ on horizontal vectors and by
$\eta=\frac12dz+\beta_q$. By Lemma~\ref{lem:charts}(iii)--(iv), and since
$h_q\ge\frac12g_0$ by \eqref{eq:hp-lower}, its coefficients, those of its inverse, and all
their derivatives are bounded independently of $q$. Hence $|\overline\nabla\zeta_q|$ and
$|\overline\nabla^2\zeta_q|$ are bounded independently of $q$.

\subsection{The projected area ratio bound}\label{subsec:D}
The following property is not an additional assumption: by Lemma~\ref{lem:D} it follows
from (S) and (F).
\begin{enumerate}
\item[(D)] \emph{Projected area ratio bound.} With $r_0,z_0$ from Lemma~\ref{lem:charts},
$\mu_t\big(F_t^{-1}(U_p\cap\{|w|<R\})\big)\le D_0R^n$ for all $p\in M$, $0<R\le r_0$ and
$t\in[0,T)$, the heights being unrestricted inside $U_p$.
\end{enumerate}

\begin{lemma}\label{lem:D}
Assume (S) and (F). Then (D) holds, with $D_0$ depending only on $(M,\eta,g)$, $T$,
$C_\theta$ and the initial immersion $F_0$.
\end{lemma}
\begin{proof}
By \cite[Lemma~4.15]{ChangWuZhang2026ContactBlowup} there are $C_A$ and $s_0>0$ with
$\mu_t(F_t^{-1}(B_s(x)))\le C_As^n$ for all $x\in M$, $0<s\le s_0$ and $t\in[0,T)$, where
$B_s(x)$ is the geodesic ball. Its proof uses only the flow equation
$\partial_tF=H-2\theta\xi$, the bound $|\theta|\le C_\theta$, an isometric embedding of
the compact manifold $M$, and the closedness of $\Sigma$, all of which are available under
(S) and (F). The same is true of \cite[Proposition~4.16]{ChangWuZhang2026ContactBlowup},
whose proof uses in addition only \cite[Lemma~4.15]{ChangWuZhang2026ContactBlowup} and the
pointwise computation of \cite[Proposition~4.5]{ChangWuZhang2026ContactBlowup}. Every
point of $M$ is the centre of a flow box $U\cong\{|w|<2r\}\times(-2\delta,2\delta)$ in the
coordinates of \cite[Lemma~4.1]{ChangWuZhang2026ContactBlowup}. By the last assertion of
\cite[Proposition~4.16]{ChangWuZhang2026ContactBlowup} there are $C_P$ and $r_1>0$ such
that, for all $t\in[0,T)$, the part of $\Sigma_t$ lying in $\{|w|<r,\ |z|<\delta\}$,
pushed forward under $w$, has mass at most $C_Pr'^n$ in every ball $B_{r'}(w_*)$ with
$r'\le r_1$ and $|w_*|\le r/2$; here all heights $|z|<\delta$ are counted.

To compare these flow boxes with the charts of Lemma~\ref{lem:charts}, cover $M$ by the
inner boxes $\{|w|<r/4,\ |z|<\delta/4\}$ of finitely many such flow boxes, and fix $C_P$,
$r_1$, $r$, $\delta$ for all of them. Let $p\in M$ lie in the inner box of the flow box
$U$. Since the coefficients of $g$ in the chart $\Upsilon_p$ are bounded independently of
$p$, the set $U_p$ lies in the $g$-ball of radius $C'(r_0+z_0)$ about $p$, with $C'$
independent of $p$; and the closed inner box of each of the finitely many flow boxes has a
positive $g$-distance from the complement of that flow box. After a first decrease of
$r_0$ and $z_0$, admissible for the reasons given below, we may therefore assume
$U_p\subset U$, so that $w\circ w_p^{-1}$ is defined on $\{|w_p|<r_0\}$. The horizontal
coordinates $w_p$ of the chart $\Upsilon_p$ of Lemma~\ref{lem:charts}
and the coordinates $w$ of $U$ are both holomorphic coordinates of the transverse
structure in which the transverse metric has uniformly bounded coefficients, together with
their inverses and derivatives (Lemma~\ref{lem:charts}(iii)--(iv) and
\eqref{eq:hp-lower}). Hence $w\circ w_p^{-1}$ is Lipschitz on the convex set
$\{|w_p|<r_0\}$, with a constant $C\ge1$ independent of $p$. The two height coordinates
differ by $z(p)$ and a change of slice of
size at most $Cr_0$. After decreasing $r_0$ and $z_0$ in Lemma~\ref{lem:charts}, which
preserves its conclusions, keeps $r_0,z_0$ dependent only on $(M,\eta,g)$, and can be done
before $\zeta_*$ is fixed, we may therefore assume $Cr_0\le r/4$ and
$z_0+Cr_0\le\delta/2$.

Then $U_p\subset\{|w|<r,\ |z|<\delta\}$, and for $R\le r_0$ the set
$U_p\cap\{|w_p|<R\}$ lies over the ball $B_{CR}(w(p))$, with $|w(p)|<r/4$ and $CR\le r/4$.
If $CR\le r_1$, this gives $\mu_t(F_t^{-1}(U_p\cap\{|w_p|<R\}))\le C_PC^nR^n$. If
$CR>r_1$, the left-hand side is at most $\mu_0(\Sigma)\le\mu_0(\Sigma)(C/r_1)^nR^n$. This
is (D) with $D_0=\max\{C_PC^n,\mu_0(\Sigma)(C/r_1)^n\}$.
\end{proof}

\subsection{The model identities, moving centres and the geometric
phase}\label{subsec:model-identities}

\begin{lemma}\label{lem:model}
In $U_p$ put $f_p:=|w|^2/16$ and $z_p:=z$. For a Legendrian plane $S\subset\cD$ let $P_S$
be the orthogonal projection onto $S$, and for a function $\gamma$ put
$\nabla^S\gamma:=P_S\overline\nabla\gamma$; for $S=dF(T_x\Sigma)$ this is the gradient of
$\gamma\circ F$ on $\Sigma$. There is a constant $C$, depending only on $(M,\eta,g)$ and
not on $p$, such that, for horizontal $U,W$ at points of $U_p$ and every Legendrian plane
$S\subset\cD$ there:
\begin{enumerate}[(a)]
\item $\xi f_p=0$, $\xi z_p=2$, and $\overline\nabla f_p$ is horizontal;
\item $|\overline\nabla^2f_p(U,W)-\frac12g(U,W)|\le C|w|^2|U||W|$ and
$\big||\overline\nabla f_p|^2-f_p\big|\le C|w|^4$;
\item $|\overline\nabla^2z_p(U,W)|\le C|w|^2|U||W|$;
\item the horizontal part $\overline\nabla^{\mathcal H}z_p$ of $\overline\nabla z_p$
satisfies $\overline\nabla^{\mathcal H}z_p=4\Phi\overline\nabla f_p+e_H$ with $|e_H|\le
C|w|^3$. Hence $\nabla^Sz_p=4\Phi\big((\overline\nabla f_p)^{\perp_S}\big)+P_Se_H$,
$\big||\nabla^Sz_p|^2-16|(\overline\nabla f_p)^{\perp_S}|^2\big|\le C|w|^4$, and
$|\overline\nabla^{\mathcal H}z_p|\le C|w|$. In the model ($e_H=0$) one has
$\nabla^Sz=4\overline\Phi\big((\overline\nabla f)^{\perp_S}\big)$ exactly.
\item Along a flow as in (F), on $F_t^{-1}(U_p)$, put
$\rho=\rho_{p,T}:=(4\pi\tau)^{-n/2}e^{-f_p/\tau}$ and $\mathcal S_p:=H+(\overline\nabla
f_p)^\perp/\tau$. Then
\begin{gather*}
(\partial_t+\Delta-|H|^2)\rho=-|\mathcal S_p|^2\rho+e_\rho\,\rho,\qquad |e_\rho|\le
C\Big(\frac{|w|^2}\tau+\frac{|w|^4}{\tau^2}\Big),\\
(\partial_t-\Delta)z_p=-4\theta+e_z,\qquad |e_z|\le C|w|^2 .
\end{gather*}
\item Along a flow as in (F), for $\ell\ge0$ and $0<\tau\le1$:
$\int_{F_t^{-1}(U_p)}(|w|^2/\tau)^\ell\rho\,d\mu\le C(n,\ell)D_0$, with $D_0$ from (D).
\end{enumerate}
\end{lemma}
\begin{proof}
\emph{The submersion.} The projection $(w,z)\mapsto w$ is a Riemannian submersion onto
$(\{|w|<r_0\},h_p)$ with totally geodesic fibres. Let $\nabla^Z$ be the Levi-Civita
connection of $h_p$, and $U^h$ the horizontal lift of a vector $U$ of the base. By
O'Neill's formula, as in \cite[\S2.3 and Lemma~2.3]{ChangWuZhang2026ContactBlowup}, for
horizontal lifts of base fields
\[
\overline\nabla_{U^h}W^h=(\nabla^Z_UW)^h+\tfrac12[U^h,W^h]^v,\qquad
\tfrac12[U^h,W^h]^v=\tfrac12\eta([U^h,W^h])\,\xi=-\tfrac12d\eta(U^h,W^h)\,\xi ;
\]
the last term is the O'Neill tensor $A_{U^h}W^h$ of
\cite[\S2.3]{ChangWuZhang2026ContactBlowup}. Below we write $U$ for $U^h$. All the error
terms below are estimated through the expansions $h_p=g_0+O(|w|^2)$, $\partial
h_p=O(|w|)$, $\beta_p-\overline\beta=O(|w|^3)$ of Lemma~\ref{lem:charts}, together with
\eqref{eq:hp-lower}, whose constants depend only on $(M,\eta,g)$; this gives the
independence of $C$ from $p$.

\emph{(a), (b).} (a) is clear. $f_p$ is basic, so its Hessian on horizontal vectors is the
$h_p$-Hessian of $|w|^2/16$, namely
$\frac18\delta_{ab}-\Gamma^c_{ab}\partial_cf_p=\frac12(h_p)_{ab}+O(|w|^2)$, because
$\Gamma=O(|w|)$ and $\partial f_p=O(|w|)$. Also $|\overline\nabla
f_p|^2=h_p^{ab}\partial_af_p\partial_bf_p=(4\delta^{ab}+O(|w|^2))\frac{w_aw_b}{64}=f_p+O(|w|^4)$.

\emph{(c).} For horizontal $W$, $\eta(W)=0$ gives $Wz_p=-2\beta_p(W)$. Hence
\[
\overline\nabla^2z_p(U,W)=-2U(\beta_p(W))+2\beta_p(\nabla^Z_UW)-\tfrac12\eta([U,W])\,\xi
z_p=-2(\nabla^Z_U\beta_p)(W)+d\eta(U,W).
\]
With $d\eta=d\beta_p$, this is $-(\nabla^Z_U\beta_p)(W)-(\nabla^Z_W\beta_p)(U)$, the
symmetrised derivative. For $\overline\beta$ and the flat metric $g_0$ it vanishes. Since
$\partial\beta_p=\partial\overline\beta+O(|w|^2)$, $\Gamma=O(|w|)$ and $\beta_p=O(|w|)$,
it is $O(|w|^2)$.

\emph{(d).} Since $Wz_p=-2\beta_p(W)$ for horizontal $W$, $\overline\nabla^{\mathcal
H}z_p$ is the horizontal lift of the $h_p$-dual of $-2\beta_p$. In the model,
$g_0(4J_0\nabla^{g_0}f,W)=-2\overline\beta(W)$, with $\nabla^{g_0}$ the $g_0$-gradient;
equivalently, $\overline\nabla z=4\overline\Phi\,\overline\nabla f+2\overline\xi$, as in
the proof of \cite[Lemma~2.5]{ChangWuZhang2026ContactBlowup}. Here $J=J_0$,
$h_p=g_0+O(|w|^2)$, and $\overline\nabla f_p$ is the horizontal lift of
$\nabla^{h_p}f_p=\nabla^{g_0}f+O(|w|^3)$; hence
$h_p(4J\nabla^{h_p}f_p,W)=-2\overline\beta(W)+O(|w|^3)|W|=-2\beta_p(W)+O(|w|^3)|W|$, by
Lemma~\ref{lem:charts}(iv). For a Legendrian plane, $\Phi S=S^\perp\cap\cD$ and $\Phi$ is
an isometry of $\cD$, so $P_S\Phi=\Phi P_{S^\perp}$ on $\cD$. The remaining statements
follow from $|\overline\nabla f_p|\le C|w|$.

\emph{(e).} We have $\partial_t(f_p\circ F)=\langle\overline\nabla f_p,H\rangle$ by (a),
and $\Delta(f_p\circ F)=\tr_\Sigma\overline\nabla^2f_p+\langle\overline\nabla
f_p,H\rangle$ (the mean curvature vector of a Legendrian is horizontal). Hence
\[
(\partial_t+\Delta-|H|^2)\rho=\rho\Big[\frac{n/2-\tr_\Sigma\overline\nabla^2f_p}{\tau}-\frac{f_p-|\overline\nabla
f_p|^2}{\tau^2}-\Big|H+\frac{(\overline\nabla f_p)^\perp}\tau\Big|^2\Big],
\]
using $f_p-|\nabla f_p|^2=(f_p-|\overline\nabla f_p|^2)+|(\overline\nabla f_p)^\perp|^2$,
where $\nabla f_p$ is the gradient on $\Sigma$. Apply (b). In the model the first two
terms of the bracket vanish; this is the computation in the proof of
\cite[Proposition~4.3]{ChangWuZhang2026ContactBlowup}. Similarly,
$\partial_tz_p=dz_p(H)-2\theta\,\xi z_p=dz_p(H)-4\theta$ and $\Delta
z_p=\tr_\Sigma\overline\nabla^2z_p+dz_p(H)$; apply (c).

\emph{(f).} Bound $(|w|^2/\tau)^\ell e^{-|w|^2/(16\tau)}\le C_\ell e^{-|w|^2/(32\tau)}$,
and use the layer-cake formula with $m(R)=\mu_t(F_t^{-1}(U_p\cap\{|w|<R\}))\le
D_0\min(R,r_0)^n$.
\end{proof}

\begin{remark}[The round sphere $S^3$]\label{rem:S3}
For the standard Sasakian $S^3$ the base is $\C P^1$. By homogeneity every $p$ has a
contact normal chart with
\[
h_p=\frac{|dw|^2}{4(1+c_S|w|^2)^2},\qquad d\eta=\omega_Z=\frac{dx\wedge
dy}{2(1+c_S|w|^2)^2},
\]
where $c_S>0$ is fixed by the normalisation (the base has constant curvature $16c_S$). The
radial gauge gives $\beta_p=\overline\beta/(1+c_S|w|^2)$. A direct computation gives:
\begin{itemize}
\item $\overline\nabla^2f_p-\frac12h_p=O(|w|^2)$ on horizontal vectors, with leading term
of exact order $c_S|w|^2$;
\item $|\overline\nabla f_p|^2-f_p=\frac{c_S|w|^4}{16}(2+c_S|w|^2)$;
\item the eigenvalues of $h_p^{-1}\overline\nabla^2z_p$ (horizontal Hessian) are
$\pm2c_S|w|^2$;
\item $|\overline\nabla^{\mathcal H}z_p-4\Phi\overline\nabla f_p|=c_S|w|^3$ exactly.
\end{itemize}
All of these vanish for $c_S=0$. The orders of Lemma~\ref{lem:model}(b)--(d) are therefore
attained. The verification is given in Appendix~\ref{app:B}.
\end{remark}

\begin{lemma}\label{lem:moving}
Let $\Upsilon_p$ be as in Lemma~\ref{lem:charts}, and $\kappa=\lambda^2$.
\begin{enumerate}[(a)]
\item \emph{Naturality.} The map $\delta_\lambda\circ\Upsilon_p$ is a contact normal chart
at $p$ for the $\cD$-homothetic structure $(\eta_\kappa,\xi_\kappa,\Phi,g_\kappa)$ of
Section~\ref{subsec:sasakian}, with $g_\kappa=\kappa g^T+\kappa^2\eta\otimes\eta$, defined
on $\{|w|<\lambda r_0,\ |z|<\lambda^2z_0\}$.
\item \emph{Closeness to the model.} In these coordinates the structure
$\sigma_\lambda:=(\delta_\lambda)_*(\eta_\kappa,\xi_\kappa,\Phi,g_\kappa)$ has Reeb field
$\overline\xi=2\partial_z$, and differs from $(\overline\eta,\overline\Phi,\overline g)$
by at most $C_k(R)\lambda^{-2}$ in $C^k(\{|w|<R,|z|<R^2\})$, for $\lambda\ge
R\max\{r_0^{-1},z_0^{-1/2}\}$.
\item \emph{Moving centres, modulo $U(n)$.} Let $R\ge1$, let $q\in U_p$ with $\hat y=(\hat
w,\hat z):=\delta_\lambda\Upsilon_p(q)$, $|\hat w|+|\hat z|^{1/2}\le R$, and let
$\lambda\ge\lambda_*(R)$. Then there is $U\in U(n)$, depending on $p,q,\lambda$ and on the
choice of the charts, such that
\[
\big\|\delta_\lambda\Upsilon_q\Upsilon_p^{-1}\delta_\lambda^{-1}-\tilde U\circ\ell_{\hat
y}^{-1}\big\|_{C^k(\{|w|<R,|z|<R^2\})}\le C_k(R)\,\lambda^{-2},\qquad \tilde
U(w,z):=(Uw,z).
\]
Consequently the rescaled functions
\[
f_q^{(\lambda)}:=\lambda^2\,f_q\circ\Upsilon_p^{-1}\circ\delta_\lambda^{-1},\qquad
z_q^{(\lambda)}:=\lambda^2\,z_q\circ\Upsilon_p^{-1}\circ\delta_\lambda^{-1},\qquad
\psi_q^{(\lambda)}:=\psi_q\circ\Upsilon_p^{-1}\circ\delta_\lambda^{-1},
\]
defined on $\delta_\lambda\Upsilon_p(U_p\cap U_q)$, where $\psi_q$ is the cut-off defined
after Lemma~\ref{lem:charts}, converge in $C^k_{\rm loc}$, uniformly in such $\hat y$, to
the model functions $f_{\hat y}$, $z_{\hat y}$ and to $1$. The functions $f_q$ and $z_q$
carry the factor $\lambda^2$ because $f$ and $z$ have weight $2$ under $\delta_\lambda$;
the cut-off has weight $0$.
\end{enumerate}
\end{lemma}
\begin{proof}
\emph{(a).} If $w$ are K\"ahler normal coordinates for $h$ with $h(0)=g_0$, then $\lambda
w$ are K\"ahler normal coordinates for $\kappa h$ with the same normalisation. The radial
primitive commutes with linear scalings, $z\mapsto\lambda^2z$ matches
$\eta_\kappa=\lambda^2\eta$, and
$(\delta_\lambda)_*\xi_\kappa=\lambda^{-2}(\delta_\lambda)_*(2\partial_z)=2\partial_z$.

\emph{(b).} Let $\delta^{\rm b}_\lambda(w)=\lambda w$ be the base part of
$\delta_\lambda$, and put $h^{(\lambda)}_p:=\kappa\big((\delta^{\rm
b}_\lambda)^{-1}\big)^*h_p$ and $\beta^{(\lambda)}_p:=\kappa\big((\delta^{\rm
b}_\lambda)^{-1}\big)^*\beta_p$; in coordinates
$(h^{(\lambda)}_p)_{ab}(w)=(h_p)_{ab}(w/\lambda)$. By Lemma~\ref{lem:charts}(iii)--(iv),
$h_p=g_0+O(|w|^2)$ and $\beta_p=\overline\beta+O(|w|^3)$, so
$h^{(\lambda)}_p=g_0+O(\lambda^{-2}|w|^2)$ and
$\beta^{(\lambda)}_p=\overline\beta+O(\lambda^{-2}|w|^3)$, with the same gain for every
derivative. For the second derivatives of $\beta_p$ we use
$|\partial^2(\beta_p-\overline\beta)|\le C|w|$, which follows from
$|\beta_p-\overline\beta|\le C|w|^3$ and the uniform $C^3$ bound of
Lemma~\ref{lem:charts}(iv) by Taylor's theorem; since a $k$-th derivative acquires the
factor $\lambda^{1-k}$ under the rescaling of $\beta_p$, this gives the rate
$\lambda^{-2}$ also for the second derivatives on fixed compact sets. In the rescaled
chart, $\sigma_\lambda$ is given by
$\frac12dz+\beta^{(\lambda)}_p$, $2\partial_z$, the horizontal lift of $J_0$, and
$\pi^*h^{(\lambda)}_p+\eta_\kappa\otimes\eta_\kappa$, and the horizontal lifts
$v-2\beta^{(\lambda)}_p(v)\partial_z$ are those of the model up to the error in
$\beta^{(\lambda)}_p$. This is the quantitative form of the last assertion of
\cite[Lemma~4.1]{ChangWuZhang2026ContactBlowup}.

\emph{(c).} K\"ahler normal coordinates are not unique: they can be changed by $U(n)$ and
by holomorphic terms of order $\ge3$. Independently chosen charts at $p$ and $q$ may
therefore differ by an $O(1)$ unitary rotation, and the transition map is not close to
$\ell_{\hat y}^{-1}$ in general; this is the reason for the factor $\tilde U$. The
function $f$, the form $\overline\beta$ and hence the vertical coordinate are
$U(n)$-invariant, which is why the scalar consequences hold.

For $\lambda\ge\lambda_*(R)$ the transition $\delta_\lambda\Upsilon_q\Upsilon_p^{-1}
\delta_\lambda^{-1}$ is defined on $\{|w|<3R,\ |z|<9R^2\}$: by Lemma~\ref{lem:charts} and
\eqref{eq:hp-lower} the domain of $\Upsilon_q$ contains a $g$-ball about $q$ of radius
independent of $q$, while the preimage of $\{|w|<3R,\ |z|<9R^2\}$ under
$\delta_\lambda\Upsilon_p$ contains $q$ and has $g$-diameter $O(R/\lambda)$. For its base
part, let $h^{(\lambda)}_p$ be the rescaled base metric of the $p$-chart, as in (b), and
define $h^{(\lambda)}_q$ in the same way. By (b), applied at $p$ and at $q$ with the
uniform constants of Lemma~\ref{lem:charts}, both are $O(\lambda^{-2})$-close to $g_0$ in
$C^k$ on balls of radius $3R$. The base part $w\mapsto\tilde w(w)$ of the transition is
holomorphic, maps $\hat w$ to $0$, and is an isometry from $h^{(\lambda)}_p$ to
$h^{(\lambda)}_q$, since both describe the same transverse metric. Its differential at
$\hat w$ is complex linear and maps $h^{(\lambda)}_p(\hat w)=g_0+O(\lambda^{-2})$ to
$h^{(\lambda)}_q(0)=g_0+O(\lambda^{-2})$; by polar decomposition it equals
$U+O(\lambda^{-2})$ for some $U\in U(n)$. The isometry equation
\[
\partial_a\partial_b\tilde w^c=\Gamma[h^{(\lambda)}_p]^e_{ab}\,\partial_e\tilde
w^c-\Gamma[h^{(\lambda)}_q]^c_{de}(\tilde w)\,\partial_a\tilde w^d\,\partial_b\tilde w^e
\]
has right-hand side $O(\lambda^{-2})$, because both Christoffel symbols are
$O(\lambda^{-2})$ and $|\partial\tilde w|$ is bounded (isometry). Integrating along
segments from $\hat w$ gives $\|\tilde w-U(\cdot-\hat w)\|_{C^1(B_{2R}(\hat w))}\le
C(R)\lambda^{-2}$, and differentiating the equation gives the same for all $C^k$ norms.

For the vertical coordinate we use that in both rescaled charts the Reeb field of
$\eta_\kappa$ is $2\partial_z$, so that the transition has the form $(w,z)\mapsto(\tilde
w(w),\,\tilde z)$, $\tilde z=z-\hat z+2G(w)$, with $G(\hat w)=0$. Let $\tilde\beta_q$ be
the form $\beta^{(\lambda)}_q$ of the rescaled $q$-chart. Comparing
$\eta_\kappa=\frac12dz+\beta^{(\lambda)}_p=\frac12d\tilde z+\tilde w^*\tilde\beta_q$ gives
$dG=\beta^{(\lambda)}_p-\tilde w^*\tilde\beta_q$. By (b),
$\tilde\beta_q=\overline\beta+O(\lambda^{-2})$ and
$\beta^{(\lambda)}_p=\overline\beta+O(\lambda^{-2})$ in $C^k$. Since $\tilde
w=U(\cdot-\hat w)+O(\lambda^{-2})$ and $U^*\overline\beta=\overline\beta$,
\[
dG=\overline\beta-\overline\beta(\cdot-\hat w)+O(\lambda^{-2})=\tfrac14\sum_j(-\hat
b_j\,dx_j+\hat a_j\,dy_j)+O(\lambda^{-2}),\qquad \hat w=(\hat a,\hat b),
\]
where $\overline\beta(\cdot-\hat w)$ is the pull-back of $\overline\beta$ under $w\mapsto
w-\hat w$. Integrating from $\hat w$ gives $2G=-\frac12(\hat b\cdot x-\hat a\cdot
y)+O(\lambda^{-2})$, so the vertical component of the transition is the vertical component
of $\ell_{\hat y}^{-1}$ up to $O(\lambda^{-2})$ in $C^k$. Together with the base part this
proves the estimate. For the scalar consequences, note that in the rescaled
$q$-coordinates $(\tilde w,\tilde z)=\delta_\lambda\Upsilon_q$ one has
$\lambda^2f_q=|\tilde w|^2/16$ and $\lambda^2z_q=\tilde z$, since $f_q=|w_q|^2/16$ and
$z_q$ is the height in the $q$-chart. Composing with the transition, and using $|Uw|=|w|$,
gives $f_q^{(\lambda)}=f_{\hat y}+O(\lambda^{-2})$ and $z_q^{(\lambda)}=z_{\hat
y}+O(\lambda^{-2})$ in $C^k$ on compact sets. The cut-off $\psi_q$ equals $1$ on
$\{|w|<\frac{r_0}2,\ |z|<\frac{z_0}2\}$ in the $q$-chart, that is, on $\{|w|<\frac{\lambda
r_0}2,\ |z|<\frac{\lambda^2z_0}2\}$ in the rescaled $q$-chart. A compact set lies in
$\{|w|<R,\ |z|<R^2\}$ for some $R$, and by the $C^0$ part of the estimate the transition
maps it into $\{|w|<C(R),\ |z|<C(R)\}$. For $\lambda$ large this set lies in the box
above, so $\psi_q^{(\lambda)}=1$ on the compact set, uniformly in $\hat y$.
\end{proof}

In the next lemma, the chart is the chart $\Upsilon_p$ on $U_p$, and an immersion $F$ into
$M$ is considered on $F^{-1}(U_p)$, where $\pi F$ is defined; all its assertions are
pointwise and refer to this set.

\begin{lemma}\label{lem:phase}
In a contact normal chart, the coordinates $w$ of Lemma~\ref{lem:charts}(iii) are
$J$-holomorphic in the sense of Section~\ref{subsec:sasakian}, and $\Omega_0$ denotes the
form of Section~\ref{subsec:sasakian} in these coordinates. For a function $u$ put
$d^cu:=-du\circ J$. Let $\pi$ be the projection $(w,z)\mapsto w$ of the chart.

Let $F$ be an oriented Legendrian immersion with $H=\Phi\nabla\theta$, and define the
geometric phase $\vartheta_p\colon\Sigma\to\R/2\pi\Z$ by
\[
(\pi F)^*\Omega_0=\big(|\Omega_0|_{h_p}\circ\pi F\big)\,e^{i\vartheta_p}\,\mathrm{vol},
\]
where $\mathrm{vol}$ is the volume form of the induced metric and the orientation. Then
\[
d(\theta-\vartheta_p)=-(\pi F)^*d^c\log|\Omega_0|_{h_p},\qquad\text{hence}\qquad
|d(\theta-\vartheta_p)|\le C|w| .
\]
After contact rescaling by $\lambda$ the bound becomes $C\lambda^{-2}|w'|$. The lemma
controls the differential only; the locally constant part of $\theta-\vartheta_p$ is not
determined (Example~\ref{ex:potential-not-phase}).
\end{lemma}
\begin{proof}
By \cite[Lemma~2.3]{ChangWuZhang2026ContactBlowup}, whose proof is pointwise and applies
in the chart \cite[\S4.6]{ChangWuZhang2026ContactBlowup}, $\pi F$ is an isometric
Lagrangian immersion into $(\{|w|<r_0\},h_p)$ with second fundamental form $\pi_*B$ and
mean curvature vector $\pi_*H$, and $\pi_*$ intertwines $\Phi$ with $J$. The Levi-Civita
connection $\nabla^Z$ of the K\"ahler metric $h_p$ induces the Chern connection on the
canonical bundle, and for the holomorphic section $\Omega_0$ this gives
$\nabla^Z_U\Omega_0=2\partial\log|\Omega_0|_{h_p}(U)\,\Omega_0=\big(d\log|\Omega_0|_{h_p}(U)+i\,d^c\log|\Omega_0|_{h_p}(U)\big)\Omega_0$.

Let $E_1,\dots,E_n$ be an oriented orthonormal frame of $T(\pi F)$, a unitary frame of
$\C^n$. Then $\Omega_0(E)=|\Omega_0|_{h_p}e^{i\vartheta_p}$ and
\[
U(\Omega_0(E))=(\nabla^Z_U\Omega_0)(E)+\sum_j\Omega_0(E_1,\dots,\nabla^Z_UE_j,\dots,E_n).
\]
The tangential parts of $\nabla^Z_UE_j$ cancel, by antisymmetry of the connection matrix.
The normal part contributes
$i\sum_j\langle\nabla^Z_UE_j,JE_j\rangle\Omega_0(E)=i\langle\pi_*H,JU\rangle\Omega_0(E)$,
by the symmetry of $\langle\pi_*B(\cdot,\cdot),J\cdot\rangle$. Comparing imaginary parts
gives
\[
d\vartheta_p=\langle\pi_*H,J\cdot\rangle+d^c\log|\Omega_0|_{h_p},\qquad
\langle\pi_*H,J\pi_*U\rangle=\langle\Phi\nabla\theta,\Phi
U\rangle=\langle\nabla\theta,U\rangle .
\]
Finally, in K\"ahler normal coordinates and by \eqref{eq:hp-lower}, $\det h_p=\det
g_0\,(1+O(|w|^2))$, so $\log|\Omega_0|_{h_p}=O(|w|^2)$ and $|d\log|\Omega_0|_{h_p}|\le
C|w|$.
\end{proof}

\begin{remark}[Relation to
{\cite[Lemma~4.1]{ChangWuZhang2026ContactBlowup}}]\label{rem:charts-M} At a single point,
the coordinates of Lemma~\ref{lem:charts} are those of
\cite[Lemma~4.1]{ChangWuZhang2026ContactBlowup}, and Lemma~\ref{lem:moving}(a), (b) are
the quantitative form, with the rate $\lambda^{-2}$, of the convergence
$\sigma_\lambda\to(\overline\eta,\overline\xi,\overline\Phi,\overline g)$ stated there.
Lemma~\ref{lem:charts} adds the uniformity in the centre, on boxes of fixed size
$r_0,z_0$; it is used for moving centres (Lemma~\ref{lem:moving}(c)), for the cut-offs
$\psi_q$, and for constants that do not depend on the centre. In the model
($e_H=e_\rho=e_z=0$) the identities (a)--(e) of Lemma~\ref{lem:model} are those of
\cite{ChangWuZhang2026ContactBlowup}: (a) follows from \cite[(2.5),
(2.6)]{ChangWuZhang2026ContactBlowup} and $\overline\xi=2\partial_z$, (b) is
\cite[(2.6)]{ChangWuZhang2026ContactBlowup}, (c) and (d) are in the proofs of
\cite[Lemmas~2.4 and~2.5]{ChangWuZhang2026ContactBlowup}, and (e) is
\cite[Proposition~4.3]{ChangWuZhang2026ContactBlowup} together with the evolution
$(\partial_t-\Delta)z=-4\theta$ in the proof of
\cite[Lemma~2.5]{ChangWuZhang2026ContactBlowup}. The estimate (f) is proved above from
(D).

In a curved manifold, \cite{ChangWuZhang2026ContactBlowup} uses the horizontal
potential of an isometric embedding of a local leaf space, which is constant along the
Reeb orbits. For the vertical cut-off in the proof of
\cite[Proposition~4.16]{ChangWuZhang2026ContactBlowup}, bounds on $|\nabla z|$ and
$(\partial_t-\Delta)z$ suffice, and no vanishing order at the centre is required. The
weighted kernel of
Section~\ref{sec:mono} depends on $z_p$, so the errors in (c), (d) and the term $e_z$ in
(e) have to vanish at the centre to the orders stated; this is where the radial gauge
$\iota_{\fr}\beta_p=0$ of Lemma~\ref{lem:charts}(ii) enters. The bound on $e_\rho$ in (e)
follows from (b), which involves only $f_p$ and the transverse metric $h_p$ in K\"ahler
normal coordinates, and does not depend on the choice of gauge.
\end{remark}

\subsection{Contact parabolic rescalings}\label{subsec:rescalings}
For $\lambda_i\to\infty$ define the contact parabolic rescalings at $(p,T)$
\cite[\S4.1]{ChangWuZhang2026ContactBlowup}, for $s\in[-S_0,0)$, where $S_0>0$ is fixed
and $i$ is so large that $S_0\lambda_i^{-2}\le\tau_0$, with
\begin{equation}\label{eq:tau0}
\tau_0:=\min\{1,\,T,\,|K+2|^{-1},\,z_0/(32C_\theta+1)\}
\end{equation}
(the term $|K+2|^{-1}$ omitted when $K=-2$); this is the constant $\tau_0$ of
Theorem~\ref{thm:A}, used throughout Section~\ref{sec:mono}:
\[
F^i_s:=\delta_{\lambda_i}\circ\Upsilon_p\circ F_{T+\lambda_i^{-2}s}\quad\text{on
}\Sigma^i_s:=F_{T+\lambda_i^{-2}s}^{-1}(U_p),\qquad
\theta^i_s:=\theta(\cdot,T+\lambda_i^{-2}s)\quad\text{on }\Sigma^i_s,
\]
so that $F^i_s$ is a Legendrian immersion of the open subset $\Sigma^i_s\subset\Sigma$,
which depends on $s$, into the rescaled box
$\delta_{\lambda_i}\Upsilon_p(U_p)=\{|w|<\lambda_ir_0,\ |z|<\lambda_i^2z_0\}$, with the
rescaled structures $\sigma_i:=\sigma_{\lambda_i}=(\eta_i,\xi_i,\Phi_i,g_i)$ of
Lemma~\ref{lem:moving}(b): the $\cD$-homothetic structure with $\kappa=\lambda_i^2$,
transported by $\delta_{\lambda_i}\Upsilon_p$. Thus
$\eta_i=\lambda_i^2(\delta_{\lambda_i}\Upsilon_p)_*\eta$, $\xi_i=2\partial_z$ and
$g_i=\lambda_i^2(\delta_{\lambda_i}\Upsilon_p)_*g^T+\eta_i\otimes\eta_i$. They converge to
$(\overline\eta,\overline\xi,\overline\Phi,\overline g)$ in $C^\infty_{\rm loc}(\Hn)$ by
Lemma~\ref{lem:moving}(b). In the notation of \cite[\S4.1]{ChangWuZhang2026ContactBlowup},
$F^i=F^{\lambda_i}$, written in the chart $\Upsilon_p$, which at the single point $p$ is a
chart as in \cite[Lemma~4.1]{ChangWuZhang2026ContactBlowup}. The Legendrian angle is not
rescaled, $\theta_\kappa=\theta$ \cite[\S2.3]{ChangWuZhang2026ContactBlowup}.

The flow $F^i$ is not a flow of a closed manifold, but it satisfies the local equations of
(F) for the structure $\sigma_i$, with Legendrian angle $\theta^i$ and with
$(K+2)\lambda_i^{-2}$ in place of $K+2$. Each $F^i_s$ is a normally parametrised
Legendrian immersion of $\Sigma^i_s$ into $(\delta_{\lambda_i}\Upsilon_p(U_p),\sigma_i)$,
with mean curvature vector $H^i=\Phi_i\nabla^i\theta^i$,
$\partial_sF^i=H^i-2\theta^i\xi_i$ and
$(\partial_s-\Delta^i)\theta^i=(K+2)\lambda_i^{-2}\theta^i$. Under (I), moreover,
$|B^i|^2\le C_I/|s|$. All statements about $F^i$ below are local in $\Hn$: integrals over
the rescaled flow are taken over $\Sigma^i_s$ or over the indicated subsets of it. Test
functions have compact support in $\Hn$; for $i$ large this support lies inside the
rescaled box, since the boxes exhaust $\Hn$.

\section{Height-localised almost monotonicity}\label{sec:mono}

This section is devoted to the almost monotonicity of the weighted kernel. A sheet on
which $\theta$ is constant moves vertically with speed $-4\theta$ in the coordinate $z$,
and the caloric height $u=z-4\tau\theta$ is the height corrected for this drift; along a
solution of \eqref{eq:LMCF} in the model, $u$ solves the heat equation (for the Lagrangian
counterpart, see Section~\ref{sec:weak}). The height
variable $v=u/\tau$ is dimensionless for the parabolic scaling $(w,z,\tau)\mapsto(\lambda
w,\lambda^2z,\lambda^2\tau)$. The weight enters through the inequality $\chi+16\chi''\ge
c_\chi\chi$ with $c_\chi=1-32a>0$, valid for $a<\frac1{32}$. It turns the cross terms
between the horizontal and the vertical direction into a non-negative defect. Hence in the
model the kernel is exactly monotone (Proposition~\ref{prop:model}), and in a curved
manifold the almost monotonicity uses the positive constant $c_\chi$ to absorb its error
terms (Proposition~\ref{prop:mono}).

The only term of the error which is not controlled by the area bounds, the one coming from
the term $(K+2)\theta$ in the evolution of the Legendrian angle, is controlled by a second
almost monotone quantity, $\int\theta^2\rho$, whose dissipation is $\int|H|^2\rho$
(Proposition~\ref{prop:theta2}). This is the Legendrian form of the argument by which
Neves \cite[Lemma~5.4]{nev07} obtains stationary limits for the Lagrangian mean curvature
flow with bounded Lagrangian angle, and in Section~\ref{sec:weak} it makes the limits
stationary.

\emph{Notation.} Write $\tau=T-t$, and put
\[
a_K(\tau)=\frac4{K+2}\big(e^{(K+2)\tau}-1\big)\quad(=4\tau\ \text{if }K+2=0),\qquad
\eps_K(\tau)=\frac{a_K(\tau)-4\tau}{4\tau}.
\]
Then $a_K'=4+(K+2)a_K$, and for $|K+2|\tau\le1$ one has $a_K\le8\tau$ and
$|\eps_K|\le2|K+2|\tau$. The weight is the Gaussian $\chi(v)=e^{-av^2}$ with
$0<a<\frac1{32}$, so that
\begin{equation}\label{eq:chi}
\chi+16\chi''=(1-32a+64a^2v^2)\chi\ \ge\ c_\chi\,\chi,\qquad c_\chi:=1-32a>0 .
\end{equation}
Moreover $\sup|\chi'|$, $\sup|\chi''|$ and $\sup\chi''^2/\chi$ are finite,
$\chi'^2/\chi=4a^2v^2e^{-av^2}\le4a/e$, and $v\chi'(v)<0$ for $v\ne0$. Unless stated
otherwise, a constant $C$ depends only on $n$, $a$, $K$, $C_\theta$, $D_0$,
$\mu_0(\Sigma)$ and $(M,\eta,g)$.

\subsection{The kernel in the model}\label{subsec:model}
In the Heisenberg group with its structure
$(\overline\eta,\overline\xi,\overline\Phi,\overline g)$ the horizontal potential is
$f=|w|^2/16$, and we write $\rho=\rho_{0,T}:=(4\pi\tau)^{-n/2}e^{-f/\tau}$ for the
transverse Gaussian kernel centred at $(0,T)$. The model structure is $\eta$-Einstein with
$K=-2$, since its transverse metric $g_0$ is flat; so hypothesis (F) for a flow in the
model reads
\begin{equation}\label{eq:flow-model}
H=\overline\Phi\nabla\theta,\qquad \partial_tF=H-2\theta\overline\xi,\qquad
(\partial_t-\Delta)\theta=0 ,
\end{equation}
and $a_K(\tau)=4\tau$. The \emph{caloric height} and the \emph{normalised height} of a
flow \eqref{eq:flow-model} are
\[
u:=z-4\tau\theta,\qquad v:=\frac u\tau ,
\]
and the \emph{contact Gaussian kernel} centred at $(0,T)$ is $\rho^\chi:=\rho\,\chi(v)$.
For a centre $\hat y\in\Hn$ one replaces $f$ and $z$ by $f_{\hat y}=f\circ\ell_{\hat
y}^{-1}$ and $z_{\hat y}=z\circ\ell_{\hat y}^{-1}$. We also write $\cS:=H+(\overline\nabla
f)^\perp/\tau$ for the defect of self-similarity.

\begin{lemma}\label{lem:model-exact}
Let $S\subset\overline{\cD}$ be a Legendrian plane at a point of $\Hn$, and for a function
$\gamma$ let $\nabla^S\gamma$ be the orthogonal projection of $\overline\nabla\gamma$ onto
$S$.
\begin{enumerate}[(a)]
\item $\overline\nabla f=\frac12X^{\cH}$ is horizontal, $|\overline\nabla f|^2=f$,
$\overline\nabla^2f(U,W)=\frac12\overline g(U,W)$ for horizontal $U,W$, and
$\tr_S\overline\nabla^2f=\frac n2$.
\item $\overline\nabla^2z(U,W)=0$ for horizontal $U,W$; in particular
$\tr_S\overline\nabla^2z=0$.
\item $\overline\nabla z=4\overline\Phi\,\overline\nabla f+2\overline\xi$. Hence
$\nabla^Sz=4\overline\Phi\big((\overline\nabla f)^{\perp_S}\big)$ and
$|\nabla^Sz|^2=16|(\overline\nabla f)^{\perp_S}|^2$, where $\perp_S$ is the component
orthogonal to $S$.
\item Along a flow \eqref{eq:flow-model},
\[
(\partial_t+\Delta-|H|^2)\rho=-|\cS|^2\rho,\qquad (\partial_t-\Delta)z=-4\theta,\qquad
(\partial_t-\Delta)u=0,\qquad \nabla u=4\tau\,\overline\Phi\cS ,
\]
so that $|\nabla v|^2=16|\cS|^2$.
\end{enumerate}
\end{lemma}
\begin{proof}
(a), (b) and the first identity of (c) are the model cases of
Lemma~\ref{lem:model}(a)--(d), where the errors $e_H$, $e_\rho$, $e_z$ vanish; they are
the identities \cite[(2.6)]{ChangWuZhang2026ContactBlowup} and those in the proofs of
\cite[Lemmas~2.4 and~2.5]{ChangWuZhang2026ContactBlowup}. For a Legendrian plane
$\overline\Phi S=S^\perp\cap\overline{\cD}$, and $\overline\Phi$ is an isometry of
$\overline{\cD}$, so $P_S\overline\Phi=\overline\Phi P_{S^\perp}$ on $\overline{\cD}$ and
the remaining statements of (c) follow. For (d), the first two identities are the model
case of Lemma~\ref{lem:model}(e). Since
$(\partial_t-\Delta)(4\tau\theta)=-4\theta+4\tau(\partial_t-\Delta)\theta=-4\theta$, the
function $u$ is caloric. Finally $\nabla u=\nabla^Sz-4\tau\nabla\theta$ with $S$ the
tangent plane, and $-\nabla\theta=\overline\Phi H$ because $H=\overline\Phi\nabla\theta$
and $\overline\Phi^2=-I$ on $\overline{\cD}$; so $\nabla u=4\overline\Phi(\overline\nabla
f)^\perp+4\tau\overline\Phi H=4\tau\overline\Phi\cS$.
\end{proof}

In the model the kernel is exactly monotone, as the backward heat kernel is along the mean
curvature flow in Euclidean space \cite{hu90}. The following statement is proved at the
end of Section~\ref{subsec:almost}, as the error-free case of the computation there.

\begin{proposition}\label{prop:model}
Let $\Sigma$ be an $n$-manifold without boundary and let
$F\colon\Sigma\times[t_0,T)\to\Hn$ be a smooth family of Legendrian immersions satisfying
\eqref{eq:flow-model}, such that $(x,t)\mapsto(F(x,t),t)$ is a proper map on
$\Sigma\times[t_0,T']$ for every $T'<T$, with $\sup|\theta|<\infty$ and
$\mu_t\big(F_t^{-1}\{|w|<R\}\big)\le D_0R^n$ for all $R>0$ and $t\in[t_0,T)$. Then
$\int\rho^\chi\,d\mu_t<\infty$, and for $t_0\le t_1<t_2<T$
\[
\int\rho^\chi\,d\mu_{t_2}-\int\rho^\chi\,d\mu_{t_1}=-\int_{t_1}^{t_2}\!\!\int
D^\chi\rho\,d\mu_t\,dt,\qquad
D^\chi:=\frac{|v\chi'(v)|}\tau+(\chi+16\chi'')(v)\,|\cS|^2\ \ge0 .
\]
In particular $t\mapsto\int\rho^\chi\,d\mu_t$ is non-increasing. It is constant on
$[t_1,t_2]$ if and only if $z=4\tau\theta$ and $H=-(\overline\nabla f)^\perp/\tau$ on
$\Sigma_t$ for $t\in[t_1,t_2]$, that is, the flow is an exactly centred self-shrinker
there. The same holds for every centre $(\hat y,T)$, $\hat y\in\Hn$, that is, for
$\ell_{\hat y}^{-1}\circ F$ in place of $F$.
\end{proposition}

\begin{example}\label{ex:two-sheets}
Let $L\subset\C^n$ be a Lagrangian subspace and $c\ne0$. The horizontal planes
$P_1=L\times\{0\}$ and $P_2=L\times\{c\}$ are Legendrian, since
$\overline\beta=\frac12\iota_{\fr}\omega_0$ vanishes on $L$, and minimal; with
$\theta\equiv0$ they are static solutions of \eqref{eq:flow-model}. For the centre $(0,T)$
and $\tau>0$, $v=0$ on $P_1$ and $v=c/\tau$ on $P_2$, and the Gaussian integral of $\rho$
over each of $P_1$ and $P_2$ is $1$: on both planes $\overline g$ restricts to $g_0$, both
project onto $L$, and for every linear $n$-dimensional subspace $P\subset\C^n$,
$\int_P(4\pi\tau)^{-n/2}e^{-|w|^2/(16\tau)}\,d\cH^n_{g_0}
=(4\pi\tau)^{-n/2}\,2^{-n}(16\pi\tau)^{n/2}=1$. Hence
\[
\int_{P_1\cup P_2}\rho^\chi\,d\cH^n=1+e^{-ac^2/\tau^2}\longrightarrow1\quad(\tau\to0),
\]
whereas the height-blind integral $\int_{P_1\cup P_2}\rho\,d\cH^n$ equals $2$ for every
$\tau$. If instead the two planes carry constant Legendrian angles
$\bar\theta_1,\bar\theta_2$, they move vertically with the speeds $-4\bar\theta_j$.
If the kernel is centred at a point of the first plane at the time $t_*$, the second plane
is at the height $c+4\bar\theta_2\sigma$ at the time $t_*-\sigma$, so that $v=c/\sigma$ on
it. The drift of each sheet is compensated by its own Legendrian angle, and the weighted
integral is again $1+e^{-ac^2/\sigma^2}$. Thus the weight separates two sheets only at
scales $\sigma$ small compared with their height difference, and at such scales it sees
only the sheet through the centre. On $P_1$ the defect vanishes, $v=0$ and $\cS=0$; on
$P_2$ the term $|v\chi'(v)|/\tau$ is positive, in accordance with
Proposition~\ref{prop:model}.
\end{example}

\subsection{Local almost monotonicity}\label{subsec:almost}
Fix $p\in M$, and let $\zeta:=\zeta_p$ and $\psi:=\psi_p=\zeta^2$ be the cut-offs defined
after Lemma~\ref{lem:charts}; so $\zeta=1$ on
$U_p^{1/2}:=\{|w|<\frac{r_0}2,|z_p|<\frac{z_0}2\}$. The statements of this section are
written for the centre $(p,T)$. They hold verbatim for any centre $(q,t_*)$, $q\in M$,
$0<t_*\le T$, with $\psi_q$ in place of $\psi$, with $\tau=t_*-t$, on the interval
$[t_*-\min\{\tau_0,t_*\},t_*)$, and with the same constants. The lemmas of
Section~\ref{subsec:density} are stated for general centres from the outset, and
Lemma~\ref{lem:tails} is proved directly from (D), without restriction on the centre time.
Indeed, the proofs use only (S), (F), the property (D) of Lemma~\ref{lem:D} on that
interval and the error bounds of Lemma~\ref{lem:model}, whose constants do not depend on
the centre. They use the cut-off only through bounds for $|\overline\nabla\zeta_q|$ and
$|\overline\nabla^2\zeta_q|$ and through the sets where $\zeta_q=1$ and where
$\zeta_q\ne0$, which are the preimages under $\Upsilon_q$ of sets fixed by $\zeta_*$; the
bounds and these coordinate sets do not depend on $q$.

\begin{proposition}\label{prop:theta2}
Assume (S), (F). There is a constant $C_H$ such that for $T-\min\{1,T\}\le t<T$
\[
\frac{d}{dt}\int\psi\,\theta^2\rho\,d\mu\le-\int\psi|H|^2\rho\,d\mu+C_H
,\qquad\text{hence}\qquad \int_{T-\min\{1,T\}}^T\!\!\int\psi|H|^2\rho\,d\mu\,dt<\infty .
\]
In particular both statements hold on $[T-\tau_0,T)$ for the constant
$\tau_0\le\min\{1,T\}$ of \eqref{eq:tau0}, which is the time interval used from
Proposition~\ref{prop:mono} onwards.
\end{proposition}
\begin{proof}
For $T-\min\{1,T\}\le t<T$ we have $\tau\le1$, so that Lemma~\ref{lem:model}(f) applies.
On the closed manifold $\Sigma_t$, integration by parts gives
\[
\frac{d}{dt}\int\psi\theta^2\rho=\int(\partial_t-\Delta)(\psi\theta^2)\,\rho+\int\psi\theta^2(\partial_t+\Delta-|H|^2)\rho
.
\]
Here $(\partial_t-\Delta)\theta^2=2(K+2)\theta^2-2|\nabla\theta|^2=2(K+2)\theta^2-2|H|^2$,
since $\Phi$ is an isometry on $\cD$, and
$(\partial_t-\Delta)\psi=-2\theta\,\xi\psi-\tr_\Sigma\overline\nabla^2\psi$ is bounded,
the
terms $d\psi(H)$ cancelling. The cross term from $(\partial_t-\Delta)(\psi\theta^2)$ is
absorbed by half of the good term:
$4|\theta||\nabla\psi||\nabla\theta|\le\psi|H|^2+16\theta^2|\nabla\zeta|^2$. With
Lemma~\ref{lem:model}(e)--(f), the right-hand side is at most
$-2\int\psi|H|^2\rho+\int\psi|H|^2\rho+C\int_{F^{-1}(U_p)}(1+|e_\rho|)\rho\le-\int\psi|H|^2\rho+C$.
Integrating, and using $\int\psi\theta^2\rho\ge0$, gives the second assertion.
\end{proof}

\begin{proposition}\label{prop:mono}
Assume (S), (F). Put
\begin{gather*}
u:=z_p-a_K(\tau)\theta,\qquad v:=u/\tau,\qquad \rho^\chi=\rho^\chi_{p,T}:=\rho\,\chi(v),\\
D^\chi:=\frac{|v\chi'(v)|}{\tau}+(\chi+16\chi'')(v)\,|\mathcal
S_p|^2\ \ge\ \frac{|v\chi'(v)|}\tau+c_\chi\chi(v)|\mathcal S_p|^2 .
\end{gather*}
Then there are $\tau_0>0$ and $C_2$ such that for $T-\tau_0\le t<T$
\begin{equation}\label{eq:mono}
\frac{d}{dt}\int\psi\rho^\chi\,d\mu\le-\frac34\int\psi D^\chi\rho\,d\mu+E(t),\qquad
E(t):=C_2\Big(1+\tau^2\int\psi|H|^2\rho\,d\mu\Big).
\end{equation}
By Proposition~\ref{prop:theta2}, $E$ is integrable on $[T-\tau_0,T)$. Hence the function
\[
t\longmapsto\int\psi\rho^\chi\,d\mu_t-\int_{T-\tau_0}^tE
\]
is non-increasing on $[T-\tau_0,T)$ and bounded below by $-\int_{T-\tau_0}^TE$.
Consequently $\lim_{t\to T}\int\psi\rho^\chi\,d\mu_t$ exists and
$\int_{T-\tau_0}^T\!\int\psi D^\chi\rho<\infty$. The height-blind version ($\chi\equiv1$,
$D=|\mathcal S_p|^2$) holds in the same way:
$\frac{d}{dt}\int\psi\rho\le-\int\psi|\mathcal S_p|^2\rho+C$.
\end{proposition}
\begin{proof}
Let $\tau_0$ be the constant \eqref{eq:tau0}. The two constraints beyond
$\tau_0\le\min\{1,T\}$ are used below in the bounds $a_K(\tau)\le8\tau$,
$|\eps_K|\le2|K+2|\tau$ and in the estimate of $R_4$.

Since $(\partial_t-\Delta)(a_K\theta)=-a_K'\theta+(K+2)a_K\theta=-4\theta$,
Lemma~\ref{lem:model}(e) gives $(\partial_t-\Delta)u=e_z$, hence
$(\partial_t-\Delta)v=v/\tau+e_z/\tau$. For the gradient, Lemma~\ref{lem:model}(d) and
$-\nabla\theta=\Phi H$ give
\begin{equation}\label{eq:grad-u}
\nabla u=4\Phi(\overline\nabla f_p)^\perp+P_\Sigma e_H+a_K\Phi H=4\tau\,\Phi(\mathcal
S_p+e_u),\qquad e_u:=\eps_KH-\Phi P_\Sigma e_H/(4\tau),
\end{equation}
so $|\nabla v|^2=16|\mathcal S_p+e_u|^2$ and $|e_u|^2\le
C((K+2)^2\tau^2|H|^2+|w|^6/\tau^2)$. By the chain rule,
\[
\frac{d}{dt}\int\psi\chi\rho=\int\chi(\partial_t-\Delta)\psi\,\rho+\int\psi(\partial_t-\Delta)\chi(v)\,\rho-2\int\chi'(v)\langle\nabla\psi,\nabla
v\rangle\rho+\int\psi\chi(-|\mathcal S_p|^2+e_\rho)\rho,
\]
with
\[
(\partial_t-\Delta)\chi(v)=\chi'(v)\Big(\frac
v\tau+\frac{e_z}\tau\Big)-16\chi''(v)|\mathcal S_p+e_u|^2 .
\]
The principal terms combine to $-\int\psi D^\chi\rho$, because $\chi'(v)v\le0$. Among
the remainders, $R_1$, $R_2$ and $R_3$ come from the curvature errors $e_z$, $e_u$ and
$e_\rho$, and are bounded with Lemma~\ref{lem:model}(e)--(f) and \eqref{eq:eu-bound},
while $R_4$ and $R_5$ come from the
cut-off and are controlled by the decay of the kernel and of the weight on its support.
The
remainders are:
\begin{itemize}
\item $R_1=\int\psi\chi'e_z\tau^{-1}\rho$, with $|R_1|\le C\int(|w|^2/\tau)\rho\le C$ by
Lemma~\ref{lem:model}(f).
\item $R_2=-16\int\psi\chi''(2\langle\mathcal
S_p,e_u\rangle+|e_u|^2)\rho\le\frac{c_\chi}8\int\psi\chi|\mathcal
S_p|^2\rho+C\int\psi|e_u|^2\rho$, by $\sup\chi''^2/\chi<\infty$.
\item $R_3=\int\psi\chi e_\rho\rho$, with $|R_3|\le C$ by Lemma~\ref{lem:model}(e)--(f).
\item $R_4=\int\chi(\partial_t-\Delta)\psi\,\rho$. Its integrand vanishes outside
$\spt\overline\nabla\psi$, and there either $|w|\ge r_0/2$, so that
$\rho\le(4\pi\tau)^{-n/2}e^{-r_0^2/(64\tau)}$; or $|z_p|\ge z_0/2$, so that
$|v|\ge\frac{z_0}{2\tau}-8C_\theta\ge\frac{z_0}{4\tau}$ and $\chi(v)\le
e^{-az_0^2/(16\tau^2)}$. Hence $|R_4|\le C\mu_0(\Sigma)$.
\item $R_5=-2\int\chi'\langle\nabla\psi,\nabla v\rangle\rho$. Since
$|\nabla\psi|^2/\psi=4|\nabla\zeta|^2$, we get
$|R_5|\le\frac{c_\chi}8\int\psi\chi(|\mathcal
S_p|^2+|e_u|^2)\rho+C\int_{\spt\overline\nabla\psi}\frac{\chi'^2}{\chi}\rho$. The last
integral is bounded as for $R_4$, since $\chi'^2/\chi=4a^2v^2e^{-av^2}$.
\end{itemize}
Collecting terms, and using $c_\chi\chi|\mathcal S_p|^2\le D^\chi$ together with
\begin{equation}\label{eq:eu-bound}
\int\psi|e_u|^2\rho\le C(K+2)^2\tau^2\int\psi|H|^2\rho+C\tau\int(|w|^2/\tau)^3\rho\le
C\tau^2\int\psi|H|^2\rho+C\tau D_0,
\end{equation}
the last step by Lemma~\ref{lem:model}(f), gives \eqref{eq:mono}. Thus
$\frac{d}{dt}\big(\int\psi\rho^\chi-\int_{T-\tau_0}^tE\big)\le-\frac34\int\psi
D^\chi\rho\le0$, and the quantity in brackets is bounded below because
$\int\psi\rho^\chi\ge0$ and $\int_{T-\tau_0}^TE<\infty$. A non-increasing function bounded
below has a limit, and integrating the inequality gives the bound on $\int\!\int\psi
D^\chi\rho$. For $\chi\equiv1$ the terms $R_1,R_2,R_5$ vanish and $|R_4|\le
C\int_{F^{-1}(U_p)}\rho\le CD_0$, which gives the height-blind version.
\end{proof}

\begin{proof}[Proof of Proposition~\ref{prop:model}]
Fix $t_0\le t_1<t_2<T$ and put $\tau_2:=T-t_2$, $\tau_1:=T-t_1$, so that
$\tau\in[\tau_2,\tau_1]$ on $[t_1,t_2]$. The finiteness of
$\int\rho^\chi\,d\mu_t\le\int\rho\,d\mu_t$ follows from the area bound and the layer-cake
formula, as in Lemma~\ref{lem:model}(f). The cut-off identity below is used twice: first
to show that the dissipation is integrable, and then, with this integrability, to show
that
the cut-off remainder $R_5$ disappears in the limit, which gives the identity.

Fix $\zeta\in C^\infty_c(\R^N)$, radially non-increasing, with $0\le\zeta\le1$, $\zeta=1$
on the ball of radius $\frac12$ and $\zeta=0$ outside the unit ball, for $N=2n$ and $N=1$,
and for $R,Z>0$ put $\psi=\psi_{R,Z}(w,z):=\zeta(w/R)^2\zeta(z/Z)^2$. Since
$(x,t)\mapsto(F(x,t),t)$ is proper on $\Sigma\times[t_1,t_2]$, the set of $x\in\Sigma$
with $F(x,t)\in\spt\psi$ for some $t\in[t_1,t_2]$ is relatively compact. Hence $\psi\circ
F_t$ is supported in a fixed compact subset of $\Sigma$ for $t\in[t_1,t_2]$, and the
differentiation under the integral sign below is justified. Moreover, the expansion in the
proof of Proposition~\ref{prop:mono} applies to $\psi$ with $K+2=0$, $a_K=4\tau$ and
$e_H=e_\rho=e_z=e_u=0$ (Lemma~\ref{lem:model-exact}). The remainders $R_1$, $R_2$, $R_3$
there vanish identically, so that
\begin{equation}\label{eq:model-cutoff}
\frac{d}{dt}\int\psi\,\rho^\chi\,d\mu_t=-\int\psi\,D^\chi\rho\,d\mu_t+R_4+R_5,\qquad
R_4=\int\chi\,(\partial_t-\Delta)\psi\,\rho,\quad R_5=-2\int\chi'\langle\nabla\psi,\nabla
v\rangle\rho .
\end{equation}
The support of $\overline\nabla\psi$ lies in $A_1\cup A_2$, where $A_1=\{\frac R2\le|w|\le
R\}$ and $A_2=\{|w|\le R,\ \frac Z2\le|z|\le Z\}$. By Lemma~\ref{lem:model-exact}(c),
$|\overline\nabla z|\le C(1+|w|)$, so on $A_1\cup A_2$
\[
\frac{|\nabla\psi|^2}{\psi}\le C\Big(\frac1{R^2}+\frac{1+R^2}{Z^2}\Big)=:C\,\kappa_{R,Z} ,
\]
and, computing
$(\partial_t-\Delta)\psi=-2\theta\,\overline\xi\psi-\tr_\Sigma\overline\nabla^2\psi$ with
the connection formula in the proof of Lemma~\ref{lem:model},
$|(\partial_t-\Delta)\psi|\le C(\kappa_{R,Z}+Z^{-1})$ there. On $A_1$ we use $\chi\le1$
and $\chi'^2/\chi\le4a/e$; on $A_2$ we have $|v|\ge Z/(2\tau)-4C_\theta$, so that $\chi$
and $\chi'^2/\chi$ are at most $C\exp\big(-\tfrac a2(Z/(2\tau)-4C_\theta)_+^2\big)$, while
$\rho\le(4\pi\tau)^{-n/2}$ and $\mu_t(F_t^{-1}(A_2))\le D_0R^n$. By the layer-cake formula
and the area bound, $\int_{F_t^{-1}(A_1)}\rho\,d\mu_t\le C(n)D_0e^{-R^2/(128\tau)}$.
Hence, with
\[
\varepsilon_{R,Z}:=C(\kappa_{R,Z}+Z^{-1})\Big[C(n)D_0e^{-R^2/(128\tau_1)}+(4\pi\tau_2)^{-n/2}D_0R^n\exp\big(-\tfrac
a2(Z/(2\tau_1)-4C_\theta)_+^2\big)\Big],
\]
which tends to $0$ as first $Z\to\infty$ and then $R\to\infty$, we have on $[t_1,t_2]$
\[
|R_4|\le\varepsilon_{R,Z},\qquad \int_{F_t^{-1}(A_1\cup
A_2)}\frac{\chi'^2}{\chi}\frac{|\nabla\psi|^2}{\psi}\rho\,d\mu_t\le\varepsilon_{R,Z} .
\]
Since $|\nabla v|=4|\cS|$, the Cauchy--Schwarz inequality gives
\[
|R_5|\le8\Big(\int_{A_1\cup
A_2}\frac{\chi'^2}{\chi}\frac{|\nabla\psi|^2}{\psi}\rho\Big)^{1/2}\Big(\int\psi\chi|\cS|^2\rho\Big)^{1/2},
\]
and by Young's inequality
$|R_5|\le\frac{c_\chi}8\int\psi\chi|\cS|^2\rho+C\varepsilon_{R,Z}$.

Integrating \eqref{eq:model-cutoff} over $[t_1,t_2]$ and using Young's form of the bound
for $R_5$ together with $c_\chi\chi|\cS|^2\le D^\chi$, we obtain
\[
\int\psi\rho^\chi\,d\mu_{t_2}-\int\psi\rho^\chi\,d\mu_{t_1}\le-\frac78\int_{t_1}^{t_2}\!\!\int\psi
D^\chi\rho\,d\mu_t\,dt+C(t_2-t_1)\varepsilon_{R,Z} .
\]
Let $Z\to\infty$ and then $R\to\infty$: $\psi\uparrow1$, the left-hand side converges to
$\int\rho^\chi d\mu_{t_2}-\int\rho^\chi d\mu_{t_1}$ by monotone convergence, and
$\int\!\!\int\psi D^\chi\rho\uparrow\int\!\!\int D^\chi\rho$. Since the left-hand side is
bounded below by $-C(n)D_0$, this shows $\int_{t_1}^{t_2}\!\!\int
D^\chi\rho\,d\mu_t\,dt<\infty$, and in particular
\begin{equation}\label{eq:model-finite}
\int_{t_1}^{t_2}\!\!\int\chi|\cS|^2\rho\,d\mu_t\,dt<\infty .
\end{equation}
For the identity itself we integrate \eqref{eq:model-cutoff} once more, now using the
Cauchy--Schwarz form of the bound for $R_5$ and \eqref{eq:model-finite}:
\[
\int_{t_1}^{t_2}|R_4|+|R_5|\,dt\le(t_2-t_1)\varepsilon_{R,Z}+8\Big((t_2-t_1)\varepsilon_{R,Z}\Big)^{1/2}\Big(\int_{t_1}^{t_2}\!\!\int\chi|\cS|^2\rho\Big)^{1/2}\longrightarrow0
.
\]
Letting $Z\to\infty$ and $R\to\infty$ as before gives the identity. It shows that
$t\mapsto\int\rho^\chi d\mu_t$ is non-increasing, and that it is constant on $[t_1,t_2]$
if and only if $D^\chi=0$ on $\Sigma_t$ for $t\in[t_1,t_2]$. By \eqref{eq:chi},
$(\chi+16\chi'')(v)>0$, and $v\chi'(v)=0$ only at $v=0$; so $D^\chi=0$ means $v=0$ and
$\cS=0$, that is, $z=4\tau\theta$ and $H=-(\overline\nabla f)^\perp/\tau$.

For an arbitrary centre, let $\hat y\in\Hn$ have horizontal coordinate $\hat w$. The left
translation $\ell_{\hat y}$ preserves $(\overline\eta,\overline\xi,\overline\Phi,\overline
g)$, so $G:=\ell_{\hat y}^{-1}\circ F$ satisfies \eqref{eq:flow-model} with the same
$\theta$, and $(x,t)\mapsto(G(x,t),t)$ is proper on $\Sigma\times[t_0,T']$ for $T'<T$.
Since $\ell_{\hat y}^{-1}$ maps $\{|w-\hat w|<R\}$ onto $\{|w|<R\}$, the area bound for
$G$ takes the form $\mu_t\big(G_t^{-1}\{|w|<R\}\big)\le\mu_t\big(F_t^{-1}\{|w|<R+|\hat
w|\}\big)\le D_0(R+|\hat w|)^n$. The argument above uses the area bound only for the bound
$\mu_t(F_t^{-1}(A_2))\le D_0R^n$ with $R$ large, and in the layer-cake estimates for
$\int\rho\,d\mu_t$ and $\int_{F_t^{-1}(A_1)}\rho\,d\mu_t$. With $D_0(R+|\hat w|)^n$ in
place of $D_0R^n$ these estimates hold with constants depending also on $|\hat w|$ and
$\tau_2$, the quantity corresponding to $\varepsilon_{R,Z}$ still tends to $0$ as first
$Z\to\infty$ and then $R\to\infty$, and the argument above applies to $G$.
\end{proof}

\subsection{Contact densities, tails and clearing-out}\label{subsec:density}

In this section we show that the localised contact density at a point of the flow tends,
as the scale goes to $0$, to the number of sheets through that point. We then bound the
parts of the density far from the centre and deduce a clearing-out lemma.

\begin{lemma}\label{lem:density}
Assume (S), (F). For $0<t<T$ and $q\in M$ put
\[
\Theta^\chi_{q,t}(\sigma):=\int\psi_q\,\rho^\chi_{q,t}\,d\mu_{t-\sigma},\qquad
0<\sigma\le\min\{\tau_0,t\},
\]
the localised contact density, with the kernel of Proposition~\ref{prop:mono} centred at
$(q,t)$. Then
\[
\lim_{\sigma\to0}\Theta^\chi_{q,t}(\sigma)=N(q,t):=\#F_t^{-1}(q)\qquad\text{(counted with
multiplicity)},
\]
and $\Theta^\chi_{q,t}(\sigma)\ge N(q,t)-C_5\sigma$ for $0<\sigma\le\min\{\tau_0,t\}$,
with $C_5$ independent of $(q,t)$.
\end{lemma}
\begin{proof}
We work in the chart $\Upsilon_q$, write $w$ and $z=z_q$ for its coordinates, and fix a
compact set $Q\subset\{|w|<r_0,\ |z|<z_0\}$ whose interior contains $\spt\zeta_*$. The
flow
is smooth on $\Sigma\times[t-\sigma_0,t]$ with $\sigma_0:=\frac12\min\{\tau_0,t\}$, so
$\partial_tF$, $H$ and $\partial_t\theta$ are bounded there; the constants below may
depend
on these bounds, which is enough for the limit. We take $\sigma\le\sigma_0$.

A horizontal vector $U$ with $dw(U)=0$ is a multiple of $\partial_z$ with
$\eta(U)=\frac12dz(U)=0$, hence $U=0$. Thus $w\circ F_t$ is an immersion of
$F_t^{-1}(U_q)$, and the compact set $\{x\in F_t^{-1}(Q):\ w(F_t(x))=0\}$ is discrete,
hence finite. It consists of the $N=N(q,t)$ points $x_1,\dots,x_N$ with $F_t(x_k)=q$ and
of points $y_1,\dots,y_L$ with $\bar z_l:=z(F_t(y_l))\neq0$. Let $0<\eps<\frac14r_0$, and
choose pairwise disjoint open neighbourhoods $V_k$ of $x_k$ and $W_l$ of $y_l$ such that
$w\circ F_t$ is an embedding of each $V_k$ with $|w\circ F_t|\le\eps$ and $|z\circ
F_t|<\frac14z_0$ on $V_k$, and such that $|z\circ F_t-\bar z_l|<\frac14|\bar z_l|$ on
$W_l$. On the compact set $A:=F_t^{-1}(Q)\setminus\big(\bigcup_kV_k\cup\bigcup_lW_l\big)$
we have $|w\circ F_t|\ge2\delta$ for some $\delta>0$. Since $|\partial_tF|$ is bounded,
for
$\sigma$ small every $x$ with $F_{t-\sigma}(x)\in\spt\psi_q$ lies in $F_t^{-1}(Q)$;
moreover $|w\circ F_{t-\sigma}-w\circ F_t|<\delta$ and $|z\circ F_{t-\sigma}-z\circ
F_t|<\frac14\min_l|\bar z_l|$ on $F_t^{-1}(Q)$, and $\psi_q\circ F_{t-\sigma}=1$ on
$\bigcup_kV_k$.

The sheets that do not pass through $q$ do not contribute in the limit, for two different
reasons. On $A$ the horizontal distance stays positive: $|w\circ F_{t-\sigma}|\ge\delta$,
so $\psi_q\rho^\chi_{q,t}\le(4\pi\sigma)^{-n/2}e^{-\delta^2/(16\sigma)}$ there. On $W_l$
the horizontal coordinate may vanish, but the height does not: $|z\circ
F_{t-\sigma}|\ge\frac12|\bar z_l|$, and since $a_K(\sigma)\le8\sigma$ and $|\theta|\le
C_\theta$, $|v|\ge|\bar z_l|/(2\sigma)-8C_\theta$; so
$\psi_q\rho^\chi_{q,t}\le(4\pi\sigma)^{-n/2}\exp\big(-a(|\bar
z_l|/(2\sigma)-8C_\theta)_+^2\big)$ there. As $\mu_{t-\sigma}(\Sigma)\le\mu_0(\Sigma)$,
the contributions of $A$ and of the $W_l$ to $\Theta^\chi_{q,t}(\sigma)$ tend to $0$ as
$\sigma\to0$.

On a sheet through $q$ we first compute the integral without the weight. Fix $k$ and put
$G_\sigma:=w\circ F_{t-\sigma}$ on $V_k$. On horizontal vectors $g$ is $h_q$ in the
coordinates $w$, so $\mu_{t-\sigma}$ on $V_k$ is the $h_q$-area of the immersion
$G_\sigma$. As $F$ is smooth, $|G_\sigma-G_0|\le C\sigma$ and
$d\mu_{t-\sigma}=(1+O(\sigma))\,d\mu_t$ on $V_k$. Since $|G_0|\le\eps$ on $V_k$,
$\big||G_\sigma|^2-|G_0|^2\big|\le C\sigma(\eps+\sigma)$, and therefore
\[
e^{-C(\eps+\sigma)}\,e^{-|G_0|^2/(16\sigma)}\le e^{-|G_\sigma|^2/(16\sigma)}\le
e^{C(\eps+\sigma)}\,e^{-|G_0|^2/(16\sigma)}\qquad\text{on }V_k .
\]
The set $G_0(V_k)$ is an embedded $n$-dimensional submanifold of $\C^n$ through $0$, with
tangent plane $P_k$ at $0$; its $h_q$-area in $\{|w|<r\}$ is at most $Cr^n$, the dilations
$\lambda G_0(V_k)$ converge to $P_k$ smoothly on compact sets as $\lambda\to\infty$, and
$h_q=g_0+O(|w|^2)$. Substituting $w=\sqrt\sigma\,y$ and bounding the part $\{|y|>R\}$ by
$Ce^{-R^2/32}$, we obtain
\[
\lim_{\sigma\to0}\int_{V_k}(4\pi\sigma)^{-n/2}e^{-|G_0|^2/(16\sigma)}\,d\mu_t=\int_{P_k}(4\pi)^{-n/2}e^{-|y|^2/16}\,d\cH^n_{g_0}(y)=1
,
\]
the last integral being the one of Example~\ref{ex:two-sheets}. Hence every limit point of
$\int_{V_k}\psi_q\rho_{q,t}\,d\mu_{t-\sigma}$ as $\sigma\to0$ lies in
$[e^{-C\eps},e^{C\eps}]$.

It remains to show that the height weight does not change this limit. This is the point
where the contact density differs from the usual Gaussian density: the weight is close to
$1$ on a sheet through $q$ only because the vertical drift of the sheet is compensated by
its Legendrian angle.

\begin{claim}
For each $k$,
\[
\lim_{\sigma\to0}\Big(\int_{V_k}\psi_q\rho^\chi_{q,t}\,d\mu_{t-\sigma}-\int_{V_k}\psi_q\rho_{q,t}\,d\mu_{t-\sigma}\Big)=0 .
\]
\end{claim}

Note that $\chi(v)$ need not tend to $1$ uniformly on $V_k$; the claim is obtained by
splitting $V_k$ at a rescaled radius $R$. By Lemma~\ref{lem:model}(d) and the Legendrian
condition, $|\nabla(z\circ F_t)|\le4|(\overline\nabla f_q)^\perp|+C|w|^3\le C|w|^2$ on
$V_k$, the normal part of the radial field $\overline\nabla f_q$ vanishing to second order
along the embedded sheet through $0$. As $z(F_t(x_k))=0$, this gives $|z\circ F_t|\le
C|G_0|^3$ on $V_k$. Along the flow $\partial_t(z\circ F)=dz(H)-4\theta$ with $|dz(H)|\le
C|w||H|$ (Lemma~\ref{lem:model}(d)), and $|\theta(\cdot,t-s)-\theta(\cdot,t-\sigma)|\le
C\sigma$ for $0\le s\le\sigma$. Let $R\ge1$. On $V_k\cap\{|G_\sigma|\le R\sqrt\sigma\}$ we
have $|w\circ F_{t-s}|\le R\sqrt\sigma+C\sigma$ for $0\le s\le\sigma$, hence
\[
z\circ F_{t-\sigma}=4\sigma\theta(\cdot,t-\sigma)+O(R^3\sigma^{3/2}),\qquad v=\frac{z\circ
F_{t-\sigma}-a_K(\sigma)\theta}{\sigma}=O(R^3\sqrt\sigma)+O(|K+2|\sigma),
\]
using $a_K(\sigma)=4\sigma(1+\eps_K(\sigma))$ with $|\eps_K|\le2|K+2|\sigma$. Thus
$\chi(v)\to1$ uniformly on $V_k\cap\{|G_\sigma|\le R\sqrt\sigma\}$. On
$V_k\cap\{|G_\sigma|>
R\sqrt\sigma\}$ we use $0\le\chi\le1$ and Lemma~\ref{lem:model}(f) with $\ell=1$, which
gives $\int_{\{|w|>R\sqrt\sigma\}}\rho_{q,t}\,d\mu_{t-\sigma}\le C(n)D_0R^{-2}$. Hence
\[
\limsup_{\sigma\to0}\Big|\int_{V_k}\psi_q\rho^\chi_{q,t}\,d\mu_{t-\sigma}-\int_{V_k}\psi_q\rho_{q,t}\,d\mu_{t-\sigma}\Big|\le
C(n)D_0R^{-2}
\]
for every $R\ge1$, and the claim follows by letting $R\to\infty$.

Together with the preceding estimates, the claim shows that every limit point of
$\Theta^\chi_{q,t}(\sigma)$ as $\sigma\to0$ lies in $[Ne^{-C\eps},Ne^{C\eps}]$. Since
$\eps$ is arbitrary, $\Theta^\chi_{q,t}(\sigma)\to N(q,t)$.

\emph{Lower bound.} Proposition~\ref{prop:mono} at the centre $(q,t)$ gives
$\lim_{\sigma'\to0}\Theta^\chi_{q,t}(\sigma')\le\Theta^\chi_{q,t}(\sigma)+\int_{t-\sigma}^tE$.
Moreover $\int_{t-\sigma}^tE\le
C_2\sigma+C_2\sigma^2\int_{t-\sigma}^t\!\int\psi|H|^2\rho\le C_5\sigma$ with
$C_5:=C_2(1+C_\theta^2CD_0+C_H)$, because the last integral is bounded by
$C_\theta^2CD_0+C_H\sigma$ (Proposition~\ref{prop:theta2} at the same centre, and
Lemma~\ref{lem:model}(f)), and $\sigma\le1$. By the beginning of
Section~\ref{subsec:almost}, $C_H$ and $C_2$, hence $C_5$, do not depend on the centre;
the smoothness bounds used for the limit above are not needed here.
\end{proof}

\begin{lemma}\label{lem:tails}
Assume (S), (F). Let $q\in M$, $t>0$ and $0<\sigma\le\min\{\tau_0,t\}$ with $t-\sigma<T$;
the centre time $t$ may be $\ge T$, and
$\Theta^\chi_{q,t}(\sigma):=\int\psi_q\rho^\chi_{q,t}\,d\mu_{t-\sigma}$ as in
Lemma~\ref{lem:density}. For $R,Z>0$ the part of $\Theta^\chi_{q,t}(\sigma)$ coming from
$\{f_q\ge R^2\}\cup\{|z_q|\ge Z\}$ is at most
\[
\tail(R,Z;\sigma):=C(n)D_0\,e^{-R^2/(2\sigma)}+(4\pi\sigma)^{-n/2}D_0(4R)^n\exp\big(-a(Z/\sigma-8C_\theta)_+^2\big).
\]
For fixed $0<\sigma_1\le\sigma_2$ and $\eps>0$ one can first choose $R$ and then $Z=Z(R)$
with $\tail(R,Z;\sigma)\le\eps$ for all $\sigma\in[\sigma_1,\sigma_2]$.
\end{lemma}
\begin{proof}
Only the time slice $t-\sigma$ enters: the proof uses (D) at that time, $|\theta|\le
C_\theta$ and $a_K(\sigma)\le8\sigma$ (as $|K+2|\sigma\le1$), and not $t<T$. On $\{f_q\ge
R^2\}=\{|w|\ge4R\}$ bound $\chi\le1$ and use the layer-cake formula with $m(r)\le D_0r^n$
($r\le r_0$; the cut-off $\psi_q$ is supported in the chart):
\[
\int_{\{|w|\ge4R\}}(4\pi\sigma)^{-n/2}e^{-|w|^2/(16\sigma)}d\mu\le
D_0(4\pi)^{-n/2}\,2\cdot4^n\int_{R/\sqrt\sigma}^\infty u^{n+1}e^{-u^2}\,du\le
C(n)D_0e^{-R^2/(2\sigma)} .
\]
On $\{f_q<R^2,\ |z_q|\ge Z\}$ we have $|v|\ge Z/\sigma-a_K(\sigma)|\theta|/\sigma\ge
Z/\sigma-8C_\theta$, the kernel is at most
$(4\pi\sigma)^{-n/2}e^{-a(Z/\sigma-8C_\theta)_+^2}$, and the mass is at most $D_0(4R)^n$.
\end{proof}

\begin{lemma}\label{lem:clearing}
Assume (S), (F). There are $c_{\rm cl},R_*,Z_*>0$, depending only on $n,a,C_\theta,D_0$,
with the following property. Let $0<t<T$, $q\in M$ and
$0<\sigma\le\min\{\tau_0,t,1/(4C_5)\}$, and put $\mathcal
Q_\sigma(q):=\Upsilon_q^{-1}\{f_q<R_*^2\sigma,\ |z_q|<Z_*\sigma\}$. If
\[
\mu_{t-\sigma}\big(F_{t-\sigma}^{-1}(\mathcal Q_\sigma(q))\big)\le c_{\rm
cl}\,\sigma^{n/2},
\]
then $q$ is not in the image of $F_t$.
\end{lemma}
\begin{proof}
Suppose $q=F_t(x)$. By Lemma~\ref{lem:density},
$\Theta^\chi_{q,t}(\sigma)\ge1-C_5\sigma\ge\frac34$. On the other hand
$\rho^\chi_{q,t}\le(4\pi\sigma)^{-n/2}$, so the part of $\Theta^\chi_{q,t}(\sigma)$ coming
from $\mathcal Q_\sigma(q)$ is at most $(4\pi)^{-n/2}c_{\rm cl}$. By Lemma~\ref{lem:tails}
with $R=R_*\sqrt\sigma$ and $Z=Z_*\sigma$, the remaining part is at most
\[
C(n)D_0e^{-R_*^2/2}+(4\pi)^{-n/2}D_0(4R_*)^n\exp\big(-a(Z_*-8C_\theta)_+^2\big),
\]
which does not depend on $\sigma$. Choose $R_*$ and then $Z_*$ such that this is
$\le\frac18$, and $c_{\rm cl}$ with $(4\pi)^{-n/2}c_{\rm cl}\le\frac18$. Then
$\Theta^\chi_{q,t}(\sigma)\le\frac14$, a contradiction.
\end{proof}

\begin{remark}[Sheets at other heights and the transverse
density]\label{rem:mono-M}\leavevmode
\begin{enumerate}[(1)]
\item For $\chi\equiv1$, $\int\psi\rho\,d\mu$ is a chart version of the localised
transverse density $\Theta_{\rm loc}$ of
\cite[Proposition~4.16]{ChangWuZhang2026ContactBlowup}, with the cut-off $\psi$ of a
contact normal chart in place of the product cut-off there. The height-blind case of
Proposition~\ref{prop:mono} has the same strength as
\cite[Proposition~4.16]{ChangWuZhang2026ContactBlowup}. The error there, $E(t)\le
C_E(1+(T-t)^{-1/2})$, comes from the vertical cut-off, on whose support the transverse
kernel does not decay; here the corresponding term $R_4$ is bounded by means of (D), which
is itself derived from the last assertion of
\cite[Proposition~4.16]{ChangWuZhang2026ContactBlowup}.
\item With the height weight, the vertical cut-off is invisible at small scales: where
$|z_p|\ge z_0/2$ one has $|v|\ge z_0/(4\tau)$, so $\chi(v)\le e^{-az_0^2/(16\tau^2)}$.
Moreover, sheets at a height $\bar z\ne0$ above the centre do not contribute as
$\tau\to0$: Lemma~\ref{lem:density} gives $\Theta^\chi_{q,t}(\sigma)\to N(q,t)$ at every
point of the flow. The transverse density, by contrast, counts all sheets which approach
the Reeb orbit of the centre \cite[Remark~4.13]{ChangWuZhang2026ContactBlowup}, and its
localisation weights the sheets above a Reeb arc by the cut-off
\cite[Remark~4.18]{ChangWuZhang2026ContactBlowup}. Lemmas~\ref{lem:density}
and~\ref{lem:clearing} concern a density localised in the Reeb direction, at the parabolic
scale $|z_q|\lesssim\sigma$, of the kind that
\cite[Remark~4.13]{ChangWuZhang2026ContactBlowup} identifies as missing for a local
regularity theorem; such a theorem for the height-localised density is not proved here.
\item \cite[Remark~4.4]{ChangWuZhang2026ContactBlowup} explains why the full generator $X$
of the dilations cannot be used in a Gaussian kernel: $\frac14|X|^2=f+\frac14z^2$ is not
constant along the Reeb orbits, the vertical velocity $-2\theta\xi$ enters its evolution,
and the resulting expression is not a square. The kernel $\rho^\chi$ avoids this
difficulty. The vertical velocity enters only through the caloric height
$u=z_p-a_K(\tau)\theta$, for which $(\partial_t-\Delta)u=e_z$, and the height enters only
through $\chi(v)$; the inequality \eqref{eq:chi} turns the cross terms into the
non-negative defect $D^\chi$.
\end{enumerate}
\end{remark}

\begin{remark}[The static case]\label{rem:static}
Let $F_{t_*}$ be a minimal Legendrian immersion, $H=0$, so that $\theta$ is locally
constant, and consider the static flow $F_t:=\mathrm{Fl}^\xi_{-2\theta(t-t_*)}\circ
F_{t_*}$, the rotation along the Reeb field, componentwise. Since the Reeb flow preserves
$(\eta,\xi,\Phi,g)$, each $F_t$ is a minimal Legendrian immersion with the same $\theta$,
so $F$ satisfies (F) with $K=-2$, and $\mu_t=\mu_{t_*}$. In the chart of a point $q$,
where $\xi=2\partial_z$ and $f_q$ does not depend on $z$, the caloric height of
$F_{t_*-\sigma}$ with respect to the centre $(q,t_*)$ is $z_q\circ
F_{t_*-\sigma}-4\sigma\theta=z_q\circ F_{t_*}$. Hence the kernel $\rho^\chi_{q,t_*}$
evaluated on $F_{t_*-\sigma}$ is the static kernel
$(4\pi\sigma)^{-n/2}e^{-f_q/\sigma}\chi(z_q/\sigma)$ evaluated on the minimal Legendrian
$F_{t_*}$, up to the cut-off $\psi_q$, which is evaluated at points moved along the Reeb
orbits by at most $4C_\theta\sigma$. Proposition~\ref{prop:mono}, in which $E(t)$ is now a
constant, and Lemma~\ref{lem:density} give an almost monotone static density, with error
$O(\sigma)$, whose limit as $\sigma\to0$ is the number of sheets of $F_{t_*}$ through $q$.
Sheets which cross the Reeb orbit of $q$ at other heights contribute $0$ to the limit.

In the model $\Hn$, where a closed minimal Legendrian immersion does not exist, the same
construction applies to a proper minimal Legendrian immersion $F_{t_*}$ of a manifold
without boundary with $\sup|\theta|<\infty$ and $\mu_{t_*}(F_{t_*}^{-1}\{|w|<R\})\le
D_0R^n$ for all $R$. The static flow $F_t=\tau_{-4\theta(t-t_*)}\circ F_{t_*}$,
componentwise, satisfies the hypotheses of Proposition~\ref{prop:model} on every interval
$[t_0,t_*)$, with $t_*$ in the role of $T$. Indeed, the vertical displacement is at most
$4\sup|\theta|(t_*-t_0)$, which keeps $(x,t)\mapsto(F(x,t),t)$ proper, and $w\circ F_t$
and $\mu_t$ do not depend on $t$, so the angle and area bounds are inherited. For this
flow the monotonicity is exact. The area growth bound is not automatic: countably many
parallel planes $\{y=0,\ z=k\}$, $k\in\Z$, form a proper minimal Legendrian immersion with
$\theta=0$ for which it fails.
\end{remark}

\subsection{Rescaled variables}\label{subsec:scaling}
Let $F^i$ be the rescalings of Section~\ref{subsec:rescalings} at $(p,T)$, with
$\lambda=\lambda_i$, and write $(w',z',s)=(\lambda w,\lambda^2z_p,\lambda^2(t-T))$ for the
rescaled coordinates. Let $\mu'$ be the area measure of $F^i_s$ for $g_i$, $H^i$ its mean
curvature vector, and put
\[
\cS^i:=H^i+\frac{(\overline\nabla f)^\perp}{|s|},\qquad v^i:=\frac{z'}{|s|}-\hat
a_i(s)\theta^i,\qquad \hat a_i(s):=\frac{a_K(\lambda_i^{-2}|s|)}{\lambda_i^{-2}|s|},
\]
where $f=|w'|^2/16$ and the normal component, the gradient and the norms are those of
$g_i$; by Lemma~\ref{lem:moving}(a), $f$ is the rescaled potential $f_p$.

\begin{lemma}\label{lem:scaling}
With $\rho'=(4\pi|s|)^{-n/2}e^{-|w'|^2/(16|s|)}$,
\[
\rho_{p,T}\,d\mu_t=\rho'\,d\mu',\qquad v=v^i,\qquad |H|^2dt=|H^i|^2ds,\qquad
|\cS_p|^2dt=|\cS^i|^2ds,
\]
and $D^\chi dt=\big(|v^i\chi'(v^i)|/|s|+(\chi+16\chi'')(v^i)|\cS^i|^2\big)ds$. Moreover
$|\hat a_i(s)-4|\le8|K+2|\lambda_i^{-2}|s|$, so $\hat a_i\to4$ uniformly on $[-S_0,0)$,
and the rescaled cut-offs $\psi_i:=\psi_p\circ(\delta_{\lambda_i}\Upsilon_p)^{-1}$ equal
$1$ on $\{|w'|<\lambda_ir_0/2,\ |z'|<\lambda_i^2z_0/2\}$, hence converge to $1$ locally
uniformly.
\end{lemma}
\begin{proof}
On horizontal vectors $g_i=\lambda^2g^T$, so $d\mu'=\lambda^nd\mu$ and the $g_i$-norm of a
horizontal vector is $\lambda$ times its $g$-norm. The mean curvature vector for $g_i$ is
$H^i=\lambda^{-2}H$, and the two terms of $\cS^i$ and $\cS_p$ scale in the same way:
$(\overline\nabla f)^\perp/|s|=\lambda^{-2}(\overline\nabla f_p)^\perp/\tau$. Indeed,
$f=\lambda^2f_p$, $|s|=\lambda^2\tau$, the gradient of $\lambda^2f_p$ for $\lambda^2g^T$
is the gradient of $f_p$ for $g^T$, and the normal projections agree because the two
metrics are proportional on the horizontal distribution, which contains the tangent
planes. Hence $\cS^i=\lambda^{-2}\cS_p$ as vectors, $|H^i|=\lambda^{-1}|H|$,
$|\cS^i|=\lambda^{-1}|\cS_p|$, and $ds=\lambda^2dt$. With $\tau=\lambda^{-2}|s|$ and
$f_p=\lambda^{-2}f$ one has $\rho_{p,T}=\lambda^n\rho'$, which gives the first identity,
and $v=(\lambda^{-2}z'-a_K(\lambda^{-2}|s|)\theta)/(\lambda^{-2}|s|)=v^i$. The identity
for $D^\chi dt$ follows. The bound for $\hat a_i-4=4\eps_K(\lambda_i^{-2}|s|)$ is the
bound for $\eps_K$ stated at the beginning of this section, and the statement about
$\psi_i$ is the definition of the cut-off.
\end{proof}

\section{Weak contact tangent flows}\label{sec:weak}

In this section neither Type~I nor an orientation is assumed, except in
Theorem~\ref{thm:weak}(vi), and $n\ge1$ is arbitrary. We use the rescalings $F^i_s$,
$\theta^i_s$ of Section~\ref{subsec:rescalings} at a point $(p,T)$, for $s\in[-S_0,0)$,
and the notation of Section~\ref{subsec:scaling}.

The analysis follows the scheme of Neves \cite{nev07}. The dissipation of
Propositions~\ref{prop:theta2} and~\ref{prop:mono} makes the limits stationary
(Theorem~\ref{thm:weak}(i)) and conical (Theorem~\ref{thm:weak}(ii), (iv)). The defect
$|v\chi'(v)|/\tau$ pins the Legendrian angle to the height (Theorem~\ref{thm:weak}(iii)),
and the evolution of $\theta$ gives the transport of Proposition~\ref{prop:transport} and
Corollary~\ref{cor:angle-measure}, the Legendrian form of the last assertion of
\cite[Theorem~A]{nev07}.

We now explain the role of the Reeb height. In the model $\Hn$ the contact form is
$\overline\eta=\frac12dz+\overline\beta$ with
$\overline\beta=-\frac14\sum_j(y_j\,dx_j-x_j\,dy_j)$, so on a Legendrian immersion
$dz=-2\overline\beta$: the function $z/2$ is a primitive of the Liouville form
$-\overline\beta$ of $\C^n$ on every Legendrian. In a contact normal chart
$\eta=\frac12dz_p+\beta_p$, so that $d(z_p\circ F_t)=-2F_t^*\beta_p$ on a Legendrian
immersion, and $\beta_p$ agrees with $\overline\beta$ only to leading order,
$\beta_p-\overline\beta=O(|w|^3)$ (Lemma~\ref{lem:charts}). In this sense the Legendrian
flow is locally exact, and the Reeb height plays the role of Neves' primitive $\beta_t$ of
the Liouville form, which exists for exact Lagrangians. Indeed, along a flow
\eqref{eq:flow-model} the identity $(\partial_t-\Delta)z=-4\theta$ of
Lemma~\ref{lem:model-exact}(d) is the evolution formula of \cite[Proposition~6.1]{nev07}
with the primitive $z/2$ in place of $\beta_t$, and the caloric height $u=z-4\tau\theta$
is, up to the factor $2$, the combination $\beta_t+2(t-T)\theta_t$ of
\cite[Corollary~6.3]{nev07}, which evolves by the heat equation. In
Section~\ref{sec:mono} the caloric height enters the transverse Gaussian kernel through
the weight $\chi(v)$, $v=u/\tau$.

In the sketch of the proof of \cite[Theorem~B]{nev07} given in \cite[Section~3]{nev07},
for an exact Lagrangian under the hypotheses of that theorem, the primitive converges,
after passing to a subsequence, on the $j$-th cone of a tangent flow to
$\bar\beta_j-2s\bar\theta_j$, where the constants $\bar\beta_j$ are not controlled and may
differ from cone to cone. In $\C^n$ this is harmless, since a constant only changes the
primitive; in $\Hn$ a change of the constant is a vertical translation, and the transverse
density cannot detect it. In Theorem~\ref{thm:weak}(v), by contrast, the identity
$z=4|s|\bar\theta_j$ on the $j$-th lifted cone says that these constants vanish for the
components that are reached without vertical recentring. Components with a different
constant lie at other heights of the same Reeb orbit, and the weight separates them.

\subsection{Compactness and transport}\label{subsec:compactness}

\begin{definition}\label{def:weak}
Let $F^i_s,\theta^i_s$ be the rescalings of Section~\ref{subsec:rescalings}, for
$s\in[-S_0,0)$, defined on $\Sigma^i_s=F_{T+\lambda_i^{-2}s}^{-1}(U_p)$. In this section
and in Section~\ref{sec:applications}, $\mu$ denotes the area measure of $F^i_s$ for the
rescaled metric $g_i$, which is the measure $\mu'$ of Section~\ref{subsec:scaling}:
$\mu=\lambda_i^n\,\mu_{T+\lambda_i^{-2}s}$ on $\Sigma^i_s$, where $\mu_t$ is the area
measure of $F_t$ for $g$. (Where a proof passes from the original flow to its rescaling,
the original area measure is written $\mu_t$.) Put $\mu^i_s:=(F^i_s)_\#\mu$, a Radon
measure on $\Hn$ counted with multiplicity, and $\mathfrak
m^i_s:=(F^i_s,\theta^i_s)_\#\mu$ on $\Hn\times[-C_\theta,C_\theta]$. A \emph{weak contact
tangent flow} at $(p,T)$ is a limit $(\mu_s,\mathfrak m_s)_{s\in[-S_0,0)}$ with
$\mu^i_s\to\mu_s$ and $\mathfrak m^i_s\to\mathfrak m_s$ weakly, \emph{for every} $s$,
along one subsequence. The first marginal of $\mathfrak m_s$ is $\mu_s$, and $\mathfrak
m_s=\int\mathfrak m_{s,x}\,d\mu_s(x)$ is its disintegration (the \emph{Young measure of
the Legendrian angle}).

A \emph{weak contact tangent flow for all negative times} is a family $(\mu_s,\mathfrak
m_s)_{s<0}$ together with one subsequence of $\lambda_i$ such that $\mu^i_s\to\mu_s$ and
$\mathfrak m^i_s\to\mathfrak m_s$ weakly for every $s<0$ along that subsequence. For every
$k\in\N$ its restriction to $[-k,0)$ is then a weak contact tangent flow with $S_0=k$.
Such families are produced by the diagonal argument after the proof of
Theorem~\ref{thm:weak}; the statements of this paper made for every $s<0$ refer to them.
\end{definition}

\begin{proposition}\label{prop:compactness}
Assume (S), (F).
\begin{enumerate}[(a)]
\item $\mu^i_s(\{|w|<R\})\le D_0R^n$ for $R\le\lambda_ir_0$, all heights included: this is
(D) rescaled, since $g_i=\lambda_i^2(\delta_{\lambda_i}\Upsilon_p)_*g$ on horizontal
vectors, so that the $g_i$-area of a Legendrian piece is $\lambda_i^n$ times its $g$-area,
and $\{|w|<R\}$ in the rescaled chart is $\{|w|<R/\lambda_i\}$ in $U_p$. In particular
$\mu^i_s(\mathcal K)\le C(\mathcal K)$ for compact $\mathcal K\subset\Hn$.
\item $\int_{-S_0}^{-\delta}\int_{(F^i_s)^{-1}(\mathcal K)}|H^i|^2\,d\mu\,ds\le C(\mathcal
K,S_0,\delta)$. This follows from Proposition~\ref{prop:theta2} and the scale invariance
of $|H|^2\rho\,d\mu\,dt$, since the rescaled cut-off equals $1$ on $\mathcal K$ for $i$
large, and $\rho'\ge(4\pi S_0)^{-n/2}e^{-R^2/(16\delta)}$ on $\mathcal
K\times[-S_0,-\delta]$ if $\mathcal K\subset\{|w|\le R\}$.
\item There is a subsequence along which $\mu^i_s\to\mu_s$ and $\mathfrak
m^i_s\to\mathfrak m_s$ for every $s\in[-S_0,0)$. Along a further subsequence, independent
of $s$, there is a set $\mathcal G\subset[-S_0,0)$ of full measure such that for every
$s\in\mathcal G$ and every compact $\mathcal K$,
\[
\int_{(F^i_s)^{-1}(\mathcal K)}\big(|H^i|^2+|\mathcal
S^i|^2+|v^i\chi'(v^i)|\big)\,d\mu\longrightarrow0 .
\]
For $s\in\mathcal G$ there is an integral Legendrian varifold $V_s$ in $(\Hn,\overline g)$
with weight $\mu_s$; the varifolds of $F^i_s$ in $(\Hn,g_i)$ converge to $V_s$, as Radon
measures on the Grassmann bundle $G_n(T\Hn)$, along a further, $s$-dependent subsequence.
Since an integral varifold is determined by its weight, $V_s$ does not depend on that
subsequence.
\item The projection does not commute with the limit:
$\pi_\#\mu_s\le\liminf_i\pi_\#\mu^i_s$, with strict inequality in general; for the two
sheets of Example~\ref{ex:two-sheets}, rescaled at a point of the lower one, the upper
sheet escapes vertically while its projection stays. With $L=\R^n$ and
$0<c<\min\{z_0,\frac12\}$ the two sheets descend to disjoint static Legendrian tori in the
compact quotient $M$ of Example~\ref{ex:potential-not-phase}; this flow satisfies (S) and
(F) with $\theta\equiv0$. For regular $M$, $\pi_\#\mu^i_s$ is at most the rescaled
projected measure $\nu^{\lambda_i}_s$ of \cite[\S4.2]{ChangWuZhang2026ContactBlowup},
which counts all sheets over the transverse chart.
\end{enumerate}
Here $\mathcal S^i=H^i+(\overline\nabla f)^\perp/|s|$ and $v^i=z/|s|-\hat a_i(s)\theta^i$,
computed in the rescaled structure.
\end{proposition}

\begin{proposition}\label{prop:energy}
Assume (S), (F). For every $\lambda_i\to\infty$, every compact $\mathcal K\subset\Hn$ and
$0<\delta<S_0$,
\[
\int_{-S_0}^{-\delta}\!\!\int_{(F^i_s)^{-1}(\mathcal K)}\Big(|H^i|^2+|\mathcal
S^i|^2+|v^i\chi'(v^i)|\Big)\,d\mu\,ds\longrightarrow0 .
\]
\end{proposition}
\begin{proof}
By Propositions~\ref{prop:theta2} and~\ref{prop:mono}, including the height-blind version,
the functions
\[
\mathcal M_\chi(t):=\int\psi\rho^\chi\,d\mu_t-\int_{T-\tau_0}^tE,\qquad \mathcal
M_2(t):=\int\psi\theta^2\rho\,d\mu_t-C_Ht,\qquad \mathcal M_0(t):=\int\psi\rho\,d\mu_t-Ct
\]
are non-increasing and bounded below on $[T-\tau_0,T)$, hence convergent as $t\to T$;
neither Type~I nor any compactness hypothesis enters. Let
$I_i=[t_i^-,t_i^+]:=[T-S_0\lambda_i^{-2},T-\delta\lambda_i^{-2}]$. Integrating the three
differential inequalities over $I_i$ gives
\[
\int_{I_i}\!\!\int\psi\Big(\tfrac34D^\chi+|H|^2+|\mathcal
S_p|^2\Big)\rho\,d\mu_t\,dt\le\sum_{\bullet\in\{\chi,2,0\}}\big(\mathcal
M_\bullet(t_i^-)-\mathcal M_\bullet(t_i^+)\big)\longrightarrow0 .
\]
In the rescaled variables, by Lemma~\ref{lem:scaling}, $\rho\,d\mu_t=\rho'\,d\mu'$ with
$\mu'=\mu$ the rescaled area measure, $|H|^2dt=|H^i|^2ds$, $|\mathcal S_p|^2dt=|\mathcal
S^i|^2ds$ and $D^\chi dt\ge|v^i\chi'(v^i)|\,ds/|s|$. The rescaled cut-off equals $1$ on
$\mathcal K$ for $i$ large, $\rho'=(4\pi|s|)^{-n/2}e^{-|w|^2/(16|s|)}\ge(4\pi
S_0)^{-n/2}e^{-R^2/(16\delta)}$ on $\mathcal K\times[-S_0,-\delta]$ if $\mathcal
K\subset\{|w|\le R\}$, and $1/|s|\ge1/S_0$. This gives the claim.
\end{proof}

\begin{proof}[Proof of Proposition~\ref{prop:compactness}(c)]
Fix compact sets $\mathcal K_1\subset\mathcal K_2\subset\dots$ exhausting $\Hn$, with
$\mathcal K_m\subset\mathrm{int}\,\mathcal K_{m+1}$. Two compactness problems have to be
separated. The measures $\mathfrak m^i_s$ converge at every time along one subsequence,
by a bound on their time derivative; the varifolds of $F^i_s$ converge only at good times,
where the energy of Proposition~\ref{prop:energy} vanishes in the limit, and along a
subsequence that may depend on the time.

\emph{Convergence at every time.} Let $Q$ be a continuous map from
$[-S_0,0)$ to $C^2_c(\Hn\times\R)$ which is $C^1$ as a map to $C^0_c(\Hn\times\R)$, with
$x$-supports in a fixed compact set $\mathcal K'$. Write $Q$ as a function of
$(x,\theta,s)$. Since $\partial_sF^i=H^i-2\theta^i\xi_i$, $\partial_s\,d\mu=-|H^i|^2d\mu$,
$\partial_s\theta^i=\Delta\theta^i+\lambda_i^{-2}(K+2)\theta^i$ and $\xi_i=2\partial_z$,
integration by parts of the $\Delta\theta$ term gives
\begin{multline}\label{eq:dnu}
\frac{d}{ds}\mathfrak m^i_s(Q_s)=\mathfrak
m^i_s\big(\partial_sQ-4\theta\,\partial_zQ\big)+\int\Big[d_xQ(H^i)-d_x\partial_\theta
Q(\nabla\theta^i)-\partial_\theta^2Q\,|\nabla\theta^i|^2\\
+\lambda_i^{-2}(K+2)\,\theta^i\partial_\theta Q-Q|H^i|^2\Big]d\mu .
\end{multline}
Since $|\nabla\theta^i|=|H^i|$ and $|\theta^i|\le C_\theta$, the last integral is bounded
by $C(Q)\int_{(F^i_s)^{-1}(\mathcal K')}\big(|H^i|+|H^i|^2+\lambda_i^{-2}\big)d\mu$. Here
the integration by parts takes place on $\Sigma^i_s$, where the integrand has compact
support, since the $x$-support of $Q$ lies in $\mathcal K'$, which is contained in the
rescaled box for $i$ large. The formula extends to functions $Q$ that are merely $C^2$ in
$\theta$, with $x$-supports in $\mathcal K'$ for $|\theta|\le C_\theta+2$. To see this,
fix $\kappa\in C^\infty_c((-C_\theta-2,C_\theta+2))$ with $\kappa=1$ on
$[-C_\theta-1,C_\theta+1]$ and apply \eqref{eq:dnu} to $\kappa(\theta)Q$, which belongs to
the class above. Both sides of \eqref{eq:dnu} only involve the values and the
$\theta$-derivatives of $Q$ at $\theta=\theta^i_s\in[-C_\theta,C_\theta]$, where $\kappa
Q$ and $Q$ agree with all their derivatives. This is how \eqref{eq:dnu} is used below and
in the proof of Proposition~\ref{prop:transport}.

For $Q=\gamma\otimes\varpi$ independent of $s$, \eqref{eq:dnu}, (a), (b) and
Cauchy--Schwarz show that the total variation of $s\mapsto\mathfrak
m^i_s(\gamma\otimes\varpi)$ on $[-S_0,-\delta]$ is bounded uniformly in $i$, and
$\sup_s|\mathfrak
m^i_s(\gamma\otimes\varpi)|\le\|\gamma\|_\infty\|\varpi\|_{C^0[-C_\theta,C_\theta]}C(\spt\gamma)$.
Let $\{\gamma_k\}\subset C^2_c(\Hn)$ be countable and dense in $C_c(\Hn)$ in the
inductive-limit sense, and let $\varpi_l(\theta)=\theta^l$, $l\ge0$. For each $(k,l,m)$,
Helly's selection theorem applied to $s\mapsto\mathfrak m^i_s(\gamma_k\otimes\varpi_l)$ on
$[-S_0,-1/m]$ yields a subsequence converging at \emph{every} $s$, and a diagonal argument
gives one subsequence along which $\mathfrak m^i_s(\gamma_k\otimes\varpi_l)$ converges for
all $k,l$ and all $s\in[-S_0,0)$. To identify the limit, fix $s$. The $\mathfrak m^i_s$
are supported in $\Hn\times[-C_\theta,C_\theta]$, with mass bounded on compact sets. By
Stone--Weierstrass on $\mathcal K_m\times[-C_\theta,C_\theta]$ the span of
$\{\gamma_k\otimes\varpi_l\}$ is dense, so all subsequential weak limits agree, and
$\mathfrak m^i_s\to\mathfrak m_s$ for every $s$. No mass escapes in the
$\theta$-direction, so $\mu^i_s=(\mathrm{pr}_1)_\#\mathfrak
m^i_s\to\mu_s:=(\mathrm{pr}_1)_\#\mathfrak m_s$.

\emph{A common set of good times.} By Proposition~\ref{prop:energy} we may pass to a
further subsequence $(i_k)$ with
\[
\int_{-S_0}^{-1/m}\!\!\int_{(F^{i_k}_s)^{-1}(\mathcal K_m)}\big(|H|^2+|\mathcal
S|^2+|v\chi'(v)|\big)\,d\mu\,ds\le2^{-k}\qquad\text{for all }m\le k .
\]
For each fixed $m$ the sum over $k$ of these integrals is finite. By Fubini, for
a.e.\ $s\in[-S_0,-1/m]$ the sum over $k$ of the inner integrals is finite, so they tend to
$0$.
Let $\mathcal G$ be the set of $s$ for which this holds for every $m$ with $s\le-1/m$; it
has full measure, and the subsequence does not depend on $s$. The convergence $\mathfrak
m^i_s\to\mathfrak m_s$ at every time persists along it. From now on further subsequences
may depend on the fixed time $s\in\mathcal G$.

\emph{The varifold limit at good times.} To construct $V_s$ for $s\in\mathcal G$, we embed
the rescaled metrics isometrically into a fixed Euclidean space, so that Allard's
compactness theorem for varifolds in a fixed Euclidean space applies. Fix $m$ and
Euclidean balls $\mathcal K_{m+1}\subset\mathbb
B_1\Subset\mathbb B_2\Subset\mathbb B$ in the coordinates of $\R^{2n+1}$; after a linear
change of coordinates $\mathbb B$ is the open unit ball, and $\mathbb B_1,\mathbb
B_2,\mathbb B$ play the roles of $B_1,B_2,B$ in \cite[Satz~5.5]{gu89}. Choose a closed
manifold $N_{\mathbb B}\supset\overline{\mathbb B}$ (for instance a flat torus
$\R^{2n+1}/L\Z^{2n+1}$ with $L$ large) and a smooth metric on $N_{\mathbb B}$ that equals
$\overline g$ on $\overline{\mathbb B}$. By \cite[Satz~1.2 and Theorem~1.1]{gu89}, which
give \cite[Folgerung~1.3]{gu89}, $N_{\mathbb B}$ has a free isometric embedding into some
$\R^{n_G}$; let $u_0\in C^\infty(\overline{\mathbb B},\R^{n_G})$ be its restriction to
$\overline{\mathbb B}$. Let $\zeta_{\mathbb B}\in C^\infty_c(\mathbb B_1)$ with
$\zeta_{\mathbb B}=1$ near $\mathcal K_{m+1}$, and $\hat g_i:=\zeta_{\mathbb
B}\,(g_i-\overline g)$. Then $\spt\hat g_i\subset\mathbb B_1$ and $\hat g_i\to0$ in every
$C^{k,1/2}$, since $g_i\to \overline g$ in $C^\infty_{\rm loc}$. The smallness hypothesis
of \cite[Satz~5.5]{gu89} is that $\|E(u_0)\|_{\rm op}\,\|E(u_0)(0,\hat g_i)\|_{C^{2,1/2}}$
is at most a constant that does not depend on $u_0$; here $E(u_0)$ is the linear operator
of \cite[(5.5)]{gu89} determined by $u_0$, and $\|E(u_0)\|_{\rm op}$ is its norm as an
operator on $C^{2,1/2}$. As $u_0$ is fixed, this product is at most $\|E(u_0)\|_{\rm
op}^2\|\hat g_i\|_{C^{2,1/2}}\to0$, so the hypothesis holds for $i$ large. By
\cite[Satz~5.5, (5.14)--(5.15)]{gu89} there are then $u_i\in C^\infty(\overline{\mathbb
B},\R^{n_G})$ with $\spt u_i\subset\mathbb B_2$, $\|u_i\|_{C^{2,1/2}}\to0$, and
$(u_0+u_i)^*g_{\rm euc}=\overline g+\hat g_i$ on $\mathbb B$, where $g_{\rm euc}$ is the
Euclidean metric of $\R^{n_G}$. Thus $\tilde u_i:=u_0+u_i$ is an isometric immersion of a
neighbourhood of $(\mathcal K_{m+1},g_i)$ into the fixed $\R^{n_G}$, and $\tilde u_i\to
u_0$ in $C^{2,1/2}$.

Now fix $s\in\mathcal G$, and put $\mathcal O:=\mathrm{int}\,\mathcal
K_{m+1}\supset\mathcal K_m$. Let $\mathcal U\subset\R^{n_G}$ be the open $\delta_{\mathcal
U}$-neighbourhood of $u_0(\mathcal K_m)$, where $2\delta_{\mathcal
U}<\mathrm{dist}\big(u_0(\mathcal K_m),u_0(\mathcal K_{m+1}\setminus\mathcal O)\big)$; the
distance is positive because $u_0$ is injective. Put $\mathcal O_{\mathcal U}:=\mathcal
O\cap u_0^{-1}(\mathcal U)$, an open set containing $\mathcal K_m$. Then
\begin{equation}\label{eq:OU}
\mathcal U\cap u_0(\mathcal K_{m+1})=u_0(\mathcal O_{\mathcal U}),
\end{equation}
so that for every compact $\mathcal K'\subset\mathcal U$ the set $\mathcal K_{m+1}\cap
u_0^{-1}(\mathcal K')$ is a compact subset of $\mathcal O_{\mathcal U}$; and, since
$\tilde u_i\to u_0$ uniformly, $\tilde u_i(\mathcal K_{m+1}\setminus\mathcal
O)\cap\mathcal U=\emptyset$ for $i$ large. Let $V^i$ be the varifold of $F^i_s$ in
$(\Hn,g_i)$, that is, a Radon measure on the Grassmann bundle $G_n(T\Hn)$ of $n$-planes.
Let $W^i$ be the push-forward under $\tilde u_i$ of $V^i\llcorner G_n(T\mathcal O)$,
i.e.\ the varifold of $\tilde u_i\circ F^i_s$ on $(F^i_s)^{-1}(\mathcal O)$. For $Y\in
C^1_c(\mathcal U,\R^{n_G})$, the function $Y\circ\tilde u_i\circ F^i_s$ has compact
support in the open set $(F^i_s)^{-1}(\mathcal O)$ of the closed manifold $\Sigma$, by the
last observation. Thus the first variation has no boundary term:
\[
\delta W^i(Y)=-\int\big\langle H^i+\tr_S\mathrm{II}_{\tilde u_i},Y\big\rangle\,d\mu ,
\]
where $S$ is the tangent plane. Hence $\|W^i\|(\mathcal U')+\|\delta W^i\|(\mathcal U')\le
C$ for every $\mathcal U'\Subset\mathcal U$, by (a) and the choice of $\mathcal G$.
Allard's compactness theorem for integer multiplicity varifolds \cite[42.7, 42.8]{si83}
gives a further subsequence with $W^i\to W$ in $\mathcal U$, $W$ integral with locally
bounded first variation.

Let $\hat u_i,\hat u_0\colon G_n(T\mathcal K_{m+1})\to G_n(\R^{n_G})$ be the induced maps
of Grassmann bundles, $\hat u_i(x,S)=(\tilde u_i(x),d\tilde u_i(S))$ and $\hat
u_0(x,S)=(u_0(x),du_0(S))$. Then $\hat u_i\to\hat u_0$ uniformly, and $\hat u_0$ is a
homeomorphism onto its image, because $u_0$ is an embedding. The varifold $W^i$ is carried
by $\hat u_i(G_n(T\mathcal K_{m+1}))$. Hence the support of $W$ lies in $\hat
u_0(G_n(T\mathcal K_{m+1}))\cap G_n(\mathcal U)=\hat u_0(G_n(T\mathcal O_{\mathcal U}))$,
by \eqref{eq:OU}, and
\[
W=(u_0)_\#V_s\quad\text{on }\mathcal U,\qquad V_s:=(\hat u_0^{-1})_\#W ,
\]
a varifold in $\mathcal O_{\mathcal U}$. There is no Jacobian factor, since $u_0$ is
isometric for $\overline g$. $V_s$ is integral, because $W$ is and $u_0$ is a smooth
embedding. Its weight is $\mu_s$: for $\gamma\in C_c(\mathcal U)$, the function
$\gamma\circ u_0$ vanishes on $\mathcal K_{m+1}$ outside a compact subset of $\mathcal
O_{\mathcal U}$, by \eqref{eq:OU}, and $\gamma\circ\tilde u_i-\gamma\circ u_0\to0$
uniformly, so that
\[
\|W\|(\gamma)=\lim_i\int_{\mathcal O}\gamma\circ\tilde u_i\,d\mu^i_s=\int_{\mathcal
O_{\mathcal U}}\gamma\circ u_0\,d\mu_s ,
\]
that is, $\|V_s\|=\mu_s$ on $\mathcal O_{\mathcal U}$. Moreover $V^i\to V_s$ on
$G_n(T\mathcal O_{\mathcal U})$: for $\varphi\in C_c(G_n(T\mathcal O_{\mathcal U}))$,
Tietze's theorem gives $\tilde\varphi\in C_c(G_n(\mathcal U))$ with
$\tilde\varphi\circ\hat u_0=\varphi$ on $G_n(T\mathcal K_{m+1})$, and since $\hat
u_i\to\hat u_0$ uniformly and the masses are bounded,
$V^i(\varphi)=W^i(\tilde\varphi)+o(1)\to W(\tilde\varphi)=V_s(\varphi)$. Finally, the
tangent planes of $F^i_s$ lie in the closed set $\mathrm{Leg}(\eta_i,g_i)$ of horizontal
planes that are Lagrangian for $d\eta_i$. These sets converge in the Hausdorff sense to
$\mathrm{Leg}(\overline\eta,\overline g)$, so the approximate tangent planes of $V_s$ are
Legendrian $\|V_s\|$-a.e.

For $m<m'$ the varifolds obtained in this way on the corresponding open sets, which
contain $\mathcal K_m$ and $\mathcal K_{m'}$, are integral and have the same weight
$\mu_s$ on the intersection. Hence they agree there, since a rectifiable varifold is
determined by its weight, and they define one integral Legendrian varifold $V_s$ in
$(\Hn,\overline g)$ with $\|V_s\|=\mu_s$. A diagonal argument in $m$ gives a subsequence
along which $V^i\to V_s$ on $G_n(T\Hn)$.
\end{proof}

\begin{proposition}\label{prop:transport}
Assume (S) and (F), and let $(\mu_s,\mathfrak m_s)$, $s\in[-S_0,0)$, be a weak contact
tangent flow at $(p,T)$. For $\sigma\in\R$ define
\[
\mathrm{Sh}_\sigma\colon\Hn\times\R\to\Hn\times\R,\qquad
\mathrm{Sh}_\sigma(w,z,\theta)=(w,\ z-4\sigma\theta,\ \theta).
\]
Then $\mathfrak m_s=(\mathrm{Sh}_{s-s_0})_\#\mathfrak m_{s_0}$ for all $s,s_0\in[-S_0,0)$.
Equivalently, $\partial_s\mathfrak m+\partial_z(-4\theta\,\mathfrak m)=0$: the part of
$\mathfrak m$ with Legendrian angle $\theta$ is translated vertically with speed
$-4\theta$.
\end{proposition}
\begin{proof}
Let $\gamma\in C^2_c(\Hn)$ and $\varpi\in C^2(\R)$, and put
$Q_s(w,z,\theta):=\gamma(w,z+4(s-s_0)\theta)\,\varpi(\theta)$. For $|\theta|\le
C_\theta+2$ and $s$ in a compact subinterval $[s_1,s_2]\subset[-S_0,0)$, the $x$-supports
lie in a fixed compact $\mathcal K'$. Moreover
\[
\partial_sQ_s=4\theta\,(\partial_z\gamma)(w,z+4(s-s_0)\theta)\,\varpi(\theta)=4\theta\,\partial_zQ_s
,
\]
so the first term on the right of \eqref{eq:dnu} vanishes identically. Integrating
\eqref{eq:dnu} from $s_0$ to $s$ (both in $[s_1,s_2]$) gives
\[
\big|\mathfrak m^i_s(Q_s)-\mathfrak m^i_{s_0}(Q_{s_0})\big|\le
C\int_{s_1}^{s_2}\!\!\int_{(F^i_{s'})^{-1}(\mathcal
K')}\big(|H^i|+|H^i|^2+\lambda_i^{-2}\big)\,d\mu\,ds' .
\]
By Proposition~\ref{prop:energy}, Proposition~\ref{prop:compactness}(a) and
Cauchy--Schwarz, the right-hand side tends to $0$. By every-time convergence, $\mathfrak
m^i_s(Q_s)\to\mathfrak m_s(Q_s)$ and $\mathfrak m^i_{s_0}(Q_{s_0})\to\mathfrak
m_{s_0}(\gamma\otimes\varpi)$. Hence
\[
\mathfrak m_s\big((\gamma\otimes\varpi)\circ\mathrm{Sh}_{s_0-s}\big)=\mathfrak
m_s(Q_s)=\mathfrak m_{s_0}(\gamma\otimes\varpi).
\]
Functions $\gamma\otimes\varpi$ span a dense subspace, so
$(\mathrm{Sh}_{s_0-s})_\#\mathfrak m_s=\mathfrak m_{s_0}$. Since
$\mathrm{Sh}_a\circ\mathrm{Sh}_b=\mathrm{Sh}_{a+b}$, this is the claim.
\end{proof}

\subsection{Structure of weak contact tangent flows}\label{subsec:structure}

\begin{theorem}\label{thm:weak}
Assume (S), (F); neither Type~I nor almost calibration is assumed, and $n\ge1$ is
arbitrary. Let $(\mu_s,\mathfrak m_s)$ be a weak contact tangent flow at $(p,T)$ along the
subsequence of Proposition~\ref{prop:compactness}(c). Then for a.e.\ $s\in[-S_0,0)$:
\begin{enumerate}[(i)]
\item $V_s$ is a stationary Legendrian integral varifold in $(\Hn,\overline g)$.
\item $(\overline\nabla f)^{\perp}=0$ for $V_s$-a.e.\ $(x,S)$: the horizontal dilation
field is tangent.
\item \emph{The Legendrian angle is pinned to the height:} $\mathfrak
m_{s,x}=\delta_{z(x)/(4|s|)}$ for $\mu_s$-a.e.\ $x$. In particular
$\spt\mu_s\subset\{|z|\le4|s|C_\theta\}$.
\item There are finitely many values $\bar\theta_1<\dots<\bar\theta_N$ in
$[-C_\theta,C_\theta]$, with $N\le CD_0$, and nonzero stationary integral Lagrangian cones
$\mathcal C_1,\dots,\mathcal C_N$ in $(\C^n,g_0)$ such that
\[
V_s=\sum_j\mathcal C_j\times\{4|s|\bar\theta_j\}.
\]
Each lifted cone is Legendrian, and on the $j$-th slice the Legendrian angle equals
$\bar\theta_j$.
\end{enumerate}
Moreover:
\begin{enumerate}[(i)]\setcounter{enumi}{4}
\item \emph{All times.} The values $\bar\theta_j$ and the cones $\mathcal C_j$ do not
depend on $s$, and for \emph{every} $s\in[-S_0,0)$
\[
\mathfrak m_s=\sum_j\big(\|\mathcal
C_j\|\times\delta_{4|s|\bar\theta_j}\big)\otimes\delta_{\bar\theta_j},\qquad
\mu_s=\sum_j\|\mathcal C_j\|\times\delta_{4|s|\bar\theta_j} ;
\]
moreover, for every $s\in[-S_0,0)$, $\mu_s$ is the weight of the stationary Legendrian
integral varifold $\sum_j\mathcal C_j\times\{4|s|\bar\theta_j\}$ in $(\Hn,\overline g)$,
each lifted cone being stationary by itself.
\item If $\Sigma$ is orientable, then for each $j$ there are distinct
$\bar\vartheta_{j1},\dots,\bar\vartheta_{jN_j}\in\R/2\pi\Z$, with $N_j\le C(n,D_0)$, and
nonzero integer multiplicity $n$-currents $\mathcal T_{jk}$ in $\C^n$, with locally finite
mass and $\partial \mathcal T_{jk}=0$, such that
\[
\mathcal C_j=\sum_k|\mathcal T_{jk}|\quad\text{as varifolds},
\]
each $\mathcal T_{jk}$ is calibrated by $\mathrm{Re}(e^{-i\bar\vartheta_{jk}}\Omega_0)$,
and each $\mathcal T_{jk}$ is invariant under the dilations $w\mapsto\lambda w$,
$\lambda>0$. Thus $\mathcal C_j$ is a finite sum of special Lagrangian cones with
\emph{geometric} phases $\bar\vartheta_{jk}$.
\end{enumerate}
\end{theorem}
\begin{proof}
Let $\mathcal G$ be the set of Proposition~\ref{prop:compactness}(c). At good times we
first show that $V_s$ is stationary, that the horizontal dilation field is tangent and
that the Legendrian angle is determined by the height, (i)--(iii). These properties force
a decomposition into finitely many lifted cones (iv), which the transport of
Proposition~\ref{prop:transport} extends to every time (v); the orientation is used only
for the decomposition by geometric phases in (vi).

\emph{(i).} Fix $s\in\mathcal G$, and let $V^i$ be the varifolds of $F^i_s$ in
$(\Hn,g_i)$, along a subsequence with $V^i\to V_s$ on $G_n(T\Hn)$ (end of the proof of
Proposition~\ref{prop:compactness}(c)). Let $Y\in C^1_c(\Hn,T\Hn)$. For a metric $g$ write
$\mathrm{div}^g_SY:=\tr\big(P^g_S\circ\overline\nabla^gY\big)$, where $P^g_S$ is the
$g$-orthogonal projection onto $S$. Since $\Sigma$ is closed, the first variation formula
for the immersion $F^i_s$ in $(\Hn,g_i)$ has no boundary term:
\[
\int\mathrm{div}^{g_i}_SY\,dV^i(x,S)=-\int g_i\big(H^i,Y\circ F^i_s\big)\,d\mu .
\]
The right-hand side tends to $0$, by Cauchy--Schwarz, (a) and $s\in\mathcal G$. Since
$g_i\to \overline g$ in $C^1$ on $\spt Y$, the integrands on the left converge uniformly
on $G_n(T\spt Y)$ to $\mathrm{div}^{\overline g}_SY$, and the masses $V^i(G_n(T\spt Y))$
are bounded. Thus the left-hand side tends to $\delta^{\overline
g}V_s(Y)=\int\mathrm{div}^{\overline g}_SY\,dV_s$. Hence $\delta^{\overline g}V_s(Y)=0$,
and $V_s$ is stationary in $(\Hn,\overline g)$.

\emph{(ii).} For $s\in\mathcal G$, $\int_{(F^i_s)^{-1}(\mathcal K)}|(\overline\nabla
f)^\perp|^2\le2|s|^2\int(|\mathcal S^i|^2+|H^i|^2)\to0$. The integrand
$|P_{S^\perp}\overline\nabla f|^2$ is continuous on the Grassmann bundle (with $g_i\to
\overline g$), so (ii) follows from varifold convergence.

\emph{(iii).} Let $\gamma\in C_c(\Hn)$, $\gamma\ge0$, and
$Q_i(x,\theta,s):=|v\chi'(v)|\gamma(x)$ with $v=z/|s|-\hat a_i(s)\theta$; let $Q$ be the
same with $\hat a_i$ replaced by $4$. Then $Q_i\to Q$ uniformly on compact sets. For each
$s$, $\mathfrak m^i_s(Q_i)\to\mathfrak m_s(Q)$, and by Proposition~\ref{prop:energy} and
Fatou's lemma in $s$,
\[
\int_{-S_0}^{-\delta}\mathfrak m_s(Q)\,ds\le\liminf_i\int_{-S_0}^{-\delta}\mathfrak
m^i_s(Q_i)\,ds=0 .
\]
Taking $\gamma$ in a countable dense family, $v=0$ $\mathfrak m_s$-a.e.\ for a.e.\ $s$.
Since $v\chi'(v)=0$ only at $v=0$, this means $\theta=z/(4|s|)$, i.e.\ $\mathfrak
m_{s,x}=\delta_{z(x)/(4|s|)}$. As $|\theta|\le C_\theta$, the support lies in the slab.

\emph{(iv).} Fix $s$ for which (i)--(iii) hold. The mechanism is the following: the
restriction of $V_s$ to any height interval is again stationary, its projection is a
stationary integral cone, and a nonzero such cone carries Gaussian mass at least $1$;
since the total Gaussian mass is bounded, only finitely many heights can occur. For
$V_s$-a.e.\ $(x,S)$ the plane $S$ is
horizontal Legendrian and $(\overline\nabla f)^{\perp_S}=0$, so the exact model identity
of Lemma~\ref{lem:model}(d) gives $\nabla^Sz=0$. Hence, for $\varpi\in C^1(\R)$ and $Y\in
C^1_c$,
\[
\delta V_s\big((\varpi\circ
z)Y\big)=\int\varpi(z)\,\mathrm{div}_SY\,dV_s+\int\varpi'(z)\langle\nabla^Sz,Y\rangle\,dV_s=\int\varpi(z)\,\mathrm{div}_SY\,dV_s
,
\]
and the left side vanishes by (i). Approximating the indicator of an interval $I$ by
smooth $\varpi$ (bounded pointwise convergence) shows that
$V_I:=V_s\llcorner\{(x,S):z(x)\in I\}$ is stationary, for every interval $I$.

The projection $\pi\colon(\Hn,\overline g)\to(\C^n,g_0)$ is a Riemannian submersion with
totally geodesic fibres. For a horizontal plane $S$ and the horizontal lift $\tilde Y$ of
$Y\in C^1_c(\C^n)$, O'Neill's formula gives $\mathrm{div}_S\tilde
Y=(\mathrm{div}_{\pi_*S}Y)\circ \pi$. The support of $V_I$ over $\spt Y$ is compact: it
lies in the slab, and its mass is bounded by (a). Multiplying $\tilde Y$ by a function of
$z$ equal to $1$ on the slab, which is allowed since $\nabla^Sz=0$, we obtain
$\delta(\pi_\#V_I)(Y)=\delta V_I(\tilde Y)=0$. The push-forward $\pi_\#V_I$ is an integral
varifold, since $\pi$ is isometric on horizontal planes and multiplicities add over
preimages. Its tangent planes are Lagrangian, since horizontal Legendrian planes project
to Lagrangian planes.

By (ii) the radial field $\mathfrak r$ (the projection of $2\overline\nabla f$) is tangent
to $\pi_\#V_I$ a.e. For a stationary rectifiable varifold $V$ in a flat space with
$\mathfrak r$ tangent a.e., and $\gamma\in C^1_c$, put $\mathcal
J(\lambda)=\lambda^{-n}\int\gamma(x/\lambda)\,d\|V\|$. With $Y(x)=\gamma(x/\lambda)x$,
\begin{equation}\label{eq:cone-homog}
\begin{aligned}
0=\delta V(Y)&=\int\Big(n\gamma(x/\lambda)+\lambda^{-1}D\gamma(x/\lambda)\cdot
P_Sx\Big)dV\\
&=\int\Big(n\gamma(x/\lambda)+\lambda^{-1}D\gamma(x/\lambda)\cdot
x\Big)d\|V\|=-\lambda^{n+1}\mathcal J'(\lambda).
\end{aligned}
\end{equation}
Hence $\mathcal J$ is constant, that is, $\|V\|(\lambda A)=\lambda^n\|V\|(A)$ for all
$\lambda>0$ and Borel sets $A$; so $\|V\|$, and hence $V$, is invariant under the varifold
push-forwards by the dilations, since a rectifiable varifold is determined by its weight.
Thus $\pi_\#V_I$ is a stationary integral Lagrangian cone with vertex $0$.

A nonzero stationary integral cone $\mathcal C$ has vertex density $\Theta(\mathcal
C,0)\ge1$, by upper semicontinuity of the density and $\Theta\ge1$ a.e., and its Gaussian
integral $\int\rho'\,d\|\mathcal C\|$ equals $\Theta(\mathcal C,0)$. Here
$\rho'=(4\pi|s|)^{-n/2}e^{-|w|^2/(16|s|)}$ is the $n$-dimensional backward Gaussian kernel
of Huisken type written in the $g_0$-distance $|w|/2$, normalised so that its integral
over any $n$-plane through $0$, with respect to $g_0$-area, equals $1$
(Example~\ref{ex:two-sheets}). On the other hand $\int\rho'\,d\mu_s\le2^nD_0/\omega_n$, by
(a) and the layer-cake formula: by (a) and lower semicontinuity on open sets,
$\mu_s(\{|w|<R\})\le D_0R^n$ for all $R>0$, and for $\rho'=\varphi(|w|)$ one has
$\int_0^\infty nR^{n-1}\varphi(R)\,dR=\omega_n^{-1}
\int_{\R^n}\varphi(|y|)\,dy=2^n/\omega_n$. If $I_1,\dots,I_N$ are disjoint open intervals
each meeting the support of $z_\#\mu_s$, then $z_\#\mu_s(I_k)>0$, so each $\pi_\#V_{I_k}$
is nonzero, and
\[
N\le\sum_k\int\rho'\,d\|\pi_\#V_{I_k}\|\le\int\rho'\,d\mu_s\le2^nD_0/\omega_n .
\]
Thus $z_\#\mu_s$ is supported on at most $2^nD_0/\omega_n$ points $\bar z_j$, and the
constant $C$ in (iv) can be taken to be $2^n/\omega_n$. Put $\mathcal C_j:=\pi_\#V_{\{\bar
z_j\}}$ and $\bar\theta_j:=\bar z_j/(4|s|)$, which lies in $[-C_\theta,C_\theta]$ by
(iii). The lift $\mathcal C_j\times\{\bar z_j\}$ is Legendrian because $\overline\beta$
vanishes on Lagrangian cones with vertex $0$: the position vector $w$ is tangent, and
$\overline\beta(U)=\frac12\omega_0(w,U)$ vanishes for tangent $U$, since the cone is
Lagrangian. By (iii), the Legendrian angle on the $j$-th slice is $\bar\theta_j$.

\emph{(v).} Fix $s_0$ for which (i)--(iv) hold. Then $\mathfrak
m_{s_0}=\sum_j\big(\|\mathcal
C_j\|\times\delta_{4|s_0|\bar\theta_j}\big)\otimes\delta_{\bar\theta_j}$ by (iii)--(iv).
By Proposition~\ref{prop:transport}, $\mathfrak m_s=(\mathrm{Sh}_{s-s_0})_\#\mathfrak
m_{s_0}$ for every $s$. The map $\mathrm{Sh}_{s-s_0}$ moves the slice at height
$4|s_0|\bar\theta_j$, with Legendrian angle $\bar\theta_j$, to height
$4|s_0|\bar\theta_j-4(s-s_0)\bar\theta_j=4|s|\bar\theta_j$, and does not change $w$. This
is the formula for $\mathfrak m_s$; $\mu_s$ is its first marginal. Applying this to two
admissible values $s_0,s_1$ shows that the data $\bar\theta_j,\mathcal C_j$ do not depend
on $s_0$. Finally, for $s_0$ as above each lifted cone $\mathcal C_j\times\{\bar z_j\}$ is
the varifold $V_{\{\bar z_j\}}$, which is stationary by the beginning of the proof of
(iv), applied to a small interval around $\bar z_j$. For an arbitrary $s$ the lifted cones
are obtained from these by the vertical translations $\tau_{4(s_0-s)\bar\theta_j}$, which
are isometries of $\overline g$, so they are stationary, and so is their sum.

\emph{(vi).} Fix $s\in\mathcal G$ for which (i)--(iv) hold; by (v) the cones $\mathcal
C_j$ do not depend on this choice. Fix $j$, and $0<\delta<\frac12\min_{k\ne l}|\bar
z_k-\bar z_l|$ (any $\delta>0$ if $N=1$). Let $B_r\subset\C^n$ denote the $g_0$-balls of
Appendix~\ref{app:A}, $B_r=\{|w|<2r\}$. We first decompose the slice of $V_s$ at height
$\bar z_j$ into calibrated currents, by applying Proposition~\ref{prop:A1} to bands of
$F^i_s$ around this height on balls of every radius and gluing the results; then we show
that each of these currents is a stationary cone. For $R\in\mathbb N$, by (a) and the
coarea formula,
\begin{multline*}
\int_{\delta/2}^{\delta}\mathcal H^{n-1}\big(\{|z\circ F^i_s-\bar
z_j|=\delta'\}\cap(\pi\circ F^i_s)^{-1}(B_{4R})\big)\,d\delta'\\
\le\sup|\nabla z|\ \mu^i_s\big(\{\tfrac\delta2\le|z-\bar
z_j|\le\delta\}\cap\pi^{-1}(\overline B_{4R})\big).
\end{multline*}
Here $\nabla z$ is the gradient of $z\circ F^i_s$, and $\sup|\nabla z|\le C(R)$ by
Lemma~\ref{lem:model}(d). The right-hand side tends to $0$, because the limit carries no
mass in the closed set on the right. Pass to a subsequence along which these right-hand
sides are summable in $i$, for every $R\in\mathbb N$. Then for
a.e.\ $\delta'\in(\frac\delta2,\delta)$ both $\bar z_j\pm\delta'$ are regular values of
every $z\circ F^i_s$ (Sard), and for every $R\in\mathbb N$ the boundary measures on the
left tend to $0$ as $i\to\infty$ (Fubini). Fix such a $\delta'$, and put
\[
L^i:=\{x:|z(F^i_s(x))-\bar z_j|\le\delta'\},\qquad L^i_R:=L^i\cap(\pi\circ
F^i_s)^{-1}(B_{4R}).
\]
$L^i$ is a manifold with boundary, oriented by (O).

Let $\vartheta_i\colon L^i\to\R/2\pi\Z$ be the geometric phase of $\pi\circ F^i_s$,
defined as in Proposition~\ref{prop:A1}, with respect to $\Omega_0$ and the rescaled base
metrics $h_i:=h^{(\lambda_i)}_p$ of the proof of Lemma~\ref{lem:moving}(b), so that
$g_i=\pi^*h_i+\eta_i\otimes\eta_i$ in the rescaled chart, and put $\hat{\mathfrak
m}^i:=(\pi\circ F^i_s,\vartheta_i)_\#(\mu\llcorner L^i)$. By (a), $\hat{\mathfrak
m}^i(B_{4R}\times\R/2\pi\Z)\le D_0(8R)^n$. Pass to a further subsequence with
$\hat{\mathfrak m}^i\to\hat{\mathfrak m}$ as Radon measures on $\C^n\times\R/2\pi\Z$.

Let $R\in\mathbb N$. For $i$ large the map $\pi\circ F^i_s\colon L^i_R\to B_{4R}$ is
proper: the preimage of a compact set $\mathcal K\subset B_{4R}$ is the preimage under
$F^i_s$ of the compact set $\pi^{-1}(\mathcal K)\cap\{|z-\bar z_j|\le\delta'\}$, which
lies in the rescaled chart. It is a Lagrangian immersion for the metrics $h_i$, which
converge to $g_0$ in $C^2(B_{4R})$ (Lemma~\ref{lem:moving}(b)). By
\cite[Lemma~2.3]{ChangWuZhang2026ContactBlowup}, whose proof is pointwise and applies in
the chart, $\pi\circ F^i_s$ is an isometric immersion for $h_i$ with mean curvature vector
$\pi_*H^i$. Hence its $h_i$-area measure is $\mu$, so that $\mu(L^i_R)\le8^nD_0R^n$; its
mean curvature has norm $|H^i|$, and $\int_{L^i_R}|H^i|^2\to0$ since $s\in\mathcal G$; and
$\mathcal H^{n-1}(\partial L^i_R)\to0$. Hence Proposition~\ref{prop:A1} applies, with
$h_0=g_0$. In its proof the Young measure is the restriction of $\hat{\mathfrak m}$ to
$B_{4R}\times\R/2\pi\Z$, because the $\hat{\mathfrak m}^i$ already converge; the further
subsequences taken there do not change it. Hence the second marginal $\mathfrak
n_R:=(\mathrm{pr}_2)_\#\big(\hat{\mathfrak m}\llcorner(B_R\times\R/2\pi\Z)\big)$ is a
finite sum of at most $N_0=N_0(n,8^nD_0)$ point masses. Let $\Lambda_R:=\spt\mathfrak
n_R$, a set with at most $N_0$ points, each of which carries positive $\mathfrak
n_R$-mass; $\Lambda_R$ is empty if and only if $\hat{\mathfrak m}(B_R\times\R/2\pi\Z)=0$.
Then
\[
\hat{\mathfrak m}\llcorner(B_R\times\R/2\pi\Z)=\sum_{\bar\vartheta\in\Lambda_R}\|\mathcal
T_{\bar\vartheta,R}\|\otimes\delta_{\bar\vartheta},
\]
where the $\mathcal T_{\bar\vartheta,R}$ are nonzero integral currents in $B_R$ with
$\partial \mathcal T_{\bar\vartheta,R}=0$ in $B_R$, calibrated by
$\mathrm{Re}(e^{-i\bar\vartheta}\Omega_0)$; they are nonzero because $\|\mathcal
T_{\bar\vartheta,R}\|(B_R)=\mathfrak n_R(\{\bar\vartheta\})>0$.

For $R<R'$ we have $\mathfrak n_R\le\mathfrak n_{R'}$, hence
$\Lambda_R\subset\Lambda_{R'}$: the sets $\Lambda_R$ increase with $R$ and have at most
$N_0$ points, so $\Lambda:=\bigcup_R\Lambda_R$ has at most $N_0$ points, and every
$\bar\vartheta\in\Lambda$ carries positive mass in some $B_R$. An integer multiplicity
current $\mathcal T$ calibrated by a form of comass $1$ is determined by its mass measure
$\|\mathcal T\|$. Indeed, $\|\mathcal T\|$ determines the rectifiable set carrying
$\mathcal T$, the multiplicity and the approximate tangent planes, and the orientation of
each plane is the one on which the form equals $1$. For $\bar\vartheta\in\Lambda$ put
$\mathcal T_{\bar\vartheta,R}:=0$ if $\bar\vartheta\notin\Lambda_R$, and
\[
\mu^{(\bar\vartheta)}:=(\mathrm{pr}_1)_\#\big(\hat{\mathfrak
m}\llcorner(\C^n\times\{\bar\vartheta\})\big).
\]
For $R<R'$, the currents $\mathcal T_{\bar\vartheta,R}$ and $\mathcal
T_{\bar\vartheta,R'}\llcorner B_R$ are calibrated by the same form and both have mass
measure $\mu^{(\bar\vartheta)}\llcorner B_R$. Hence they are equal. Thus the $\mathcal
T_{\bar\vartheta,R}$ define an integer multiplicity current $\mathcal T_{\bar\vartheta}$
in $\C^n$ with locally finite mass, $\partial \mathcal T_{\bar\vartheta}=0$, calibrated by
$\mathrm{Re}(e^{-i\bar\vartheta}\Omega_0)$, and with $\|\mathcal
T_{\bar\vartheta}\|=\mu^{(\bar\vartheta)}\ne0$.

Since the circle is compact,
\[
\sum_{\bar\vartheta\in\Lambda}\|\mathcal
T_{\bar\vartheta}\|=(\mathrm{pr}_1)_\#\hat{\mathfrak m}=\lim_i(\pi\circ
F^i_s)_\#(\mu\llcorner L^i).
\]
The limit $\mu_s$ has no mass on $\{|z-\bar z_j|=\delta'\}$, and $\pi$ is proper on
$\{|z-\bar z_j|\le\delta'\}$. Hence the right-hand side is
$\pi_\#\big(\mu_s\llcorner\{|z-\bar
z_j|\le\delta'\}\big)=\pi_\#\big(\mu_s\llcorner\{z=\bar z_j\}\big)=\|\mathcal C_j\|$. Thus
$\mathcal C_j$ and $\sum_{\bar\vartheta}|\mathcal T_{\bar\vartheta}|$ are rectifiable
varifolds with the same weight, hence equal. (The sum is of varifolds, i.e.\ of the masses
$\|\mathcal T_{\bar\vartheta}\|$; currents are not added, so no cancellation between
opposite orientations can occur. When two phases in $\Lambda$ differ by $\pi$, the
varifold $\mathcal C_j$ alone does not determine the currents $\mathcal
T_{\bar\vartheta}$; they are determined by $\hat{\mathfrak m}$.)

We now show that each $\mathcal T_{\bar\vartheta}$ is a stationary cone. Fix
$\bar\vartheta\in\Lambda$, and put $\mathcal T:=\mathcal T_{\bar\vartheta}$ and
$\Omega^{\rm Re}:=\mathrm{Re}(e^{-i\bar\vartheta}\Omega_0)$. The first variation
$\delta|\mathcal T|(Y)$ is the integral of $\mathrm{div}_SY$ against $|\mathcal T|$, so it
depends continuously on $Y$ in the $C^1$ topology when $\spt Y$ stays in a fixed compact
set. Hence, to see that $|\mathcal T|$ is stationary, it suffices to show $\delta|\mathcal
T|(Y)=0$ for $Y\in C^\infty_c(\C^n,\C^n)$. Let $\mathrm{Fl}^Y_t$ be the flow of such a
$Y$, choose $r_Y$ with $\spt Y\subset B_{r_Y}$, and let $\gamma\in C^\infty_c(\C^n)$ with
$0\le\gamma\le1$ and $\gamma=1$ near $\overline B_{r_Y}$. Since $\partial\mathcal T=0$,
the boundary of $\mathcal T\llcorner B_{r_Y}$ is supported in $\partial B_{r_Y}$. The flow
$\mathrm{Fl}^Y_t$ is the identity outside $\spt Y$ and maps $\spt Y$ onto itself. Hence
the affine homotopy $(\varsigma,x)\mapsto(1-\varsigma)x+\varsigma\,\mathrm{Fl}^Y_t(x)$,
$\varsigma\in[0,1]$, is smooth, maps the convex set $B_{r_Y}$ into itself, and does not
move the points outside $\spt Y$, in particular those near $\partial B_{r_Y}$; and
$\mathrm{Fl}^Y_t$ fixes $\mathcal T\llcorner(\C^n\setminus B_{r_Y})$. The homotopy formula
\cite[26.22]{si83}, applied to $\mathcal T\llcorner B_{r_Y}$ and this homotopy, therefore
has no boundary term and gives $(\mathrm{Fl}^Y_t)_\#\mathcal T-\mathcal T=\partial\Xi_t$
with $\spt\Xi_t\subset\overline B_{r_Y}$. As $d\Omega^{\rm Re}=0$ and $d\gamma=0$ near
$\spt\Xi_t$, we have $(\partial\Xi_t)(\gamma\,\Omega^{\rm
Re})=\Xi_t(d\gamma\wedge\Omega^{\rm Re})=0$. Since $\mathcal T$ is calibrated by
$\Omega^{\rm Re}$, which has comass $1$,
\[
\int\gamma\,d\|\mathcal T\|=\mathcal T(\gamma\,\Omega^{\rm
Re})=(\mathrm{Fl}^Y_t)_\#\mathcal T(\gamma\,\Omega^{\rm
Re})\le\int\gamma\,d\|(\mathrm{Fl}^Y_t)_\#\mathcal T\| ,
\]
as in \cite[Ch.~II, Theorem~4.2]{hl82}. Thus
$t\mapsto\int\gamma\,d\|(\mathrm{Fl}^Y_t)_\#\mathcal T\|$ is minimal at $t=0$, and its
derivative there, which is $\delta|\mathcal T|(Y)$ because $\gamma=1$ near $\spt Y$,
vanishes. Next, since $\|\mathcal T\|\le\|\mathcal C_j\|$, the approximate tangent planes
of $\mathcal T$ and of $\mathcal C_j$ agree $\|\mathcal T\|$-a.e., and $\mathfrak r$ is
tangent to $\mathcal C_j$ a.e., as shown in the proof of (iv). The computation
\eqref{eq:cone-homog} therefore applies to $|\mathcal T|$ and gives $\|\mathcal
T\|(\lambda A)=\lambda^n\|\mathcal T\|(A)$ for all $\lambda>0$. The push-forward of
$\mathcal T$ under $w\mapsto\lambda w$ is calibrated by $\Omega^{\rm Re}$, which has
constant coefficients, and it has the same mass measure as $\mathcal T$. Hence it equals
$\mathcal T$. Finally write $\Lambda=\{\bar\vartheta_{j1},\dots,\bar\vartheta_{jN_j}\}$,
with $N_j\le N_0$, and $\mathcal T_{jk}:=\mathcal T_{\bar\vartheta_{jk}}$.
\end{proof}

\emph{All negative times.} Theorem~\ref{thm:weak} refers to a fixed $S_0>0$. Weak contact
tangent flows for all negative times, in the sense of Definition~\ref{def:weak}, are
obtained by a diagonal argument. Apply Proposition~\ref{prop:compactness}(c) with
$S_0=1,2,\dots$ in turn, each time to the subsequence produced at the previous step (for
$S_0=k$ the rescalings are defined for all $i$ with $k\lambda_i^{-2}\le\tau_0$), and let
$(\lambda_{i_l})$ be the diagonal subsequence. Along it, $\mu^i_s\to\mu_s$ and
$\fm^i_s\to\fm_s$ for every $s<0$. For each $k$, from its $k$-th term on,
$(\lambda_{i_l})$ is a subsequence of the $k$-th subsequence, and the conclusions of
Proposition~\ref{prop:compactness}(c) on $[-k,0)$, with its set $\cG_k$ of good times,
persist along every further subsequence. Hence the restriction of $(\mu_s,\fm_s)$ to
$[-k,0)$ is a weak contact tangent flow in the sense of Definition~\ref{def:weak} with
$S_0=k$, to which Theorem~\ref{thm:weak} applies.

The values $\bar\theta_j$ and the cones $\cC_j$ obtained on $[-k,0)$ do not depend on $k$.
Indeed, by (v), for any fixed $s<0$ the numbers $4|s|\bar\theta_j$ are the heights of the
level sets $\{z=\mathrm{const}\}$ charged by $\mu_s$, and $\|\cC_j\|$ is the image of
$\mu_s\llcorner\{z=4|s|\bar\theta_j\}$ under the projection $(w,z)\mapsto w$, which
determines the integral varifold $\cC_j$. Hence the conclusions of
Theorem~\ref{thm:weak}(v) hold for every $s<0$ with the same $\bar\theta_j$ and $\cC_j$,
and (vi) concerns these cones. Theorems~\ref{thm:B} and~\ref{thm:D}(a),
Corollary~\ref{cor:heights} and Theorem~\ref{thm:curves} are stated for such flows. The
same argument, started from any subsequence along which a weak contact tangent flow on
$[-S_0,0)$ is given, extends that flow to all negative times along a further subsequence,
without changing it on $[-S_0,0)$, nor the $\bar\theta_j$ and $\cC_j$.

\subsection{Weighted densities and their consequences}\label{subsec:weighted}

In the next proposition, for a sequence of flows $F^i$, $\mu^i_t$ denotes the area measure
of $F^i$ at the original time $t$ for $g$, not the push-forward of
Definition~\ref{def:weak}.
\begin{proposition}\label{prop:moving-densities}
Let $F^i$ be flows satisfying (S), (F) with the same constants, including the constant
$D_0$ of (D). Let $p_i\in M$, $T_i$, $\lambda_i\to\infty$, and let $\tilde
F^i_s,\tilde{\mathfrak m}^i_s$ be the rescalings at $(p_i,T_i)$. Fix $s$ and assume
$\tilde{\mathfrak m}^i_s\to\mathfrak m_s$ weakly. Let $(q_i,t_i)$ be centres, not
necessarily on the flow and possibly with $t_i\ge T_i$, with
\[
\hat y_i:=\delta_{\lambda_i}\Upsilon_{p_i}(q_i)\to\hat y\in\Hn,\qquad \tilde
t_i:=\lambda_i^2(t_i-T_i)\to\tilde t>s,
\]
and put $\sigma_i:=\lambda_i^{-2}(\tilde t_i-s)$ and $\tilde\sigma:=\tilde t-s>0$. Let
\[
\overline\rho^\chi_{\hat y,\tilde t}(x,\theta)=(4\pi\tilde\sigma)^{-n/2}e^{-f_{\hat
y}(x)/\tilde\sigma}\chi\big(z_{\hat y}(x)/\tilde\sigma-4\theta\big)
\]
be the model kernel. Then for every $\varpi\in C(\R)$
\[
\lim_i\int\psi_{q_i}\,\varpi(\theta)\,\rho^\chi_{q_i,t_i}\,d\mu^i_{t_i-\sigma_i}=\int
\overline\rho^\chi_{\hat y,\tilde t}\,\varpi(\theta)\,d\mathfrak m_s .
\]
The convergence is uniform in the centres: for a compact set $\mathcal Y\subset\Hn$ and
$[\sigma_1,\sigma_2]\subset(0,\infty)$, the supremum of the difference of the two sides,
over all centres with $\hat y_i\in\mathcal Y$ and $\tilde t_i-s\in[\sigma_1,\sigma_2]$,
with the right side taken at $(\hat y_i,\tilde t_i)$, tends to $0$.
\end{proposition}
\begin{proof}
After rescaling, the kernels converge uniformly on a common compact set, but the
integrals are taken over the whole flow and the two charts differ; so we also need tail
bounds, uniform in $i$, outside a common compact set, and a choice of this set whose
boundary carries no mass in the limit. Uniformity in the centres then follows by
compactness of the parameters. The integrand on the left vanishes outside $U_{q_i}$. We
split the integral as
$I_i'+I_i''$, where $I_i'$ is the integral over the part of the flow mapped into
$U_{p_i}$, the domain of $\Upsilon_{p_i}$, and $I_i''$ the integral over the part mapped
into $U_{q_i}\setminus U_{p_i}$, and we rescale $I_i'$ by
$\delta_{\lambda_i}\Upsilon_{p_i}$ and $t=T_i+\lambda_i^{-2}s'$. Put $\tilde
f_i:=f_{q_i}^{(\lambda_i)}$, $\tilde z_i:=z_{q_i}^{(\lambda_i)}$ and
$\tilde\psi_i:=\psi_{q_i}^{(\lambda_i)}$ in the notation of Lemma~\ref{lem:moving}(c),
that is, $\tilde f_i=\lambda_i^2f_{q_i}$ and $\tilde z_i=\lambda_i^2z_{q_i}$ written in
the rescaled $p_i$-chart, while $\tilde\psi_i$ is $\psi_{q_i}$ so written. With
$\tilde\sigma_i:=\lambda_i^2\sigma_i=\tilde t_i-s$, and with
$d\mu'=\lambda_i^nd\mu^i_{t_i-\sigma_i}$ the rescaled area measure, the kernel transforms
as in the proof of Lemma~\ref{lem:scaling}, which is a pointwise computation that does not
use that the centre of the kernel is the centre of the rescaling:
\[
(4\pi\sigma_i)^{-n/2}e^{-f_{q_i}/\sigma_i}\chi\Big(\frac{z_{q_i}}{\sigma_i}-\frac{a_K(\sigma_i)}{\sigma_i}\theta\Big)d\mu^i_{t_i-\sigma_i}=(4\pi\tilde\sigma_i)^{-n/2}e^{-\tilde
f_i/\tilde\sigma_i}\chi\Big(\frac{\tilde
z_i}{\tilde\sigma_i}-\frac{a_K(\sigma_i)}{\sigma_i}\theta\Big)d\mu' .
\]
Hence $I_i'=\int\tilde\psi_i\,\varpi(\theta)\,\tilde\rho^\chi_i\,d\tilde{\mathfrak
m}^i_s$, with
\[
\tilde\rho^\chi_i=(4\pi\tilde\sigma_i)^{-n/2}e^{-\tilde f_i/\tilde\sigma_i}\chi\big(\tilde
z_i/\tilde\sigma_i-\hat a_i\theta\big),\qquad \tilde\sigma_i=\tilde t_i-s,\quad \hat
a_i=a_K(\sigma_i)/\sigma_i\to4;
\]
outside the image of $U_{q_i}$ we put $\tilde\psi_i\tilde\rho^\chi_i:=0$, in accordance
with $\psi_{q_i}=0$ there. The remainder $I_i''$ is a tail term and is estimated below. By
Lemma~\ref{lem:moving}(c), $\tilde f_i\to f_{\hat y}$, $\tilde z_i\to z_{\hat y}$ and
$\tilde\psi_i\to1$ uniformly on compact sets; the unitary ambiguity there does not affect
these functions. Hence
\begin{equation}\label{eq:kernel-local}
\tilde\psi_i\varpi(\theta)\tilde\rho^\chi_i\longrightarrow\varpi(\theta)\overline\rho^\chi_{\hat
y,\tilde t}\quad\text{uniformly on }\mathcal K\times[-C_\theta,C_\theta]
\end{equation}
for every compact $\mathcal K$, uniformly in $\hat y$ in compact sets.

For the tails we work in the rescaled $q_i$-chart, where $f=|w|^2/16$ exactly and (D)
holds at $q_i$. Given $\eps>0$, first choose $R$ and then $Z=Z(R)$ by
Lemma~\ref{lem:tails} at the centre $(q_i,t_i)$ with radius $R/\lambda_i$, height
$Z/\lambda_i^2$ and $\sigma=\sigma_i=\tilde\sigma_i/\lambda_i^2$ (note
$\tail(R/\lambda_i,Z/\lambda_i^2;\sigma_i)=\tail(R,Z;\tilde\sigma_i)$), such that the
kernel mass outside
\[
\mathcal K_{R,Z}:=\{\lambda_i^2f_{q_i}\le R^2,\ \lambda_i^2|z_{q_i}|\le Z\}
\]
is at most $\eps$, uniformly in $i$ and in $\tilde\sigma_i$ in a compact subset of
$(0,\infty)$. The lemma applies also if $t_i\ge T_i$, since the evaluation time
$t_i-\sigma_i=T_i+\lambda_i^{-2}s$ is a time of the flow. For $i$ large, $\mathcal
K_{R,Z}$ lies in the common domain of the two charts, and in $p_i$-coordinates it lies in
a fixed compact set $\mathcal K''$. Indeed, by Lemma~\ref{lem:moving}(c) the transition is
defined at the origin, the image of $p_i$, and maps it to a bounded point, so that
$p_i\in U_{q_i}$ and $\delta_{\lambda_i}\Upsilon_{q_i}(p_i)$ is bounded.
Lemma~\ref{lem:moving}(c) with $p_i$ and $q_i$ interchanged then gives these inclusions.
Put $\Phi:=\max\{f_{\hat y}-R^2,\,|z_{\hat y}|-Z\}$,
so that $A:=\{\Phi\le0\}$ is the model set. The uniform convergence $\tilde f_i\to f_{\hat
y}$, $\tilde z_i\to z_{\hat y}$ on $\mathcal K''$ gives, for every $\delta>0$ and $i\ge
i(\delta)$,
\begin{equation}\label{eq:KRZ-incl}
\{\Phi\le-\delta\}\subset\mathcal K_{R,Z}\subset\{\Phi\le\delta\}\qquad\text{in
$p_i$-coordinates.}
\end{equation}
Points of the $q_i$-chart outside $\mathcal K_{R,Z}$ are part of the tail. They include
all points outside the $p_i$-chart, so
$|I_i''|\le\|\varpi\|_{C^0[-C_\theta,C_\theta]}\eps$ for $i$ large; here and below only
the values of $\varpi$ on $[-C_\theta,C_\theta]$ enter, because $|\theta^i|\le C_\theta$,
so that $\tilde{\mathfrak m}^i_s$ and $\mathfrak m_s$ are supported in
$\Hn\times[-C_\theta,C_\theta]$. The limit $\mathfrak m_s$ satisfies the same tail bound
in the model, because mass on open sets is lower semicontinuous.

The set $\{\Phi=0\}$, which contains the boundary of $A$, lies in $\{f_{\hat
y}=R^2,\ |z_{\hat y}|\le Z\}\cup\{f_{\hat y}\le R^2,\ |z_{\hat y}|=Z\}$. For each integer
$Z_1$ the
sets $\{f_{\hat y}=R^2,\ |z_{\hat y}|\le Z_1\}$, $R>0$, are disjoint and $\mathfrak m_s$
is finite on bounded sets, so at most countably many of them carry mass; likewise, for
each integer $R_1$ at most countably many of the sets $\{f_{\hat y}\le R_1^2,\ |z_{\hat
y}|=Z\}$, $Z>0$, carry mass. Thus there is a countable set of radii, and a countable set
of heights, outside which $\{\Phi=0\}$ is $\mathfrak m_s$-null. The tail bound above holds
for all $R\ge R(\eps)$ and then all $Z\ge Z(R)$, so we may choose $R$ and then $Z$ outside
these countable sets. Let $g_i:=\tilde\psi_i\varpi\tilde\rho^\chi_i$ and
$g:=\varpi\overline\rho^\chi_{\hat y,\tilde t}$. By \eqref{eq:KRZ-incl}, the integrals of
$g_i$ over $\mathcal K_{R,Z}$ and over $A$ against $\tilde{\mathfrak m}^i_s$ differ by at
most $\sup_i\sup|g_i|\,\tilde{\mathfrak m}^i_s(\{|\Phi|\le\delta\})$, and
$\limsup_i\tilde{\mathfrak m}^i_s(\{|\Phi|\le\delta\})\le\mathfrak
m_s(\{|\Phi|\le\delta\})\to\mathfrak m_s(\{\Phi=0\})=0$ as $\delta\to0$. On the fixed set
$A$, \eqref{eq:kernel-local} and weak convergence, with $\mathfrak m_s(\partial
A\times\R)=0$, give convergence of the integrals of $g_i$ to that of $g$. Outside
$\mathcal K_{R,Z}$, including the remainder $I_i''$, and outside $A$ in the model, the
tail bound gives $|\cdot|\le\|\varpi\|_{C^0[-C_\theta,C_\theta]}\eps$ for both sides. This
proves the convergence.

For the uniformity in the centres, note first that the right-hand side is continuous in
$(\hat y,\tilde\sigma)\in\Hn\times(0,\infty)$: by the tail bound, applied to centres with
$\hat y_i=\hat y$, the part of it coming from the complement of $\{f_{\hat y}\le
R^2,\ |z_{\hat y}|\le Z\}$ is at most $\|\varpi\|_{C^0[-C_\theta,C_\theta]}\eps$ for all
$\hat
y$ and all $\tilde\sigma$ in a compact subset of $(0,\infty)$, since the bound of
Lemma~\ref{lem:tails} does not depend on the centre; and on compact sets
$\overline\rho^\chi_{\hat y,\tilde t}(x,\theta)$ is uniformly continuous in $(\hat
y,\tilde\sigma,x,\theta)$. If the convergence were not uniform, there would be a
subsequence and centres with $\hat y_i\in\mathcal Y$ and $\tilde
t_i-s\in[\sigma_1,\sigma_2]$ for which the difference of the two sides stays bounded away
from $0$. A further subsequence has $\hat y_i\to\hat y$ and $\tilde t_i\to\tilde t$, and
the convergence just proved, together with the continuity of the right-hand side, gives a
contradiction.
\end{proof}

\begin{corollary}\label{cor:angle-measure}
Assume (S), (F), and let $p\in M$.
\begin{enumerate}[(a)]
\item For every $\varpi\in C^2(\R)$ the limit $\cL_{p,T}(\varpi):=\lim_{t\to
T}\int\psi_p\,\varpi(\theta)\,\rho^\chi_{p,T}\,d\mu_t$ exists.
\item For every weak contact tangent flow at $(p,T)$, written as in
Theorem~\ref{thm:weak}(v), $\cL_{p,T}(\varpi)=\sum_jm_j\,\varpi(\bar\theta_j)$, where
$m_j\ge1$ is the density of $\cC_j$ at its vertex.
\end{enumerate}
Consequently the measure $\nu_{p,T}:=\sum_jm_j\delta_{\bar\theta_j}$ on $\R$, hence the
set $\{\bar\theta_1,\dots,\bar\theta_N\}$ of Legendrian angles and the vertex densities
carried by them, is the same for all weak contact tangent flows at $(p,T)$, whatever the
sequence of scales.
\end{corollary}
\begin{proof}
\emph{(a).} Let $\varpi_+\in C^2(\R)$ satisfy $\varpi_+''\ge0$ and $\varpi_+\ge1$ on
$I_\theta:=[-C_\theta,C_\theta]$; below, $\varpi_+$, $\varpi_+'$ and $\varpi_+''$ are
evaluated at $\theta\in I_\theta$, so only the values on $I_\theta$ matter, and we write
$M_+:=\max_{I_\theta}\varpi_+$. As in the proof of Proposition~\ref{prop:theta2},
\[
\frac{d}{dt}\int\psi\,\varpi_+\chi(v)\,\rho=\int(\partial_t-\Delta)(\psi\varpi_+\chi)\,\rho+\int\psi\varpi_+\chi\,(\partial_t+\Delta-|H|^2)\rho
,
\]
and, by (F) and $|\nabla\theta|=|H|$,
\begin{multline*}
(\partial_t-\Delta)(\psi\varpi_+\chi)=\varpi_+(\partial_t-\Delta)(\psi\chi)+\psi\chi\big((K+2)\theta\varpi_+'-\varpi_+''|H|^2\big)\\
-2\chi\varpi_+'\langle\nabla\theta,\nabla\psi\rangle-2\psi\varpi_+'\chi'\langle\nabla\theta,\nabla
v\rangle .
\end{multline*}
The terms with the factor $\varpi_+$ are those of the proof of
Proposition~\ref{prop:mono}. There, the principal terms combine pointwise to $-\psi
D^\chi\rho$; each remainder $R_1$--$R_5$ is bounded by a pointwise inequality for its
integrand, which is then integrated; and the absorption of the two terms
$\frac{c_\chi}8\psi\chi|\mathcal S_p|^2\rho$ into $\frac14\psi D^\chi\rho$ is pointwise as
well. Since $\varpi_+\in[1,M_+]$ on $I_\theta$, these estimates survive multiplication by
$\varpi_+$, and give
\[
\int\varpi_+(\partial_t-\Delta)(\psi\chi)\,\rho+\int\psi\varpi_+\chi\,(\partial_t+\Delta-|H|^2)\rho\le-\frac34\int\psi\varpi_+D^\chi\rho+M_+\,E(t),
\]
with $E$ as in \eqref{eq:mono}. Of the new terms, $\psi\chi(K+2)\theta\varpi_+'$ is
bounded, $-\psi\chi\varpi_+''|H|^2\le0$, and the cut-off term
$-2\chi\varpi_+'\langle\nabla\theta,\nabla\psi\rangle$ is bounded by
$\psi|H|^2+C|\nabla\zeta|^2$. The cross term
$-2\psi\varpi_+'\chi'\langle\nabla\theta,\nabla v\rangle$ needs more care. Here $|\nabla
v|=4|\mathcal S_p+e_u|$ by \eqref{eq:grad-u}, and
$|\chi'|\le(4a/e)^{1/2}\chi^{1/2}$ (Section~\ref{sec:mono}), so
\[
8|\varpi_+'||\chi'||H||\mathcal S_p+e_u|\le\tfrac{c_\chi}4\,\chi|\mathcal
S_p|^2+C\big(|H|^2+|e_u|^2\big).
\]
The first term is absorbed by $\frac34\psi\varpi_+D^\chi\ge\frac34c_\chi\psi\chi|\mathcal
S_p|^2$, since $\varpi_+\ge1$, and $\int\psi|e_u|^2\rho\le C\big(1+\int\psi|H|^2\rho\big)$
by \eqref{eq:eu-bound}. Hence $\frac{d}{dt}\int\psi\varpi_+\rho^\chi\le
C\big(1+\int\psi|H|^2\rho\big)$, which is integrable on $[T-\tau_0,T)$ by
Proposition~\ref{prop:theta2}. Hence
$\int\psi\varpi_+(\theta)\rho^\chi-\int_{T-\tau_0}^tC(1+\int\psi|H|^2\rho)$ is
non-increasing and bounded below, and the limit exists. Every $\varpi\in C^2(\R)$ is a
difference of two such functions: $\varpi=(\varpi+c_1\theta^2+c_2)-(c_1\theta^2+c_2)$ with
$c_1\ge\frac12\max_{I_\theta}|\varpi''|$ and $c_2\ge1+\max_{I_\theta}|\varpi|$, since then
both functions have nonnegative second derivative and are at least $1$ on $I_\theta$.

\emph{(b).} By Proposition~\ref{prop:moving-densities} with $q_i=p$ and $t_i=T$,
$\cL_{p,T}(\varpi)=\int \overline\rho^\chi_{0,0}\,\varpi(\theta)\,d\mathfrak m_s$ for
every $s\in[-S_0,0)$: with $q_i=p$, $t_i=T$ the left-hand side of that proposition is the
weighted integral at the times $T-\lambda_i^{-2}|s|\to T$, which converges to
$\cL_{p,T}(\varpi)$ by (a). By Theorem~\ref{thm:weak}(v) and the paragraph \emph{All
negative times} following its proof, on the $j$-th lifted cone the height is
$4|s|\bar\theta_j$ and the Legendrian angle is $\bar\theta_j$, so $v=0$ and
$\overline\rho^\chi_{0,0}=\rho'$ there. Hence
$\cL_{p,T}(\varpi)=\sum_j\varpi(\bar\theta_j)\int\rho'\,d\|\mathcal C_j\|$, and the
Gaussian integral of a cone is its vertex density, which is $\ge1$ (proof of
Theorem~\ref{thm:weak}(iv)).
\end{proof}

\begin{corollary}\label{cor:heights}
In a weak contact tangent flow for all negative times, for every $s<0$ the lifted cones
$\cC_j\times\{4|s|\bar\theta_j\}$ with different Legendrian angles lie at the mutual
height distance $4|s||\bar\theta_j-\bar\theta_k|$. If two of them have a common projected
point, the vertical segment joining them is a Reeb chord of the model of
$\overline\eta$-action $2|s||\bar\theta_j-\bar\theta_k|$, which tends to $0$ as $s\to0$.
The rate is the one of the flow: if $t\mapsto\gamma_t$ is a smooth family of Reeb chords
of $\Sigma_t$ with endpoints $F_t(x_\pm(t))$,
$F_t(x_+(t))=\mathrm{Fl}^\xi_{A(t)}(F_t(x_-(t)))$, then its $\eta$-action $A(t)$ satisfies
\begin{equation}\label{eq:chord}
A'(t)=2\big(\theta(x_-(t),t)-\theta(x_+(t),t)\big).
\end{equation}
\end{corollary}
\begin{proof}
The first two assertions are Theorem~\ref{thm:weak}(v) and the fact that a Reeb segment
has $\overline\eta$-length equal to half its $z$-length. The projections of any two lifted
cones always meet: each nonzero cone $\cC_j$ is dilation invariant with closed support, so
$0\in\spt\cC_j$, and the vertical segment over the common vertex $w=0$ is such a Reeb
chord. For \eqref{eq:chord}, differentiate
$F_t(x_+(t))=\mathrm{Fl}^\xi_{A(t)}(F_t(x_-(t)))$ and apply $\eta$: since the Reeb flow
preserves $\eta$, and $\eta$ vanishes on $dF_t(T\Sigma)$, one gets
$\eta(\partial_tF)(x_+)=\eta(\partial_tF)(x_-)+A'(t)$, and $\eta(\partial_tF)=-2\theta$ by
(F). On the limit the lower lifted cone has the Legendrian angle $\bar\theta_j$ and the
upper one $\bar\theta_k>\bar\theta_j$, and the action $2|s|(\bar\theta_k-\bar\theta_j)$
has the derivative $2(\bar\theta_j-\bar\theta_k)$ in $s$, in agreement with
\eqref{eq:chord}.
\end{proof}

\section{Applications}\label{sec:applications}

In this section we apply the results of Sections~\ref{sec:mono} and~\ref{sec:weak} in
three situations. In Section~\ref{subsec:typeI} we consider Type~I singularities, the
setting of \cite[Theorem~1.1]{ChangWuZhang2026ContactBlowup}. There the vertical position
of the limits is determined along a suitable sequence of scales by a separate argument,
and the components are identified as planes through the rigidity theorem
\cite[Theorem~A]{cqz} for Legendrian self-shrinkers. The transverse kernel cannot give
more, since it is constant along the Reeb orbits (Remark~\ref{rem:M-compare}(3)). Here the
defect $|v\chi'(v)|/\tau$ of the weighted kernel gives the exact centring along every
sequence (Proposition~\ref{prop:centring}). The $\theta^2$-monotonicity of
Proposition~\ref{prop:theta2} gives $H=0$ on the limit, so that the components are planes
by an elementary argument (Corollary~\ref{cor:planes}). The compactness of the rescalings
is taken from \cite{ChangWuZhang2026ContactBlowup} (Proposition~\ref{prop:C-eta}).

In Section~\ref{subsec:two-planes} we describe the flow near a point at which one weak
contact tangent flow consists of two planes with different Legendrian angles. In
Section~\ref{subsec:curves} we treat Legendrian curves. There the cones of
Theorem~\ref{thm:weak} are lines, and a local component of the rescaled flow has a single
height and a single Legendrian angle in the limit (Theorem~\ref{thm:curves}). This is the
Legendrian counterpart, for curves, of the theorem of Neves \cite[Theorem~B]{nev07},
without his assumptions of almost calibration and rationality.

\subsection{Type~I singularities and exact centring}\label{subsec:typeI}

\begin{lemma}\label{lem:confinement}
Assume (S), (F), (I). Then $F_t\to F_T$ uniformly as $t\to T$. There are $\tau_1>0$ and
$C_3,C_4$, depending on $n,C_I,C_\theta$ and the chart constants, such that for every
$p\in M$, \emph{every} $x\in F_T^{-1}(p)$ and $\tau\le\tau_1$,
\[
|w(F_t(x))|\le C_3\sqrt\tau,\qquad |z_p(F_t(x))|\le C_4\,\tau .
\]
\end{lemma}
\begin{proof}
Since $|\partial_tF|\le|H|+2|\theta|\le\sqrt{nC_I}\,\tau^{-1/2}+2C_\theta$ is integrable,
$F_t$ converges uniformly, and the distance from $F_t(x)$ to $F_T(x)$ is at most
$2\sqrt{nC_I\tau}+2C_\theta\tau$ for every $x$. Hence there is $\tau_1$, independent of
$x$ and $p$, such that $F_{t'}(x)\in U_p$ for $t'\in[T-\tau_1,T]$ whenever $F_T(x)=p$. In
the chart $\Upsilon_p$ the horizontal speed of $F_t(x)$ is $|H|\le\sqrt{nC_I/\tau}$;
integrating from $t$ to $T$, the $h_p$-distance to $w=0$ is at most $2\sqrt{nC_I\tau}$,
hence $|w|\le C_3\sqrt\tau$. For the height, Lemma~\ref{lem:model}(d) gives
$|\frac{d}{dt}z_p(F_t(x))|=|dz_p(H)-4\theta|\le C|w||H|+4C_\theta\le
C\,C_3\sqrt{nC_I}+4C_\theta=:C_4$; here $dz_p(H)=-2\beta_p(H)$. Integrate from $t$ to $T$,
using $z_p(p)=0$. None of the constants depends on $x$ or $p$. For a single point $x$ this
is the non-emptiness step of the proof of
\cite[Theorem~1.1]{ChangWuZhang2026ContactBlowup}.
\end{proof}

The following compactness property of the rescalings is a hypothesis in
Proposition~\ref{prop:centring} and Corollary~\ref{cor:planes}; in the $\eta$-Einstein
case it holds under Type~I (Proposition~\ref{prop:C-eta}).
\begin{enumerate}
\item[(C)] \emph{Compactness hypothesis.} After passing to a subsequence, for every
compact $\mathcal K\subset\Hn$ the parts of $F^i_s$ mapping into $\mathcal K$ converge,
after reparametrisation, smoothly on compact subsets of space-time for $s\in[-2,-1]$ (that
is, the reparametrised maps converge in $C^\infty$ on compact subsets of the domain times
$[-2,-1]$, time derivatives included) to a smooth Legendrian flow $F^\infty$ of complete,
properly immersed submanifolds of $(\Hn,\overline g)$ (possibly with several components
and multiplicities). Under the same reparametrisations, the Legendrian angles converge,
$\theta^i\to\theta^\infty$, locally uniformly. \emph{No vertical recentring is applied.}
Throughout this section, local smooth convergence of flows is understood in this
space-time sense. For such convergence the convergence of the angles is automatic, by
\eqref{eq:angle-velocity} in the proof of Proposition~\ref{prop:centring}, and
$\theta^\infty$ is the Legendrian angle of the limit flow.
\end{enumerate}

\begin{proposition}\label{prop:centring}
Assume (S), (F), (I), let $p\in M$ and $\lambda_i\to\infty$, and assume (C). Then for
every $s\in[-2,-1]$ the limit satisfies the following.
\begin{enumerate}[(a)]
\item $z=4(-s)\,\theta^\infty$ on $F^\infty_s$: every component is \emph{exactly centred}.
\item $\mathcal S^\infty:=H^\infty+(\overline\nabla f)^\perp/(-s)\equiv0$.
\item $H^\infty\equiv0$.
\item By Lemma~\ref{lem:confinement}, the rescaled images of every point $x$ with
$F_T(x)=p$ lie in the fixed compact set $\{|w|\le\sqrt2C_3,\ |z|\le2C_4\}$. In particular,
components containing such points are reached without recentring, and (a)--(c) apply to
them.
\item If, instead of (C), the vertically shifted flows $\tau_{b_i}\circ F^i_s$,
$(w,z)\mapsto(w,z+b_i)$, converge locally smoothly on compact subsets of space-time, after
reparametrisation, to a smooth Legendrian flow $\tilde F^\infty$ for a bounded sequence of
constants $b_i$ (no convergence of the angles being assumed), then $b_i\to b$ along a
subsequence, the angles converge under the same reparametrisations, and the shift
$\tau_{-b}\circ\tilde F^\infty$ satisfies (a)--(c), with $\theta^\infty$ its Legendrian
angle. Thus bounded recentrings are never needed.
\end{enumerate}
\end{proposition}
\begin{proof}
Let $I_i=[T-2\lambda_i^{-2},T-\lambda_i^{-2}]$. As in the proof of
Proposition~\ref{prop:energy}, with $S_0=2$ and $\delta=1$, the monotone quantities of
Propositions~\ref{prop:theta2} and~\ref{prop:mono} give
\[
\int_{I_i}\!\!\int\psi\Big(\tfrac34D^\chi+|H|^2+|\mathcal
S_p|^2\Big)\rho\,d\mu_t\,dt\longrightarrow0 ,
\]
and by Lemma~\ref{lem:scaling} this reads
\begin{equation}\label{eq:centring-diss}
\int_{-2}^{-1}\!\!\int\psi_i\Big(\tfrac34\Big(\frac{|v^i\chi'(v^i)|}{|s|}+(\chi+16\chi'')(v^i)|\mathcal
S^i|^2\Big)+|H^i|^2+|\mathcal S^i|^2\Big)\rho'\,d\mu'\,ds\longrightarrow0 ,
\end{equation}
where $\rho'=(4\pi|s|)^{-n/2}e^{-|w'|^2/(16|s|)}$, $v^i=z'/|s|-\hat a_i(s)\theta^i$ with
$\hat a_i\to4$ uniformly, and the cut-offs $\psi_i\to1$ locally uniformly. By (C) and
Fatou's lemma,
\[
\int_{-2}^{-1}\!\!\int\Big(\frac{|v^\infty\chi'(v^\infty)|}{|s|}+c_\chi\chi(v^\infty)|\mathcal
S^\infty|^2+|H^\infty|^2\Big)\rho'\,d\mu^\infty ds=0,\qquad
v^\infty=\frac{z}{|s|}-4\theta^\infty .
\]
Since $v\chi'(v)<0$ for $v\neq0$, $\chi(0)=1$, and everything is continuous, (a)--(c)
follow for all $s$.

Assertion (d) is Lemma~\ref{lem:confinement} in rescaled variables. For (e), write
$G_i(x,s):=\tau_{b_i}\circ F^i_s(\phi_i(x,s))$ for the reparametrised shifted flows, which
converge to $\tilde F^\infty$ in $C^\infty$ on compact subsets of space-time; the
reparametrisations $\phi_i$ may depend on $s$, the shifts $b_i$ do not. We first recover
the angles from the velocities. Since $\xi_i=2\partial_z$ and $\mathcal
L_{\xi_i}\eta_i=0$, the form $\eta_i$ does not depend on $z$, so
$\tau_{b_i}^*\eta_i=\eta_i$. By the flow equation $\partial_sF^i=H^i-2\theta^i\xi_i$, with
$\eta_i(H^i)=0$ and $\eta_i(\xi_i)=1$, and since $dF^i_s(\partial_s\phi_i)$ is tangent to
the Legendrian immersion $F^i_s$, so that $\eta_i$ annihilates it,
\begin{equation}\label{eq:angle-velocity}
\eta_i(\partial_sG_i)=\eta_i\big(\partial_sF^i_s\circ\phi_i+dF^i_s(\partial_s\phi_i)\big)=-2\,\theta^i_s\circ\phi_i
.
\end{equation}
As $\partial_sG_i\to\partial_s\tilde F^\infty$ locally uniformly and
$\eta_i\to\overline\eta$ in $C^\infty_{\rm loc}$ (Lemma~\ref{lem:moving}(b)), the angles
$\theta^i_s\circ\phi_i$ converge locally uniformly to
$\tilde\theta^\infty:=-\frac12\overline\eta(\partial_s\tilde F^\infty)$, the Legendrian
angle of $\tilde F^\infty$. Now pass to a subsequence with $b_i\to b$. The unshifted
reparametrised flows $\tau_{-b_i}\circ G_i=F^i_s\circ\phi_i$ converge in $C^\infty$ on
compact subsets of space-time to $\tau_{-b}\circ\tilde F^\infty$, and since
$\tau_{-b}^*\overline\eta=\overline\eta$ the Legendrian angle of $\tau_{-b}\circ\tilde
F^\infty$ is again $\tilde\theta^\infty$, the limit of the angles $\theta^i_s\circ\phi_i$.
Thus the unshifted flows satisfy the convergence statements of (C) used in the proof of
(a)--(c), namely local smooth convergence and convergence of the angles under the same
reparametrisations; completeness and properness of the limit were not used there. The
argument for (a)--(c), applied to them, gives (a)--(c) for $\tau_{-b}\circ\tilde
F^\infty$.
\end{proof}

\begin{corollary}\label{cor:planes}
Under the hypotheses of Proposition~\ref{prop:centring}, each connected component of
$F^\infty_s$, $s\in[-2,-1]$, is a horizontal Legendrian plane
\[
\{(w,\,4(-s)\bar\theta_j):\ w\in P_j\}
\]
over a Lagrangian subspace $P_j\subset\C^n$, possibly with multiplicity. Here
$\bar\theta_j\in\R$ is the limit of the Legendrian angle $\theta^i$ on that component. The
limit flow is static horizontally, and each plane translates vertically with speed
$-4\bar\theta_j$, reaching the origin at $s=0$. The geometric phase of $P_j$, with respect
to $\Omega_0$ and an orientation, is a constant $\bar\vartheta_j\in\R/2\pi\Z$ (see
Remark~\ref{rem:potential-phase}).
\end{corollary}
\begin{proof}
Fix $s\in[-2,-1]$ and a connected component of $F^\infty_s$. By (b) and (c),
$(\overline\nabla f)^\perp=0$, so the horizontal dilation field $2\overline\nabla
f=X^{\mathcal H}$, the horizontal lift of $\mathfrak r$, is tangent to each component. The
component is complete and properly immersed, and the flow of this field is complete, so
the projected component is invariant under the dilations $w\mapsto\lambda w$, $\lambda>0$.
By (I), $|B^\infty|\le\sqrt{C_I/|s|}$ is bounded, and by
\cite[Lemma~2.3]{ChangWuZhang2026ContactBlowup} the second fundamental form $B$ of the
projected component has the same norm. Scaling gives $|B|(\lambda w)=\lambda^{-1}|B|(w)$,
so $|B|(w)\le\lambda\sup|B|$; letting $\lambda\to0$ forces $B\equiv0$. Hence the
projection is a complete connected totally geodesic Lagrangian, invariant under dilations,
i.e.\ a linear subspace $P_j$. Its geometric phase is constant, because its tangent plane
is constant and $\Omega_0$ has constant coefficients. On $P_j$ the form
$\overline\beta=\frac12\iota_{\mathfrak r}\omega_0$ vanishes, so the Legendrian lift has a
constant height $\bar z_j$. By (a), $\theta^\infty=\bar z_j/(4|s|)=:\bar\theta_j$ on the
component, and $\theta^i\to\bar\theta_j$ there, by (C).

By (C), for every $s\in[-2,-1]$ the measures $\mathfrak m^i_s=(F^i_s,\theta^i_s)_\#\mu$ of
Definition~\ref{def:weak} converge weakly to $\mathfrak
m_s:=(F^\infty_s,\theta^\infty_s)_\#\mu^\infty$, where $\mu^\infty$ is the $\overline
g$-area measure of $F^\infty_s$ on its parameter domain, as in the proof of
Proposition~\ref{prop:centring}; thus $(F^\infty_s)_\#\mu^\infty$ is the measure of
$F^\infty_s$ on $\Hn$, counted with multiplicity. The proof of
Proposition~\ref{prop:transport} uses neither (I) nor (C), and it applies on $[-2,-1]$; it
gives $\mathfrak m_s=(\mathrm{Sh}_{s-s_0})_\#\mathfrak m_{s_0}$ for $s,s_0\in[-2,-1]$. By
the previous steps,
\[
\mathfrak m_{s_0}=\sum_jm_j\,\big((\mathcal H^n\llcorner
P_j)\times\delta_{4|s_0|\bar\theta_j}\big)\otimes\delta_{\bar\theta_j},
\]
summing over the distinct pairs $(P_j,\bar\theta_j)$, with multiplicities $m_j$ and
$\mathcal H^n$ taken for $g_0$; and $\mathfrak m_s$ has the same form at time $s$. The
terms of distinct pairs are mutually singular, so this form is unique. Since
$\mathrm{Sh}_{s-s_0}$ maps the term of a pair at time $s_0$ to the term of the same pair
at time $s$, the pairs and the multiplicities do not depend on $s$. In the statement,
``reaching the origin at $s=0$'' refers to the height function $4(-s)\bar\theta_j$, which
is linear in $s$ and vanishes at $s=0$; under (C) the flow $F^\infty$ itself is only given
for $s\in[-2,-1]$.
\end{proof}

\begin{proposition}\label{prop:C-eta}
Assume (S), (F), (I), and that $M$ is $\eta$-Einstein with the constant $K$ of (F). Let
$p\in M$ and $\lambda_i\to\infty$. Then (C) holds. More precisely, after passing to a
subsequence and reparametrising, the rescalings $F^i$ converge smoothly on compact subsets
of $\Hn\times(-\infty,0)$ to a flow $F^\infty$ of proper Legendrian immersions for
$(\overline\eta,\overline\xi,\overline\Phi,\overline g)$, each component of which is
complete. Moreover $\partial_sF^\infty=H^\infty-2\theta^\infty\overline\xi$ and
$(\partial_s-\Delta)\theta^\infty=0$, where $\theta^\infty$ is the limit of the Legendrian
angles; and the push-forward of the measure of $F^\infty_s$ under $w$ has mass at most
$C_1R^n$ in $\{|w|<R\}$, for all $R>0$ and $s<0$, with the constant $C_1$ that
\cite[Lemmas~4.8 and~4.17]{ChangWuZhang2026ContactBlowup} provide for the flow $F$ (the
area constant of \cite{ChangWuZhang2026ContactBlowup}; not to be confused with the
constants $C_H$, $C_P$ above). The limit is non-empty if and only if $p\in F_T(\Sigma)$;
in particular it is non-empty at a singular point.
\end{proposition}
\begin{proof}
By (F), $H=\Phi\nabla\theta$, so the mean curvature form $d\theta$ is exact, and $F$
solves \cite[(1.1)]{ChangWuZhang2026ContactBlowup} with $\theta$ normalised by
\cite[(2.2)]{ChangWuZhang2026ContactBlowup}; (I) is the Type~I bound of
\cite[\S4.1]{ChangWuZhang2026ContactBlowup}. Hence $F$ is in the setting of
\cite[Theorem~1.1]{ChangWuZhang2026ContactBlowup}, with three differences. There $T$ is
the maximal existence time, the flow develops a Type~I singularity at time $T$, and $p$ is
a singular point. Here, by contrast, $T<\infty$ is any time up to which the flow is
defined with the Type~I bound (I), the flow need not become singular at $T$, and $p$ is
arbitrary. The lemmas of \cite{ChangWuZhang2026ContactBlowup} used below are stated in the
situation of \cite[Theorem~1.1]{ChangWuZhang2026ContactBlowup}, and we indicate for each
why its proof uses none of these three features. The rescalings $F^i$ are the contact
parabolic rescalings $F^{\lambda_i}$ of \cite[\S4.1]{ChangWuZhang2026ContactBlowup}, in
the chart $\Upsilon_p$, which is a chart as in
\cite[Lemma~4.1]{ChangWuZhang2026ContactBlowup} (Remark~\ref{rem:charts-M}).

\cite[Lemma~4.8]{ChangWuZhang2026ContactBlowup} for regular $M$, and
\cite[Lemma~4.17]{ChangWuZhang2026ContactBlowup} for arbitrary $M$, give for every $S>1$,
$R>0$ and $m\ge0$, for $i$ large and $s\in[-S,-S^{-1}]$: $|\nabla^mB^i|\le C_m$ at the
points of $F^i_s$ with $|w|,|z|\le R$; $|\theta^i|\le C_\theta$; and the area of the part
of $F^i_s$ over $\{|w|\le R\}$ that lies in a rescaled inner box in this chart, as in
\cite[Lemma~4.17]{ChangWuZhang2026ContactBlowup}, is at most $C_1R^n$; these inner boxes
exhaust $\Hn$. The
derivative bounds come from the Type~I bound on compact time intervals of the whole closed
flow, in the $\cD$-homothetic structures, through local estimates for the projected flow
\cite[Lemma~4.8]{ChangWuZhang2026ContactBlowup} or through the evolution inequalities on
the closed manifold $\Sigma$ and the maximum principle
\cite[Lemma~4.17]{ChangWuZhang2026ContactBlowup}. The angle bound is the maximum principle
bound $C_\theta$. The area bound comes from
\cite[Lemma~4.6(ii)]{ChangWuZhang2026ContactBlowup}, respectively from the last assertion
of \cite[Proposition~4.16]{ChangWuZhang2026ContactBlowup}, at an arbitrary centre. These
rest on the almost monotonicity of \cite[Proposition~4.5]{ChangWuZhang2026ContactBlowup},
a pointwise computation along the flow, on the small-scale area bound
\cite[Lemma~4.15]{ChangWuZhang2026ContactBlowup}, which uses only the closedness of
$\Sigma$ and $|\theta|\le C_\theta$, and on the compactness of $\Sigma_0$. None of these
arguments uses that $T$ is maximal, that the flow becomes singular at $T$, or that $p$ is
a singular point. The time $T$ enters only as the centre time of the kernels and as the
end of the interval on which (I) holds. Since
\cite[Lemma~4.17]{ChangWuZhang2026ContactBlowup} applies to every compact $\eta$-Einstein
$M$, regular ones included, we take for $C_1$ in the statement the area constant of that
lemma, that is, of the last assertion of
\cite[Proposition~4.16]{ChangWuZhang2026ContactBlowup}, rescaled; for regular $M$ the
constant of \cite[Lemma~4.6(ii)]{ChangWuZhang2026ContactBlowup} would serve as well. For
arbitrary $M$ this is the estimate that gives $C_P$ in the proof of Lemma~\ref{lem:D},
applied in the flow box at $p$ rather than uniformly over a finite cover.

With these bounds, the proof of \cite[Lemma~4.9]{ChangWuZhang2026ContactBlowup}, through
the compactness theorem for proper immersions of Breuning \cite{bre} cited there, gives
the convergence with the stated properties; that proof uses only the bounds above, the
flow equation and the convergence of the structures. In
\cite[Lemma~4.9]{ChangWuZhang2026ContactBlowup} it is assumed that $F^{\lambda_i}(x_i,-1)$
lies in a fixed compact set for some $x_i$; this is used only to see that the limit is
non-empty, and here the limit may be empty. Since the convergence is smooth on compact
subsets of space-time, it is uniform in $s\in[-2,-1]$, and the angles converge with it; no
vertical recentring is applied. This is (C).

Finally, if $x\in F_T^{-1}(p)$, Lemma~\ref{lem:confinement} keeps the rescaled images of
$x$ in a fixed compact set, so the limit is non-empty; a singular point lies in
$F_T(\Sigma)$ by the non-emptiness step of the proof of
\cite[Theorem~1.1]{ChangWuZhang2026ContactBlowup}. If $p\notin F_T(\Sigma)$, then, since
$F_t\to F_T$ uniformly (Lemma~\ref{lem:confinement}) and $F_T(\Sigma)$ is compact,
$F_t(\Sigma)$ stays at a positive distance from $p$ for $t$ close to $T$, and the rescaled
images leave every compact set.
\end{proof}

\begin{theorem}\label{thm:planes-eta}
Assume (S), (F), (I), and that $M$ is $\eta$-Einstein with the constant $K$ of (F); these
hypotheses contain the setting of \cite[Theorem~1.1]{ChangWuZhang2026ContactBlowup}, in
which $T$ is the maximal time, the flow becomes singular at $T$ and $p$ is a singular
point, but none of these three properties is assumed here (see the proof of
Proposition~\ref{prop:C-eta}). Let $p\in M$ and $\lambda_i\to\infty$. Then a subsequence
of the contact parabolic rescalings at $(p,T)$ converges, smoothly on compact subsets of
$\Hn\times(-\infty,0)$, to a flow $F^\infty$ with the following properties.
\begin{enumerate}[(a)]
\item For every $s<0$, each component of $F^\infty_s$ is a horizontal Legendrian plane
$P_j\times\{4|s|\bar\theta_j\}$ over a Lagrangian subspace $P_j\subset\C^n$, where
$\bar\theta_j$ is the limit of the Legendrian angle on that component. There are at most
$2^nC_1/\omega_n$ of them, counted with multiplicity, where $C_1$ is the area constant of
Proposition~\ref{prop:C-eta} and $\omega_n$ is the volume of the Euclidean unit ball of
$\R^n$.
\item $F^\infty$ satisfies $H-2\theta^\infty\overline\xi=\frac1{2s}X^\perp$ for $s<0$, and
$F^\infty_s=\delta_{\sqrt{-s}}F^\infty_{-1}$ up to reparametrisation, where each component
$P_j\times\{4\bar\theta_j\}$ of $F^\infty_{-1}$ is a Legendrian self-shrinker in the sense
of \cite[(1.2)]{ChangWuZhang2026ContactBlowup}. Thus, if $p\in F_T(\Sigma)$, in particular
if $p$ is a singular point, the conclusion of
\cite[Theorem~1.1]{ChangWuZhang2026ContactBlowup} holds for a subsequence of \emph{every}
sequence $\lambda_i\to\infty$, with $c=0$ for every component, and the limit is a finite
union of Legendrian planes, as in \cite[Corollary~1.2]{ChangWuZhang2026ContactBlowup}. For
$p\notin F_T(\Sigma)$ the limit is empty by (c); the limits of
\cite[Theorem~1.1]{ChangWuZhang2026ContactBlowup} are non-empty, and the comparison is
made only for $p\in F_T(\Sigma)$.
\item $F^\infty$ is non-empty if and only if $p\in F_T(\Sigma)$, in particular if $p$ is a
singular point.
\end{enumerate}
\end{theorem}
\begin{proof}
\emph{(a).} By Proposition~\ref{prop:C-eta}, (C) holds along a subsequence, with
convergence on compact subsets of $\Hn\times(-\infty,0)$.
Proposition~\ref{prop:centring}(a)--(c) and Corollary~\ref{cor:planes} hold, with the same
proofs, on every compact interval $[s_1,s_2]\subset(-\infty,0)$ in place of $[-2,-1]$. In
the derivation of \eqref{eq:centring-diss} take
$I_i=[T+s_1\lambda_i^{-2},T+s_2\lambda_i^{-2}]$, that is, Proposition~\ref{prop:energy}
with $S_0=-s_1$ and $\delta=-s_2$. Moreover, Proposition~\ref{prop:transport}, used in the
proof of Corollary~\ref{cor:planes}, holds on every compact subinterval of $[-S_0,0)$.
This gives the planes. On a horizontal plane the metric $\overline g$ restricts to
$g_0=\frac14|dw|^2$, so the part of $P_j\times\{4|s|\bar\theta_j\}$ over $\{|w|<R\}$ has
area $2^{-n}\omega_nR^n$. The area bound of Proposition~\ref{prop:C-eta} therefore allows
at most $2^nC_1/\omega_n$ planes, counted with multiplicity.

\emph{(b).} On the $j$-th component, $H=0$, $\theta^\infty=\bar\theta_j$ and
$z=4|s|\bar\theta_j$. The horizontal part $X^{\mathcal H}$ of $X$ is the horizontal lift
of $\mathfrak r$, which is tangent to $P_j$, and $\overline\xi$ is normal; so
$X^\perp=z\,\overline\xi$, and
\[
H-2\theta^\infty\overline\xi=-2\bar\theta_j\,\overline\xi=-\frac{z}{2|s|}\,\overline\xi=\frac1{2s}X^\perp
.
\]
Since $P_j$ is linear,
$\delta_{\sqrt{-s}}(P_j\times\{4\bar\theta_j\})=P_j\times\{4|s|\bar\theta_j\}$. On
$P_j\times\{4\bar\theta_j\}$ one has $z=4\bar\theta_j=4\theta^\infty$, and the case $s=-1$
of the displayed identity is \cite[(1.2)]{ChangWuZhang2026ContactBlowup}. In
\cite[Theorem~1.1]{ChangWuZhang2026ContactBlowup} the components along an arbitrary
sequence are vertical translates $\tau_{-4c_j}(\delta_{\sqrt{-s}}S_j)$; here $c_j=0$ for
every $j$, and each $S_j$ is a plane.

\emph{(c).} This is the last assertion of Proposition~\ref{prop:C-eta}.
\end{proof}

\begin{remark}[Comparison with {\cite[Sections~4.3
and~4.6]{ChangWuZhang2026ContactBlowup}}]\label{rem:M-compare}\leavevmode
\begin{enumerate}[(1)]
\item In the setting of \cite[Theorem~1.1]{ChangWuZhang2026ContactBlowup},
Theorem~\ref{thm:planes-eta} is stronger in three respects.
\cite[Theorem~1.1]{ChangWuZhang2026ContactBlowup} gives $c=0$ along a suitable sequence;
along an arbitrary sequence the components are vertical translates
$\tau_{-4c_j}(\delta_{\sqrt{-s}}S_j)$ of self-similar flows.
Theorem~\ref{thm:planes-eta}(b) gives $c_j=0$ for every component, along a subsequence of
every sequence. \cite[Corollary~1.2]{ChangWuZhang2026ContactBlowup} identifies the
self-shrinkers $S_j$ as planes through \cite[Theorem~A]{cqz}, after showing that $S_j$ is
orientable and normalising the angle $\theta^{\rm cqz}$ component by component
\cite[Remark~4.10, Corollary~4.11]{ChangWuZhang2026ContactBlowup};
Theorem~\ref{thm:planes-eta}(a) uses neither. The bound $2^nC_1/\omega_n$ on the number of
components replaces the bound $N_0$ of the proof of
\cite[Theorem~1.1]{ChangWuZhang2026ContactBlowup}, and the centre $p$ need not be a
singular point.
\item The proof of \cite[Theorem~1.1]{ChangWuZhang2026ContactBlowup} has the steps:
non-emptiness; compactness; self-similarity of the horizontal part; a bound on the number
of components; and the vertical centre. Here non-emptiness is Lemma~\ref{lem:confinement},
uniform in $p$ and in $x\in F_T^{-1}(p)$; compactness is used as it stands
(Proposition~\ref{prop:C-eta}); the horizontal self-similarity, $\cS^\infty=0$, comes from
the height-blind case of Proposition~\ref{prop:mono}, the analogue of
\cite[Propositions~4.5 and~4.16]{ChangWuZhang2026ContactBlowup}. The step on the vertical
centre, which in \cite{ChangWuZhang2026ContactBlowup} writes the horizontally self-similar
flows as vertical translates of self-shrinkers
(\cite[Lemma~2.5]{ChangWuZhang2026ContactBlowup}) and passes to the limits $\sigma\to0$ of
the dilations of the closed set of all limits, is replaced by
Proposition~\ref{prop:centring}(a). On $\tau_{-4c}(\delta_{\sqrt{-s}}S)$ with
$z=4\theta^\infty$ on $S$ one has $v^\infty=-4c/|s|$, and Fatou's lemma applied to the
finite dissipation $\int\!\!\int\psi D^\chi\rho$ gives $v^\infty\chi'(v^\infty)=0$, hence
$c=0$. The identification of the components as planes is replaced by $H^\infty=0$, which
comes from the $\theta^2$-monotonicity of Proposition~\ref{prop:theta2}, the Legendrian
form of \cite[Lemma~5.4]{nev07}, and the elementary argument of
Corollary~\ref{cor:planes}.
\item The transverse kernel cannot give $c=0$: it is constant along the Reeb orbits
\cite[Remark~4.4]{ChangWuZhang2026ContactBlowup}, so its monotonicity does not distinguish
$\tau_{-4c}(\delta_{\sqrt{-s}}S)$ from $\delta_{\sqrt{-s}}S$. The defect
$|v\chi'(v)|/\tau$ of the weighted kernel vanishes only on the exactly centred flows.
\item \emph{Where Type~I enters.} Type~I is used only in (C) (in the $\eta$-Einstein case
through \cite[Lemmas~4.8 and~4.17]{ChangWuZhang2026ContactBlowup}), in
Lemma~\ref{lem:confinement}, hence in Proposition~\ref{prop:centring}(d), and in
Corollary~\ref{cor:planes} and Theorem~\ref{thm:planes-eta}, where it bounds the second
fundamental form of the limit. Propositions~\ref{prop:theta2} and~\ref{prop:mono} hold
without it, and this makes Section~\ref{sec:weak} possible. The derivative bounds of
\cite[Lemmas~4.8 and~4.17]{ChangWuZhang2026ContactBlowup} are proved for $\eta$-Einstein
manifolds. In the general Sasakian setting, (C) remains a hypothesis of
Proposition~\ref{prop:centring} and Corollary~\ref{cor:planes}.
\item \emph{What is not covered.}\label{item:recentring} By Lemma~\ref{lem:confinement},
whose constants do not depend on $x\in F_T^{-1}(p)$, the rescaled images of all points of
$F_T^{-1}(p)$ stay in a fixed compact set. Hence a component visible only after unbounded
vertical recentring consists of limits of rescaled images of points $x_i$ with
$F_T(x_i)\ne p$. Such components are not described here. \cite[Corollary~1.3 and
Lemma~4.12]{ChangWuZhang2026ContactBlowup} concern the transverse density of a regular
manifold and White's local regularity theorem for the projected flow; the corresponding
statement for the contact density is not proved here.
\end{enumerate}
\end{remark}

\subsection{Two planes with different Legendrian angles}\label{subsec:two-planes}
Fix $p\in M$ and $S_0>3$. A \emph{two-plane configuration} is a weak contact tangent flow
of the form
\[
\cP(P_1,P_2)\colon\quad \mu_s=(\cH^n\llcorner
P_1)\times\delta_{4|s|\bar\theta_1}+(\cH^n\llcorner
P_2)\times\delta_{4|s|\bar\theta_2}\qquad(s\in[-S_0,0)),
\]
with multiplicity-one Lagrangian planes $P_1,P_2\subset\C^n$, $\cH^n$ taken for $g_0$, and
Legendrian angles $\bar\theta_1\ne\bar\theta_2$ on the two sheets. We consider the
hypothesis
\begin{enumerate}
\item[(P)] some weak contact tangent flow at $(p,T)$ is a two-plane configuration
$\cP(P_1,P_2)$.
\end{enumerate}
\noindent Put $\mathsf b^2:=8C_\theta+2$ and, in the chart $\Upsilon_p$,
\[
\mathsf d(q,t):=\max\big\{|w(q)|,\ (T-t)^{1/2},\ \mathsf b^{-1}|z_p(q)|^{1/2}\big\}.
\]

\begin{proposition}\label{prop:two-planes}
Assume (S), (F) and (P).
\begin{enumerate}[(a)]
\item Every weak contact tangent flow at $(p,T)$ is a two-plane configuration
$\cP(P_1',P_2')$, at every $s\in[-S_0,0)$, with multiplicity-one Lagrangian planes and the
\emph{same} Legendrian angles $\bar\theta_1,\bar\theta_2$.
\item For every $\eps\in(0,1]$ there is $r_1>0$ such that every point $q=F_t(x)$ of the
flow with $t<T$ and $\mathsf d(q,t)\le r_1$ satisfies
\[
\min_{j=1,2}\big|z_p(q)-4(T-t)\bar\theta_j\big|\le\eps\,\mathsf d(q,t)^2 .
\]
In particular $|z_p(q)|\le(4C_\theta+1)\mathsf d^2<\mathsf b^2\mathsf d^2$, so the
vertical term never realises the maximum in $\mathsf d$: near $(p,T)$ the flow lies in the
union of the two bands $\{|z_p-4(T-t)\bar\theta_j|\le\eps\,\mathsf d^2\}$, $j=1,2$, around
the heights $4(T-t)\bar\theta_1$ and $4(T-t)\bar\theta_2$. Each band has half-width
$\eps\,\mathsf d^2$ and their centres are $4(T-t)|\bar\theta_1-\bar\theta_2|$ apart, so
the bands are disjoint where $\eps\,\mathsf d^2<2(T-t)|\bar\theta_1-\bar\theta_2|$, in
particular on $\{\mathsf d^2\le C(T-t)\}$ for any $C\ge1$ with $\eps
C<2|\bar\theta_1-\bar\theta_2|$, and there no point of the flow lies between them. Where
$\mathsf d^2$ is large compared with $T-t$, the bands may overlap.
\end{enumerate}
\end{proposition}
\begin{proof}
\emph{(a).} By (P) and Corollary~\ref{cor:angle-measure}(b),
$\cL_{p,T}(\varpi)=\varpi(\bar\theta_1)+\varpi(\bar\theta_2)$. For any other weak contact
tangent flow, the same corollary gives
$\sum_jm_j'\,\delta_{\bar\theta_j'}=\delta_{\bar\theta_1}+\delta_{\bar\theta_2}$, with
$m_j'$ the vertex densities of the cones $\mathcal C_j'$. Thus there are exactly two
cones, with Legendrian angles $\bar\theta_1,\bar\theta_2$ and vertex density $1$. A
stationary integral cone of vertex density $1$ has all density ratios at the vertex equal
to $1$, so Allard's regularity theorem \cite[24.2]{si83} applies at the vertex at every
scale. The cone is therefore a $C^{1}$ graph near the vertex, hence a multiplicity-one
plane, and it is Lagrangian by Theorem~\ref{thm:weak}(iv). By Theorem~\ref{thm:weak}(v)
the description holds at every $s$.

\emph{(b).} Suppose the assertion fails for some $\eps$. Then there are flow points
$q_k=F_{t_k}(x_k)$ violating the inequality with $r_k:=\mathsf d(q_k,t_k)\to0$. Rescale at
$(p,T)$ by $r_k^{-1}$. The rescaled points $(\hat q_k,\hat s_k)$, $\hat s_k\in[-1,0]$,
satisfy
\[
\max\{|\hat w_k|,|\hat s_k|^{1/2},\mathsf b^{-1}|\hat z_k|^{1/2}\}=1,\qquad \min_j|\hat
z_k-4|\hat s_k|\bar\theta_j|>\eps .
\]
Pass to a subsequence along which $(\hat q_k,\hat s_k)\to(\hat q,\hat s)$ and the rescaled
flows converge at every time (Proposition~\ref{prop:compactness}(c)). By part (a) the
limit is some $\mathcal P(P_1',P_2')$ at every time, and $\min_j|\hat z-4|\hat
s|\bar\theta_j|\ge\eps$.

A point $\ell_{\hat q}(w,z)$ of the closed model box $\ell_{\hat q}\{f\le
R_*^2c_B,\ |z|\le Z_*c_B\}$, with $\hat q=(\hat w,\hat z)$ and $|\hat w|\le1$, has
$|w|\le4R_*\sqrt{c_B}$ and $|z|\le Z_*c_B$, and by the group law its height is $\hat
z+z+\frac12(\hat b\cdot x-\hat a\cdot y)$ with $\hat w=\hat a+i\hat b$, so that
\begin{equation}\label{eq:box-height}
\big|\text{height}-\hat z\big|\le Z_*c_B+\tfrac12|\hat w||w|\le Z_*c_B+2R_*\sqrt{c_B}\le
C_0(R_*\sqrt{c_B}+Z_*c_B),\qquad C_0:=2 .
\end{equation}
Fix $c_B\in(0,\frac18)$ with $C_0(R_*\sqrt{c_B}+Z_*c_B)+4c_BC_\theta<\eps/2$, and apply
Lemma~\ref{lem:clearing} at $(q_k,t_k)$ with $\sigma_k:=r_k^2(c_B+\hat s_k-\hat s)$. This
is positive and at most $2c_Br_k^2$ for $k$ large, and the evaluation time $t_k-\sigma_k$
is the \emph{fixed} rescaled time $\hat s-c_B<0$. Since $q_k$ lies on the flow, the box
$\mathcal Q_{\sigma_k}(q_k)$ carries mass $>c_{\rm cl}\sigma_k^{n/2}$ at that time. By
Lemma~\ref{lem:moving}(c) the rescaled boxes converge to the closed model box
$\overline{\mathcal Q}=\ell_{\hat q}\{f\le R_*^2c_B,\ |z|\le Z_*c_B\}$; the centred box
$\{f\le R_*^2c_B,\ |z|\le Z_*c_B\}$ is $U(n)$-invariant, so the unitary ambiguity does not
matter. Weak convergence at time $\hat s-c_B$ gives $\mu_{\hat s-c_B}(\overline{\mathcal
Q})\ge c_{\rm cl}c_B^{n/2}>0$. On the other hand, $\mu_{\hat s-c_B}$ is carried by the
planes $P'_j\times\{4(|\hat s|+c_B)\bar\theta_j\}$. Every point of $\overline{\mathcal Q}$
has height within $\eps_h:=C_0(R_*\sqrt{c_B}+Z_*c_B)$ of $\hat z$ by \eqref{eq:box-height}
(recall $|\hat w|\le1$), hence at distance $\ge\eps-4c_BC_\theta-\eps_h>\eps/2$ from
$4(|\hat s|+c_B)\bar\theta_j$ for both $j$. Hence $\mu_{\hat s-c_B}(\overline{\mathcal
Q})=0$, a contradiction.

That the vertical term never realises the maximum in $\mathsf d$ was derived in the
statement from $T-t\le\mathsf d^2$ and $\eps\le1$. Equivalently, rescale a flow point $q$
with $\mathsf d(q,t)\le r_1$ at $(p,T)$ by $\mathsf d(q,t)^{-1}$, so that $|\hat s|\le1$;
if $\mathsf b^{-1}|\hat z|^{1/2}=1$, then $|\hat z-4|\hat s|\bar\theta_j|\ge\mathsf
b^2-4C_\theta=4C_\theta+2>1\ge\eps$, contrary to the inequality. The assertions on the
bands follow from the inequality by the interval computation given in the statement.
\end{proof}

By (b), near $(p,T)$ the flow lies in two disjoint height bands on the regions $\{\mathsf
d^2\le C(T-t)\}$, which shrink to the point $p$ as $t\to T$. We do not show that on a
fixed space-time neighbourhood of $(p,T)$ the flow consists of two smooth sheets, nor that
a Legendrian double point forms at $(p,T)$ without curvature blow-up.

\subsection{Legendrian curves}\label{subsec:curves}
For $n=1$ the closed manifold $\Sigma$ is a finite union of circles, so (O) holds and
Theorem~\ref{thm:weak}(vi) applies. We identify the special Lagrangian cones of $\C=\C^1$.

\begin{lemma}\label{lem:lines}
Let $\mathcal T$ be a nonzero integer multiplicity rectifiable $1$-current in $\C$ with
locally finite mass, $\partial \mathcal T=0$, calibrated by $\Omega^{\rm
Re}=\mathrm{Re}(e^{-i\bar\vartheta}d\zeta)$ for some $\bar\vartheta\in\R/2\pi\Z$, where
$\zeta=\frac12(x-iy)$ and the comass is taken for $g_0$, and invariant under the dilations
$w\mapsto\lambda w$, $\lambda>0$. Then $\mathcal T=m\,[L]$ for an integer $m\ge1$, where
$[L]$ is the line through the origin oriented by the $g_0$-unit vector
$e_{\bar\vartheta}=2(\cos\bar\vartheta\,\partial_x-\sin\bar\vartheta\,\partial_y)$.
\end{lemma}
\begin{proof}
Write $\mathcal T=\vec{\mathcal T}\wedge\|\mathcal T\|$ with the orienting unit vector
$\vec{\mathcal T}$. For a $g_0$-unit vector
$\tau=\cos\alpha\,(2\partial_x)+\sin\alpha\,(2\partial_y)$ one has
$d\zeta(\tau)=e^{-i\alpha}$ and $\Omega^{\rm Re}(\tau)=\cos(\bar\vartheta+\alpha)$, which
equals $1$ exactly for $\alpha=-\bar\vartheta$. Since $\mathcal T$ is calibrated,
$\Omega^{\rm Re}(\vec{\mathcal T})=1$ $\|\mathcal T\|$-a.e., so $\vec{\mathcal
T}=e_{\bar\vartheta}=:e$ a.e.

Then $\partial \mathcal T=0$ reads
$\int\partial_ef\,d\|\mathcal T\|=0$ for all $f\in C^\infty_c(\C)$: the distributional
derivative of the Radon measure $\|\mathcal T\|$ in the direction $e$ vanishes.
Mollifying, $\|\mathcal T\|*j_\eps$ is a smooth function with vanishing derivative in the
direction $e$, hence invariant under the translations $w\mapsto w+te$; letting $\eps\to0$,
$\|\mathcal T\|$ is invariant under these translations. Consequently $\|\mathcal
T\|=\cL^1_e\otimes\nu$, where $\cL^1_e$ is the $g_0$-length on the line $\R e$ and $\nu$
is a Radon measure on the orthogonal line $\R e^\perp$: indeed $\nu(A):=\|\mathcal
T\|([0,1)e+A)$ defines a Radon measure, and by invariance $\|\mathcal
T\|([a,b)e+A)=(b-a)\nu(A)$ for rational $a<b$, hence for all $a<b$.

The push-forward of
$\mathcal T$ under $w\mapsto\lambda w$ has mass measure
$\lambda\,(\lambda\cdot)_\#\|\mathcal T\|$, so the dilation invariance of $\mathcal T$
gives $\|\mathcal T\|(I\times A)=\lambda\|\mathcal T\|(\lambda^{-1}I\times\lambda^{-1}A)$
for intervals $I\subset\R e$ and Borel $A\subset\R e^\perp$, that is,
$\nu(A)=\nu(\lambda^{-1}A)$ for all $\lambda>0$. For $A=(0,r]e^\perp$ this gives
$\nu((0,r]e^\perp)=\nu((0,r/\lambda]e^\perp)\to\nu(\emptyset)=0$ as $\lambda\to\infty$, by
continuity from above, and likewise on the other side. Hence $\nu=m\delta_0$ and
$\|\mathcal T\|=m\cL^1_e=m\cH^1_{g_0}\llcorner L$, where $\cH^1_{g_0}$ is the
one-dimensional Hausdorff measure of $g_0$, which on a line is half the Euclidean length.
Thus $\mathcal T=m[L]$, and $m$ is a positive integer because $\mathcal T$ is a nonzero
integer multiplicity current.
\end{proof}

\begin{theorem}\label{thm:curves}
Let $n=1$ and assume (S), (F). Let $p\in M$ and $\lambda_i\to\infty$, and let
$(\mu_s,\fm_s)_{s<0}$ be a weak contact tangent flow at $(p,T)$ for all negative times,
along the diagonal subsequence constructed after the proof of Theorem~\ref{thm:weak},
described as in Theorem~\ref{thm:weak}(v) by the Legendrian angles
$\bar\theta_1<\dots<\bar\theta_N$ and the cones $\cC_1,\dots,\cC_N$; for $k\in\N$ let
$\cG_k\subset[-k,0)$ be the set of good times of Proposition~\ref{prop:compactness}(c) for
the restriction to $[-k,0)$.
\begin{enumerate}[(a)]
\item Each $\cC_j$ is a finite union of lines through the origin of $\C$, with integer
multiplicities: $\|\cC_j\|=\sum_km_{jk}\,\cH^1_{g_0}\llcorner L_{jk}$ for finitely many
lines $L_{jk}$ through $0$ and integers $m_{jk}\ge1$, where $\cH^1_{g_0}$ is the Hausdorff
measure of $g_0=\frac14|dw|^2$. Hence, for every $s<0$, $\mu_s$ is the weight of the union
of the horizontal Legendrian lines $L_{jk}\times\{4|s|\bar\theta_j\}$ through the Reeb
axis, counted with the multiplicities $m_{jk}$, and $m_j=\sum_km_{jk}$ in
Corollary~\ref{cor:angle-measure}.
\item Let $R>0$, let $s_0\in\cG_k$ be a good time for some $k$, write $F:=F^i_{s_0}$ and
let $B_r$ denote the Euclidean ball of radius $r$ about $0$ in the coordinates $(w,z)$ of
$\R^3$. Let $\Sigma_i\subset\Sigma$ be open sets with $F(\partial\Sigma_i)\cap
B_{3R}=\emptyset$, such that $\Sigma_i\cap F^{-1}(B_{2R})$ has exactly one connected
component meeting $F^{-1}(B_R)$. Then, after passing to a subsequence, there is $\bar
z\in[-R,R]$ with:
\begin{enumerate}[(i)]
\item $\sup_{\Sigma_i\cap F^{-1}(B_R)}|z\circ F-\bar z|\to0$: a single height;
\item with $\bar\theta:=\bar z/(4|s_0|)$, $\theta^i_{s_0}\to\bar\theta$ in $\mu$-measure
on $\Sigma_i\cap F^{-1}(B_R)$, and for all $\gamma\in C_c(B_R)$ and $\varpi\in C(\R)$,
\[
\int_{\Sigma_i}\varpi(\theta^i_{s_0})\,\gamma\circ
F\,d\mu-\varpi(\bar\theta)\int_{\Sigma_i}\gamma\circ F\,d\mu\longrightarrow0,
\]
a single value of the Legendrian angle; if the lengths of $\Sigma_i\cap F^{-1}(B_R)$ do
not tend to $0$, then $|\bar\theta|\le C_\theta$;
\item if moreover, for some fixed $\delta_0\in(0,R)$, the component meets
$F^{-1}(B_{R-\delta_0})$ for all $i$, then the measures $(F|_{\Sigma_i})_\#\mu\llcorner
B_R$ converge to a \emph{nonzero} measure, the weight in $B_R$ of a stationary integral
Legendrian varifold supported in $\{z=\bar z\}$.
\end{enumerate}
\end{enumerate}
In particular, a locally connected component of the rescaled flow cannot split into parts
of different heights or of different Legendrian angles in the limit.
\end{theorem}
\begin{proof}
\emph{(a).} Since $\Sigma$ is orientable, Theorem~\ref{thm:weak}(vi) writes
$\cC_j=\sum_k|\mathcal T_{jk}|$ as varifolds, where the $\mathcal T_{jk}$ are integer
multiplicity $1$-currents in $\C$ with locally finite mass, $\partial \mathcal T_{jk}=0$,
calibrated by $\mathrm{Re}(e^{-i\bar\vartheta_{jk}}\Omega_0)$ with $\Omega_0=d\zeta$, and
invariant under dilations. By Lemma~\ref{lem:lines}, $\mathcal T_{jk}=m_{jk}[L_{jk}]$ for
a line $L_{jk}$ through $0$ and an integer $m_{jk}\ge1$; distinct $k$ may give the same
unoriented line, with opposite orientations, and then the multiplicities add. The lifts
$L_{jk}\times\{\bar z\}$ are Legendrian because $\overline\beta$ vanishes on lines through
the origin. The description of $\mu_s$ for every $s<0$ is Theorem~\ref{thm:weak}(v),
applied on $[-k',0)$ for every $k'\in\N$ as explained after the proof of that theorem.
Finally the vertex density of a line through $0$ is $1$, so
$m_j=\Theta(\cC_j,0)=\sum_km_{jk}$.

\emph{(b).} Write $b_i:=z\circ F$; gradients, lengths and $\mu$ are taken for $g_i$, and
$\lambda_i$, $\hat a_i$, $\cS^i$, $v^i$ are as in Section~\ref{subsec:scaling}. Since
$s_0\in\cG_k$, along the diagonal subsequence from its $k$-th term on (which is a
subsequence of the one of Proposition~\ref{prop:compactness}(c) on $[-k,0)$, so that the
good-time conclusion persists),
\[
\eps_i:=\int_{F^{-1}(B_{3R})}\big(|H^i|^2+|\cS^i|^2+|v^i\chi'(v^i)|\big)\,d\mu\longrightarrow0
.
\]
Two estimates on $F^{-1}(B_{3R})$ are used. First, by Lemma~\ref{lem:model}(d), rescaled
with errors $O(\lambda_i^{-2})$ (Lemma~\ref{lem:moving}(b)),
\begin{equation}\label{eq:grad-z}
|\nabla b_i|^2\le16|(\overline\nabla
f)^\perp|^2+o(1)\le32|s_0|^2\big(|\cS^i|^2+|H^i|^2\big)+o(1),
\end{equation}
uniformly on $F^{-1}(B_{3R})$. Second, the metrics are comparable: for a horizontal vector
$U=(U_w,U_z)$ at a point of $B_{3R}$ one has $U_z=-2\beta^{(i)}(U_w)$, with
$\beta^{(i)}:=\beta^{(\lambda_i)}_p$ the rescaled $1$-form of the chart (proof of
Lemma~\ref{lem:moving}(b)), so $|U_w|^2\le g_{\rm euc}(U,U)\le(1+CR^2)|U_w|^2$ for the
Euclidean metric $g_{\rm euc}$ of $\R^3$. On the other hand,
$g_i(U,U)=h_i(U_w,U_w)\in[\frac18|U_w|^2,\frac12|U_w|^2]$ for $i$ large, with
$h_i=h^{(\lambda_i)}_p$ as in the proof of Theorem~\ref{thm:weak}(vi). Hence
\begin{equation}\label{eq:metric-comparison}
c(R)\,g_{\rm euc}(U,U)\le g_i(U,U)\le\tfrac12g_{\rm euc}(U,U)\qquad\text{for horizontal
}U\text{ at points of }B_{3R},\ i\text{ large}.
\end{equation}

For (i), let $\cD_i$ be the connected component of $\Sigma_i\cap F^{-1}(B_{2R})$ that
meets $F^{-1}(B_R)$. It contains $\Sigma_i\cap F^{-1}(B_R)$, and $\mu(\cD_i)\le2D_0R$ for
$i$ large, by Proposition~\ref{prop:compactness}(a). By \eqref{eq:grad-z},
$\int_{\cD_i}|\nabla b_i|^2\,d\mu\to0$. Since $\Sigma$ is a finite union of circles,
$\cD_i$ is a circle or an open arc, and any two of its points are joined by an arc in
$\cD_i$. Hence, uniformly in $x,y\in\cD_i$,
\[
|b_i(x)-b_i(y)|\le\int_{\cD_i}|\nabla b_i|\,d\mu\le(2D_0R)^{1/2}\Big(\int_{\cD_i}|\nabla
b_i|^2\,d\mu\Big)^{1/2}\longrightarrow0 .
\]
Choosing $x_i\in\Sigma_i\cap F^{-1}(B_R)$ and a subsequence with $b_i(x_i)\to\bar
z\in[-R,R]$ gives (i).

For (ii), note that on $\Sigma_i\cap F^{-1}(B_R)$ we have $v^i=b_i/|s_0|-\hat a_i\theta^i$
with $\hat a_i\to4$, so $v^i-4(\bar\theta-\theta^i)\to0$ uniformly by (i). Let $\mathsf
g(t)=|t\chi'(t)|$, which is positive for $t\ne0$. For $\delta>0$, on $\{x\in\Sigma_i\cap
F^{-1}(B_R):\ |\theta^i(x)-\bar\theta|\ge\delta\}$ we have $|v^i|\in[2\delta,C]$ for $i$
large, so the measure of this set is at most $\eps_i/\min_{[2\delta,C]}\mathsf g\to0$.
Thus $\theta^i\to\bar\theta$ in measure on $\Sigma_i\cap F^{-1}(B_R)$, which is the first
assertion of (ii). The second follows: if $|\bar\theta|\le C_\theta$, because the lengths
are bounded and $\varpi$ is uniformly continuous on $[-C_\theta,C_\theta]$; if
$|\bar\theta|>C_\theta$, because then the set above with $\delta=|\bar\theta|-C_\theta$ is
all of $\Sigma_i\cap F^{-1}(B_R)$, so the lengths tend to $0$, while $\varpi$ is bounded
on $[-C_\theta,C_\theta]$ and $\gamma$ is bounded. If the lengths do not tend to $0$,
convergence in measure of $\theta^i\in[-C_\theta,C_\theta]$ forces $|\bar\theta|\le
C_\theta$.

For (iii), let $y_i\in\cD_i$ with $F(y_i)\in B_{R-\delta_0}$. First, for $i$ large,
$\cD_i$ is not a whole component $\Sigma'$ of $\Sigma$. Otherwise $\mu(\Sigma')\le2D_0R$
and $\int_{\Sigma'}|H^i|^2d\mu\le\eps_i$. The projection $w\circ F|_{\Sigma'}$ is a closed
curve in the $w$-plane, isometrically immersed for the rescaled base metric $h_i$, and its
geodesic curvature has norm $|H^i|$ (\cite[Lemma~2.3]{ChangWuZhang2026ContactBlowup},
whose proof is pointwise and applies in the chart). Parametrise this curve by
$h_i$-arclength $\ell$, so that $d\mu=d\ell$, and let $e$ be its Euclidean unit tangent as
a function of $\ell$. Since $h_i\to g_0=\frac14|dw|^2$ in $C^1$ on $\{|w|\le2R\}$
and $g_0$ has constant coefficients, $|e'|\le2|H^i|+o(1)$ uniformly. As the curve is
closed, $\int\langle(w\circ F)',e(0)\rangle=0$; so $e$ does not stay at angle less than
$\pi/2$ from $e(0)$, and $\int_{\Sigma'}|e'|\ge\pi/2$. Since $\mu(\Sigma')\le2D_0R$, for
$i$ large this gives $\int_{\Sigma'}|H^i|\,d\mu\ge\pi/8$, and by Cauchy--Schwarz
$(\pi/8)^2\le2D_0R\,\eps_i\to0$, a contradiction.

Thus $\cD_i$ is an open arc, and we
use the margin $\delta_0$ to show that it has length bounded below inside $B_R$. Its ends
accumulate at points $x_*\notin\cD_i$ with $|F(x_*)|\ge2R$: either
$x_*\in\partial\Sigma_i$, which is mapped outside $B_{3R}$, or $x_*\in\Sigma_i$, and then
$F(x_*)\notin B_{2R}$, since otherwise $x_*\in\cD_i$. Hence $\cD_i$ contains an arc from
$y_i$ to a point of $F^{-1}(\partial B_{R-\delta_0/2})$ whose image lies in $\overline
B_{R-\delta_0/2}$ and has Euclidean length $\ge\delta_0/2$. By
\eqref{eq:metric-comparison} its $\mu$-measure is $\ge\frac12\sqrt{c(R)}\,\delta_0$. Since
the limit
measure of the closed ball $\overline B_{R-\delta_0/2}$ is at least the $\limsup$ of these
lengths, the limit is nonzero on $B_R$.

The immersions $F|_{\Sigma_i\cap F^{-1}(B_{2R})}$
have no boundary in $B_{2R}$: the boundary of this parameter set is mapped to $\partial
B_{2R}$ or outside $B_{3R}$. The argument in the part of the proof of
Proposition~\ref{prop:compactness}(c) headed ``The varifold limit at good times'', with
$\Hn$
replaced by $B_{2R}$ and compact sets $\cK_m$ exhausting $B_{2R}$, and the proof of
Theorem~\ref{thm:weak}(i), apply to them with $\int|H^i|^2\to0$. They give a stationary
integral Legendrian varifold in $B_{2R}$. By (i) its restriction to $B_R$ is supported in
$\{z=\bar z\}$. This proves (iii).
\end{proof}

\begin{remark}[The margin $\delta_0$]\label{rem:margin}
Without the margin in (iii), meeting the open ball $B_R$ does not give a lower mass bound
in $B_R$: a line at height $z_i\uparrow R$ meets $B_R$ for every $i$, but its measure in
$B_R$ tends to $0$. Such a configuration occurs at a fixed rescaled time $s_0$ for a flow
satisfying (F). Take a totally geodesic Legendrian great circle in the standard Sasakian
$S^3$, for which $K+2=4>0$ (the standard structure of $S^{2n+1}$ has $K=2n$;
cf.\ \cite[Example~2.7(a)]{s1}), moving along the Reeb field with
$\theta(t)=\theta_Te^{(K+2)(t-T)}$, $\theta_T>0$. In a chart at the final point whose
projection is a real line, $z(t)=4\theta_T(1-e^{-(K+2)\tau})/(K+2)$. Hence
$z_i=4|s_0|\theta_T-2(K+2)\theta_T|s_0|^2\lambda_i^{-2}+O(\lambda_i^{-4})$, which
increases to $R:=4|s_0|\theta_T$, because $K+2>0$. Here $H=0$, $(\overline\nabla
f)^\perp=0$ and $v=0$, so every slice is good, and the configuration occurs for every
$s_0$, with $R=4|s_0|\theta_T$. It does not contradict Theorem~\ref{thm:curves}(b),
because the margin hypothesis of (iii) fails: the line does not meet $B_{R-\delta_0}$ once
$z_i>R-\delta_0$. A version of (iii) without the margin would have to exclude, for each
$R$, a null set of times depending on $R$. The same example, with a totally geodesic
Legendrian $S^n\subset S^{2n+1}$ (for which $K+2=2n+2>0$), taken in a chart in which its
projection is the real $n$-plane, shows that in every dimension meeting $B_R$ does not
give a lower mass bound in $B_R$; a higher-dimensional analogue of (iii) would therefore
also need a margin.
\end{remark}

\appendix
\section{Finite-phase compactness for metrics close to a flat one}\label{app:A}

In this appendix we prove the finite-phase compactness used in the proof of
Theorem~\ref{thm:weak}(vi). It is a version, for metrics close to a flat one, of
\cite[Proposition~5.1]{nev07} and \cite[Proposition~3.11]{nev11}, whose proof is that of
\cite[Proposition~5.1]{nev07}. There the lower density bound for the limit comes from
Allard's compactness theorem and the monotonicity formula for the stationary limit
varifold; here it is obtained instead from near calibration of the limit currents, and
Allard's theorem is not used (Remark~\ref{rem:A2}).

Let $J_0$ and $\Omega_0$ be as in Section~\ref{subsec:sasakian}, and let $h_0$ be a
constant K\"ahler metric for $J_0$ with $|\Omega_0|_{h_0}=1$ in the sense of
Section~\ref{subsec:sasakian}, that is, $|\Omega_0(E_1,\dots,E_n)|=1$ for an $h_0$-unitary
frame; for instance $h_0=g_0$. Balls $B_r$ are $h_0$-balls centred at $0$. For every
$\bar\vartheta\in\R/2\pi\Z$ the form $\mathrm{Re}(e^{-i\bar\vartheta}\Omega_0)$ has
$h_0$-comass $1$: in $h_0$-unitary coordinates it is the special Lagrangian calibration of
\cite[Chapter~III, Theorem~1.10]{hl82}. It has constant coefficients, so it is closed, and
exact on $\C^n$.

\begin{proposition}[{cf.\ \cite[Proposition~5.1]{nev07},
\cite[Proposition~3.11]{nev11}}]\label{prop:A1}
Let $n\ge1$. Let $h_i$ be K\"ahler metrics for $J_0$ on $B_{4R}$ with
$\eps_i:=\|h_i-h_0\|_{C^2(B_{4R})}\to0$. Let $L^i$ be oriented $n$-manifolds, possibly
with boundary, and $F_i\colon L^i\to B_{4R}$ proper $\omega_{h_i}$-Lagrangian immersions.
Define the geometric phase $\vartheta_i\colon L^i\to\R/2\pi\Z$ by
$F_i^*\Omega_0=\big(|\Omega_0|_{h_i}\circ F_i\big)\,e^{i\vartheta_i}\,\mathrm{vol}_{h_i}$.
Assume
\begin{enumerate}[(a)]
\item $\mu_i(L^i)\le D_0R^n$, where $\mu_i$ is the $h_i$-area;
\item $\int_{L^i}|H_i|^2\,d\mu_i+\mathcal H^{n-1}_{h_i}(\partial L^i)\to0$.
\end{enumerate}
Then, after passing to a subsequence, there are finitely many distinct
$\bar\vartheta_1,\dots,\bar\vartheta_N\in\R/2\pi\Z$ and integral currents $\mathcal
T_1,\dots,\mathcal T_N$ in $B_R$ with $\partial \mathcal T_j=0$ in $B_R$, each $\mathcal
T_j$ calibrated by $\mathrm{Re}(e^{-i\bar\vartheta_j}\Omega_0)$ (special Lagrangian of
phase $\bar\vartheta_j$), such that for all $\gamma\in C_c(B_R)$ and $\varpi\in
C(\R/2\pi\Z)$
\[
\lim_{i\to\infty}\int_{L^i}\varpi(\vartheta_i)\,\gamma\circ
F_i\,d\mu_i=\sum_{j=1}^N\varpi(\bar\vartheta_j)\int\gamma\,d\|\mathcal T_j\| .
\]
If $\vartheta_i$ lifts to real functions $\tilde\vartheta_i$ with $|\tilde\vartheta_i|\le
D_0$, the same holds for $\varpi(\tilde\vartheta_i)$, with real phases
$\bar\vartheta_j\in[-D_0,D_0]$ and $\varpi\in C(\R)$. The number $N$ is bounded by a
constant depending only on $n$ and $D_0$.
\end{proposition}
\begin{proof}
All constants $C$ depend only on $n$ and $D_0$. We cut the phase circle into short arcs
and pass to limits of the currents carried by the corresponding pieces; these limits are
nearly calibrated. The main point is a lower density bound for these limits with a
constant independent of the arcs. Since the total mass is bounded, it allows only finitely
many phases, and the limits are then identified as calibrated currents.

The computation
in the proof of Lemma~\ref{lem:phase} gives $d\vartheta_i=\langle
H_i,J_0\cdot\rangle+F_i^*d^c\log|\Omega_0|_{h_i}$. Since $|d\log|\Omega_0|_{h_i}|\le
C\eps_i$, we get $|\nabla\vartheta_i|\le|H_i|+C\eps_i$, and hence, by (a) and (b),
\begin{equation}\label{eq:A-grad}
\int_{L^i}|\nabla\vartheta_i|\,d\mu_i\le(D_0R^n)^{1/2}\Big(\int|H_i|^2\Big)^{1/2}+C\eps_iR^n\to0
.
\end{equation}
The measures $\mathfrak m^i:=(F_i,\vartheta_i)_\#\mu_i$ on $B_{4R}\times(\R/2\pi\Z)$ have
bounded mass. Pass to a subsequence with $\mathfrak m^i\to\mathfrak m$ weakly; this Young
measure is fixed from now on. Let $\mu:=(\mathrm{pr}_1)_\#\mathfrak m$.

Fix $m\ge3$. We claim that, after passing to a further subsequence, which may depend on
$m$, the following hold for almost every $a\in[0,2\pi/m)$:
\begin{enumerate}[(i)]
\item $\mathfrak m(B_{4R}\times\{a+2\pi k/m\})=0$ for all $k$;
\item $A_i(a):=\sum_k\mathcal H^{n-1}\big(\vartheta_i^{-1}(a+2\pi k/m)\big)\to0$ along the
subsequence;
\item all $a+2\pi k/m$ are regular values of every $\vartheta_i$.
\end{enumerate}
Indeed, the second marginal of $\mathfrak m$ is a finite measure on the circle, so it has
at most countably many atoms, and (i) fails for at most countably many $a$. By the coarea
formula for the circle-valued map $\vartheta_i$ and \eqref{eq:A-grad},
$\int_0^{2\pi/m}A_i(a)\,da=\int_{L^i}|\nabla\vartheta_i|\,d\mu_i\to0$; choose the
subsequence so that these integrals are summable, then $\sum_iA_i(a)<\infty$ for almost
every $a$ by monotone convergence, which gives (ii). By Sard's theorem the critical values
of each $\vartheta_i$ form a null set, and (iii) fails only for $a$ in the countable union
over $i$ and $k$ of the translates of these null sets. The further subsequence does not
change $\mathfrak m$. Let $a$ be any value for which (i)--(iii) hold, let $I_k=(a+2\pi
k/m,\,a+2\pi(k+1)/m)$, $L^i_k=\vartheta_i^{-1}(I_k)$, and let $\bar\vartheta_k^{(m)}$ be
the midpoint of $I_k$. Then $\mathcal H^{n-1}(\partial L^i_k)\to0$, by (b) and (ii).

The currents $\mathcal T^i_k:=(F_i)_\#[L^i_k]$ are integral. They satisfy $\mathbf
M_{h_0}(\mathcal T^i_k)\le(1+C\eps_i)\mu_i(L^i_k)$ and $\mathbf M(\partial \mathcal
T^i_k\llcorner B_{4R})\to0$. By the Federer--Fleming compactness theorem
\cite[27.3]{si83}, a subsequence converges weakly in $B_{4R}$ to integral currents
$\mathcal T_k$ with $\partial \mathcal T_k=0$ in $B_{4R}$. Moreover $\mu_i\llcorner
L^i_k$, pushed forward, converges to $\mu_k:=(\mathrm{pr}_1)_\#(\mathfrak m\llcorner
B_{4R}\times I_k)$, by (i). By lower semicontinuity of mass, $\|\mathcal T_k\|\le\mu_k$.
For $\gamma\in C_c(B_{4R})$ with $\gamma\ge0$,
\[
\mathcal
T^i_k\big(\gamma\,\mathrm{Re}(e^{-i\bar\vartheta_k^{(m)}}\Omega_0)\big)=\int_{L^i_k}\gamma\circ
F_i\,|\Omega_0|_{h_i}\cos\big(\vartheta_i-\bar\vartheta_k^{(m)}\big)\,d\mu_i\ge(1-C\eps_i)\cos(\pi/m)\int_{L^i_k}\gamma\circ
F_i\,d\mu_i ,
\]
and letting $i\to\infty$ and using that the comass is $1$,
\begin{equation}\label{eq:A1}
\cos(\pi/m)\,\mu_k(\gamma)\le \mathcal
T_k\big(\gamma\,\mathrm{Re}(e^{-i\bar\vartheta_k^{(m)}}\Omega_0)\big)\le\|\mathcal
T_k\|(\gamma)\le\mu_k(\gamma).
\end{equation}

\begin{claim}
For every $k$, every $x\in\spt\|\mathcal T_k\|\cap B_{2R}$ and every $r\le R$,
\begin{equation}\label{eq:A-density}
\|\mathcal T_k\|(B_r(x))\ge c_0r^n,\qquad c_0=c_0(n)>0 .
\end{equation}
\end{claim}
The proof uses only that $\mathcal T_k$ is an integral current with $\partial\mathcal
T_k=0$ in $B_{4R}$ and the near calibration \eqref{eq:A1} with $\cos(\pi/m)\ge\frac12$;
the mechanism differs for $n\ge2$ and $n=1$. Let $\mathcal V(r)=\|\mathcal T_k\|(B_r(x))$;
then $\mathcal V>0$ on $(0,R]$, and
$\mathcal V$ is non-decreasing. For a.e.\ $r$ the slice $S_r=\partial(\mathcal
T_k\llcorner B_r(x))$ is an integral $(n-1)$-current supported in $\partial B_r(x)$ with
$\mathbf M(S_r)\le\mathcal V'(r)$, because $\partial \mathcal T_k=0$ in $B_{4R}$. Write
$\Omega^{\rm Re}_k:=\mathrm{Re}(e^{-i\bar\vartheta_k^{(m)}}\Omega_0)$. By \eqref{eq:A1}
with $\gamma\uparrow1_{B_r(x)}$, for a.e.\ $r$,
\[
(\mathcal T_k\llcorner B_r(x))(\Omega^{\rm Re}_k)\ge\cos(\pi/m)\,\mathcal
V(r)\ge\tfrac12\mathcal V(r).
\]

For $n\ge2$, the isoperimetric theorem \cite[30.1]{si83} gives an integral current $Q$
with $\partial Q=S_r$ and $\mathbf M(Q)\le c_n\mathbf M(S_r)^{n/(n-1)}$. Replacing $Q$ by
its push-forward under the nearest-point retraction onto $\overline B_r(x)$, which is
$1$-Lipschitz and the identity on $\overline B_r(x)\supset\spt S_r$, we may assume $\spt
Q\subset\overline B_r(x)$. Then $\mathcal T_k\llcorner B_r(x)-Q$ is a compactly supported
integral cycle, and $\Omega^{\rm Re}_k$ is exact, so $(\mathcal T_k\llcorner
B_r(x))(\Omega^{\rm Re}_k)=Q(\Omega^{\rm Re}_k)\le\mathbf M(Q)$. Hence
\[
\mathcal V(r)\le2c_n\,\mathcal V'(r)^{n/(n-1)}
\]
for a.e.\ $r$, so that $(\mathcal V^{1/n})'\ge c(n)$ a.e., and since $\mathcal V$ is
non-decreasing, $\mathcal V(r)\ge c_0r^n$.

For $n=1$, $\Omega^{\rm Re}_k$ is a constant
$1$-form on $\C$, $\Omega^{\rm Re}_k=d\ell$ for a linear function $\ell$, and $S_r$ is an
integral $0$-current, a finite sum of points with integer multiplicities. If $S_r=0$, then
$\mathcal T_k\llcorner B_r(x)$ is a compactly supported integral $1$-cycle, and $(\mathcal
T_k\llcorner B_r(x))(d\ell)=S_r(\ell)=0$, contradicting $(\mathcal T_k\llcorner
B_r(x))(\Omega^{\rm Re}_k)\ge\frac12\mathcal V(r)>0$. Thus $S_r\ne0$. Since $S_r$ is a
boundary, $S_r(1)=(\mathcal T_k\llcorner B_r(x))(d1)=0$, so its multiplicities sum to $0$.
A nonzero integral $0$-current with this property has mass $\ge2$. Hence $\mathcal
V'(r)\ge\mathbf M(S_r)\ge2$ for a.e.\ $r$, and $\mathcal V(r)\ge2r$ since $\mathcal V$ is
non-decreasing with $\mathcal V(0^+)=0$. This proves the claim, with $c_0=2$ when $n=1$.
In both cases, by \eqref{eq:A1}, $\spt\mu_k=\spt\|\mathcal T_k\|$ and
$\mu_k(B_r(x))\ge c_0r^n$ for such $x,r$.

By the claim, every phase that carries mass in $B_R$ uses at least $c_0R^n$ of the total
mass, which is at most $D_0R^n$ on $B_{2R}$; so there are at most
$N_0:=\lfloor D_0\,2^n/c_0\rfloor$ such phases. Indeed, let $\mathfrak
m_R:=(\mathrm{pr}_2)_\#(\mathfrak m\llcorner B_R\times\R/2\pi\Z)$. Suppose $\spt\mathfrak
m_R$ contained $N_0+1$ distinct
points $t_0,\dots,t_{N_0}$. Choose $m$ so large that $2\pi/m$ is smaller than the distance
on the circle between any two of the $t_l$, and then $a\in[0,2\pi/m)$ such that (i)--(iii)
hold and no $t_l$ lies in $a+(2\pi/m)\Z$; this excludes only a null set and finitely many
values. Then every $t_l$ lies in the interior of an arc $I_{k_l}$, and the arcs $I_{k_l}$
are distinct, since each arc is shorter than the distance between two of the $t_l$. Then
$\mu_{k_l}(B_R)>0$, so there are $x_l\in\spt\mu_{k_l}\cap\overline B_R$, and by
\eqref{eq:A-density} $\mu_{k_l}(B_R(x_l))\ge c_0R^n$. The arcs are disjoint, so
$\sum_l\mu_{k_l}\le\mu$, and $\mu(B_{2R})\le\liminf\mu_i(F_i^{-1}(B_{2R}))\le D_0R^n$.
Thus
\[
(N_0+1)c_0R^n\le\sum_l\mu_{k_l}(B_R(x_l))\le\mu(B_{2R})\le D_0R^n,
\]
which contradicts $(N_0+1)c_0>2^nD_0\ge D_0$. Hence $\spt\mathfrak
m_R\subset\{\bar\vartheta_1,\dots,\bar\vartheta_N\}$ with $N\le N_0$.

Finally, choose $m$ and $a$ as in the previous paragraph, with the $\bar\vartheta_j$ in
place of the $t_l$, so that each $\bar\vartheta_j$ lies in the interior of its own arc
$I_{k_j}$, and carry out the construction of the currents $\mathcal T_k$ above for this
choice. On $B_R\times I_{k_j}$ we have $\mathfrak
m=\mu^{(j)}\otimes\delta_{\bar\vartheta_j}$, where $\mu^{(j)}=\mu_{k_j}\llcorner B_R$.
The integrals over $L^i_{k_j}$ below are the integrals of the functions
$(w,\phi)\mapsto\gamma(w)\,|\Omega_0|_{h_i}(w)\cos(\phi-\bar\vartheta_j)$ against the measures
$\mathfrak m^i\llcorner(B_{4R}\times I_{k_j})=(F_i,\vartheta_i)_\#(\mu_i\llcorner L^i_{k_j})$.
Since these functions converge uniformly to $\gamma(w)\cos(\phi-\bar\vartheta_j)$, while these
measures have bounded mass and, by (i), converge weakly to
$\mathfrak m\llcorner(B_{4R}\times I_{k_j})$,
\begin{align*}
\mathcal
T_{k_j}\big(\gamma\,\mathrm{Re}(e^{-i\bar\vartheta_j}\Omega_0)\big)&=\lim_i\int_{L^i_{k_j}}\gamma\circ F_i\,|\Omega_0|_{h_i}\cos(\vartheta_i-\bar\vartheta_j)\,d\mu_i=\mu^{(j)}(\gamma)\\
&\ge\|\mathcal T_{k_j}\|(\gamma)\ \ge\ \mathcal
T_{k_j}\big(\gamma\,\mathrm{Re}(e^{-i\bar\vartheta_j}\Omega_0)\big)
\end{align*}
for $\gamma\ge0$ supported in $B_R$. Hence $\|\mathcal T_{k_j}\|=\mu^{(j)}$ on $B_R$, and
$\mathcal T_j:=\mathcal T_{k_j}\llcorner B_R$ is calibrated by
$\mathrm{Re}(e^{-i\bar\vartheta_j}\Omega_0)$. Since $\mathfrak
m\llcorner(B_R\times\R/2\pi\Z)=\sum_j\|\mathcal T_j\|\otimes\delta_{\bar\vartheta_j}$, the
formula follows.

If $\vartheta_i$ lifts to $\tilde\vartheta_i$ with $|\tilde\vartheta_i|\le D_0$, we repeat
the argument with intervals in $[-D_0,D_0]$ instead of arcs, cutting at the points
$a+2D_0k/m$, $k\in\Z$, for almost every $a\in[0,2D_0/m)$, with the same justification as
for (i)--(iii); this uses $|\nabla\tilde\vartheta_i|=|\nabla\vartheta_i|$. We take
$m\ge\max\{3,3D_0/\pi\}$; then the phase differences from the midpoints of the intervals
are at most $D_0/m\le\pi/3$, so the argument above applies with $\cos(D_0/m)\ge\frac12$ in
place of $\cos(\pi/m)$. The lift must be a lift of the geometric phase $\vartheta_i$, not
an arbitrary function such as the Legendrian angle $\theta$ of (F).
\end{proof}

\begin{remark}\label{rem:A2}
The proof does not use Allard's theorem. Integrality of the limit comes from
Federer--Fleming applied to the pieces $L^i_k$. The lower density bound comes from near
calibration, applied to the limit currents rather than to the pieces, together with the
isoperimetric inequality for $n\ge2$ and the boundary-mass argument for $n=1$. In
particular no Lipschitz bound on the phase is needed, such as hypothesis~(b) of the
Poincar\'e-type lemma \cite[Lemma~3.7]{nev11}, which is not available for the phase. The
proposition is a statement about static sequences; when it is applied to flows
(Theorem~\ref{thm:weak}(vi)), orientation, properness and the boundary condition have to
be checked separately.
\end{remark}

\section{The curvature errors on the round three-sphere}\label{app:B}

This appendix verifies the formulas of Remark~\ref{rem:S3}. Take $n=1$, $w=(x,y)$,
$r=|w|$, $c=c_S$ and $u=1+cr^2$; a pair $(a,b)$ stands for the vector
$a\partial_x+b\partial_y$. Then
\[
h_p=\frac{dx^2+dy^2}{4u^2},\qquad d\eta=\omega_Z=\frac{dx\wedge dy}{2u^2},\qquad
\beta_p=\frac{\overline\beta}{u}=\frac{x\,dy-y\,dx}{4u},\qquad f_p=\frac{r^2}{16}.
\]
The base has constant curvature $4u^2\Delta\log u=16c$.

Since $d\overline\beta=\frac12dx\wedge dy=\omega_0$ and
$(x\,dx+y\,dy)\wedge\overline\beta=\frac{r^2}4dx\wedge dy$,
\[
d\beta_p=-\frac{2c}{u^2}\,(x\,dx+y\,dy)\wedge\overline\beta+\frac{d\overline\beta}u=\frac{u-cr^2}{2u^2}\,dx\wedge
dy=d\eta ,
\]
and $\iota_{\mathfrak r}\beta_p=0$. Moreover $h_p=g_0+O(r^2)$, so this is a contact normal
chart as in Lemma~\ref{lem:charts}. With $x_1=x$, $x_2=y$ and $\partial_ju=2cx_j$, the
Christoffel symbols of $h_p$ are
\[
\Gamma^k_{ij}=-\frac1u\big(\delta_{ik}\partial_ju+\delta_{jk}\partial_iu-\delta_{ij}\partial_ku\big).
\]

For the first claim of Remark~\ref{rem:S3}, with $\partial_if_p=x_i/8$ one finds
\[
\overline\nabla^2f_p-\tfrac12h_p=\frac
c{2u}\,(x\,dx+y\,dy)^{\otimes2}-\frac{c^2r^4}{8u^2}\,(dx^2+dy^2)
\]
on horizontal lifts. Relative to $h_p$ its eigenvalues are $2cr^2+\frac32c^2r^4$ in the
radial direction and $-\frac12c^2r^4$ in the angular direction, so its $h_p$-norm is
$2cr^2+O(r^4)$.

For the second claim, $|\overline\nabla
f_p|^2=|\nabla^{h_p}f_p|^2_{h_p}=4u^2\cdot\frac{r^2}{64}$, so $|\overline\nabla
f_p|^2-f_p=\frac{r^2}{16}(u^2-1)=\frac{cr^4}{16}(2+cr^2)$.

For the third claim we use that, by the proof of Lemma~\ref{lem:model}(c),
$\overline\nabla^2z_p(U,W)=-(\nabla^Z_U\beta_p)(W)-(\nabla^Z_W\beta_p)(U)$ on horizontal
lifts. All data are invariant under rotations of $w$, so it suffices to compute at
$w=(r,0)$. There $\beta_p=\frac r{4u}\,dy$, the derivatives of the components are
$\partial_y(\beta_p)_x=-\frac1{4u}$, $\partial_x(\beta_p)_y=\frac1{4u}-\frac{cr^2}{2u^2}$,
$\partial_x(\beta_p)_x=\partial_y(\beta_p)_y=0$, and the only nonzero Christoffel symbols
are $\Gamma^x_{xx}=-\Gamma^x_{yy}=\Gamma^y_{xy}=\Gamma^y_{yx}=-2cr/u$. Hence, at $(r,0)$,
\[
\overline\nabla^2z_p=-\frac{cr^2}{2u^2}\,(dx\otimes dy+dy\otimes dx),
\]
and the eigenvalues of $h_p^{-1}\overline\nabla^2z_p$ are
$\pm4u^2\cdot\frac{cr^2}{2u^2}=\pm2cr^2$.

Finally, for the fourth claim, $\overline\nabla^{\mathcal H}z_p$ is the horizontal lift of
$-2\beta_p^\sharp$ (proof of Lemma~\ref{lem:model}(d)), and $\Phi\overline\nabla f_p$ is
the horizontal lift of $J_0\nabla^{h_p}f_p$, where $J_0(a,b)=(b,-a)$ by
Section~\ref{subsec:sasakian}. Since $\beta_p^\sharp=u\,(-y,x)$ and
$\nabla^{h_p}f_p=\frac{u^2}2(x,y)$,
\[
-2\beta_p^\sharp-4J_0\nabla^{h_p}f_p=2u(1-u)\,(y,-x)=-2cr^2u\,(y,-x),
\]
whose $h_p$-norm is $\frac1{2u}\cdot2cr^2u\cdot r=cr^3$. For $c=0$ this is the model
identity $-2\overline\beta^\sharp=4J_0\nabla^{g_0}f$.

\section*{Use of artificial intelligence}
ChatGPT (OpenAI) and Claude (Anthropic) assisted in developing and checking mathematical
arguments and in exploratory calculations. The authors take full responsibility for the paper.

\bibliographystyle{amsalpha}
\bibliography{HeightLocalised-arxiv-v1}

\end{document}